\documentclass[11pt]{amsart}

\usepackage{macros}
\usepackage{epigraph}

\usepackage{setspace}
\spacing{1.1}

\usepackage{parskip}
\author{Owen Gwilliam and Brian R. Williams}
\date{\today}
\title{Holomorphic field theories and higher algebra}

\renewcommand{\op}{\operatorname}
\def\lie#1{\ensuremath{\mathfrak{#1}}}

\begin{document}

\dedicatory{In memory of Yuri Manin, with gratitude and admiration.}

\maketitle

\epigraph{The mathematical language of classical physics is based upon real numbers \dots. The mathematical language of quantum physics is based upon complex numbers, and it would be natural to expect that the complex analytic and the algebraic geometry should replace the differential geometry of the classical period.}{Yuri Manin from the introduction of {\it Frobenius Manifolds,\\ Quantum Cohomology, and Moduli Spaces}}

\begin{abstract}
Aimed at complex geometers and representation theorists, this survey explores higher dimensional analogues of the rich interplay between Riemann surfaces, Virasoro and Kac-Moody Lie algebras, and conformal blocks.
We introduce a panoply of examples from physics --- field theories that are holomorphic in nature, such as holomorphic Chern-Simons theory --- and interpret them as (derived) moduli spaces in complex geometry;
no comfort with physics is presumed.
We then describe frameworks for quantizing such moduli spaces, 
offering a systematic generalization of vertex algebras and conformal blocks via factorization algebras,
and we explain how holomorphic field theories generate examples of these higher algebraic structures.
We finish by describing how the conjecture of Seiberg duality predicts a surprising relationship between holomorphic gauge theories on algebraic surfaces and how it suggests analogues of the Hori-Tong dualities already studied by algebraic geometers.
\end{abstract}

\tableofcontents

\section{What is a holomorphic field theory? Why study them?}

The most important field theories in physics involve a metric on a manifold (whether Riemannian or Lorentzian), and hence geometry has long had a clear role in physics.
In recent decades attention has expanded to include {\em topological} field theories (or TFTs), 
in which --- to simplify --- only the underlying smooth topology of the manifold matters.
There are, however, theories that depend on a complex structure on the manifold,
and as one might hope, these are particularly beautiful, 
just as complex variables is a particularly beautiful wing of analysis.
In essence, such {\em holomorphic} field theories have variational partial differential equations (i.e., ``equations of motion'') that are holomorphic.
These theories admit independent motivations from mathematics and physics,
which we sketch before giving more careful definitions.

The mathematical motivation is simple:
for most classical holomorphic field theories, the solutions to the equations of motion form a moduli space of natural interest to complex geometers,
and one might hope that the quantization of this moduli space is equally interesting.
A well-known example is holomorphic Chern-Simons theory,
which lives on a Calabi-Yau 3-fold~$X$ and involves the moduli space $\RR\Bun_G(X)$ of holomorphic principal $G$-bundles on~$X$.
(We add the prefix $\RR$ because it is best to use the {\em derived} moduli space.
We do not assume the reader is versed in derived geometry and will offer motivations, glosses, and references throughout the survey.
See Section~\ref{sec: DAG} for a start.)

As a more mundane example, let $\RR\Hol(X,Y)$ denote the derived moduli space of holomorphic maps from some complex manifold $X$ to another complex manifold $Y$.
For $Y = V$ a complex vector space, this ``derived'' enhancement is simply the entire  cohomology $H^*(X, \cO) \otimes V$, and not just the global functions $H^0(X,\cO) \otimes V$.
There are holomorphic theories whose derived moduli spaces are encoded by cohomology groups of holomorphic vector bundles.
(A physicist might view them as analogs of scalar field theories.)

A physical motivation comes from supersymmetric field theories, 
which (loosely speaking) are theories on $\RR^n$ whose symmetries include not only the isometries of $\RR^n$, but an extension to a super Lie group known as a super Poincar\'e group.
Thanks to the enhanced symmetry, such theories can often be understood in greater detail than non-supersymmetric theories.
The relevance to this survey is that supersymmetric theories admit ``twists'' (in essence, deformations) that are holomorphic field theories, 
even when they do not admit twists to TFTs \cite{CosSUSY}. 
(The famed A- and B-models of mirror symmetry are TFTs produced by twisting supersymmetric thories.)
As an example, four-dimensional minimally supersymmetric Yang-Mills theory admits no topological twist, 
but it does admit a twist to holomorphic field theory whose moduli space is described, in part, by the moduli space of holomorphic $G$-bundles~$\RR \Bun_G(X)$. 

It is natural to wonder how the myriad tools and computations deployed on a supersymmetric theory carry over to this holomorphic twist 
and, conversely, whether information can flow back from a holomorphically twisted theory to the original supersymmetric theory, 
much as TFTs have enriched physicists' understanding of supersymmetric theories.
In short, a holomorphic twist contains the ``holomorphic sector'' of the untwisted theory just as a topological twist contains its ``topological sector.''
To aid translation between mathematics and physics, 
it may be helpful to know that the role of derived geometry on the mathematical side is mirrored by the use of the BRST/BV formalism on the physics side;
they are closely related and serve similar purposes.

As a final motivational remark, we note that in one complex dimension (i.e., working on Riemann surfaces), 
holomorphic field theories appear as chiral conformal field theories (CFTs) on Riemann surfaces. 
The study of CFTs --- whose manifestations in mathematics include vertex algebras, affine Lie algebras, nonabelian theta functions, and loop groups --- have had a strong impact in mathematics,
so it is natural to explore higher-dimensional analogs.
So far, most mathematical effort in this direction has focused on (essentially) classical aspects of gauge theory,
but the quantum aspects are profoundly interesting.
We expect that these holomorphic examples will be particularly tractable and beautiful but also useful in offering insights even for non-holomorphic theories.

\subsection{The aims of this survey}

Over the last few decades there has been a vigorous dialogue between complex geometry and theoretical physics,
around topics like mirror symmetry, enumerative algebraic geometry, moduli of Calabi-Yau manifolds, the geometric Langlands program, and integrable systems, among others.
Covering all these topics would be difficult.
Our survey's focus will be narrow:
it is upon higher dimensional versions of the relationship between complex geometry and algebraic structures like Lie algebras, Poisson algebras, and vertex algebras.
We will emphasize analogies to TFTs (such as the Gerstenhaber algebras appearing in mirror symmetry) and CFTs (such as the affine Lie algebras relating to moduli of $G$-bundles on Riemann surfaces),
but our primary interest is on complex manifolds above dimension one.
Experience suggests that higher dimensions require higher algebra,
i.e., algebraic methods that involve cochain complexes.
Thus differential graded (dg) algebras and categories,
as well as fellow travelers like derived geometry, appear.

A community has formed, involving both mathematicians and physicists, to explore and exploit these relationships,
and it continues to grow as there is a lot to explore.
Examples of efforts in this direction include \cite{PTVV}, \cite{FHK,HenKap}, \cite{CLbcov, LiVertex}, \cite{CGholography,BGKWY}, \cite{CostelloPaquette}, \cite{BBBDN}, and \cite{AlfonsiYoungHC},
where our own background and interests (and unfortunately, our ignorance) have shaped this list.

The material below focuses on connections with this work,
and the examples from physics have a perturbative flavor.
We will not discuss important directions like mirror symmetry, the geometric Langlands program, the BPS/CFT correspondence, Donaldson-Thomas theory, or other areas that are more global or more directly geometric,
although we expect there are rich connections with the topics we do discuss.
Not including these topics is primarily because our core theme already swelled to more pages than we expected, 
because a concise discussion of any of those topics does not yet seem possible with respect to our theme,
and because of our lack of background in these areas.
We hope the ideas surveyed here do connect profitably to these topics.

As complementary perspectives to that offered here, we recommend
\begin{itemize}
\item Kapranov's lectures \cite{KapLect} on factorization algebras in algebraic geometry,
\item Costello's lectures \cite{CosScheim} on how holomorphic field theories emerge from supersymmetric theories, and
\item a recent survey \cite{GKW} for physicists on how higher algebra appears in perturbative QFT, with a strong emphasis on the holomorphic examples.
\end{itemize}
Reading these would give a more well-rounded view of the subject and its possibilities.

\subsection{How to read this paper}

Section~\ref{sec: HFT} is devoted to describing holomorphic field theories, primarily by offering a bounty of examples, including holomorphic analogs of gauge and gravity theories.
We outline as well how to obtain holomorphic field theories from supersymmetric field theories,
a rich source of examples to explore.

The next three sections are, in a sense, the core of the survey,
as they describe increasingly general types of higher algebras and how these capture the observables of holomorphic field theories.
Section~\ref{sec: algebra} revisits the relationship between quantum mechanics and associative algebras,
chiral conformal field theory and vertex algebras,
and holomorphic field theories and higher-dimensional generalizations of vertex algebras.
Section~\ref{sec: config} offers a useful framework for studying these higher vertex algebras,
using operads,
and it develops powerful analogies to the little $n$-disks operad $E_n$ that plays an important role in contemporary topology.
Section~\ref{sec: fact alg} introduces factorization algebras,
providing a local-to-global principle to describe the observables on complex manifolds beyond affine spaces.
In particular, by pushing forward factorization algebras,
one can relate field theories of different dimension,
i.e., a version of the physicist's ``compactification.''

A reader could read this core portion independently of the rest of the paper,
as it focuses on algebra and geometry and requires no substantial knowledge of field theory {\it per se}.

Section~\ref{sec: renorm} describes recent work that allows one to produce examples:
there is now a convenient framework for renormalization of holomorphic field theories on $\CC^n$,
and so the field theories of Section~\ref{sec: HFT} generate higher algebras of the form described in Sections~\ref{sec: algebra},~\ref{sec: config}, and~\ref{sec: fact alg}.
The reader could skip to this section directly from Section~\ref{sec: HFT}.

We conclude the paper by describing the kind of fascinating mathematical questions
that arise by taking physical conjectures (e.g., dualities of supersymmetric theories) and twisting them into conjectures about holomorphic field theories. 
Section~\ref{sec: seiberg} explores the holomorphic cousin of Seiberg duality,
which suggests a remarkable correspondence between a gauge theory with group $SL(N)$ and a gauge theory with group $SL(M)$, where $M \neq N$.
It is a kind of enhancement of Weyl's invariant theory and has close connections with topics like homological projective duality.

\subsection{About Manin and his influence}

Yuri Manin was a pioneer into this intellectual region where field theory, complex geometry, and homotopical algebra meet,
and his work --- beyond the beautiful results --- also gestured at rich possibilities for novel and powerful forms of geometry and physics and algebra.
His writing summoned these possibilities out the realm of dreams and towards their current unfolding today into derived geometry and higher operadic structures.
To whit, our epigraph comes from \cite{ManFrob}, which synthesizes stacks and operads and Frobenius structures into a beautiful and coherent view on mirror symmetry.
But we could equally well have quoted his \cite{ManGauge}, which is rich with deep reformulations of supersymmetric field theory that connect it to complex geometry,
or his \cite{BeiMan, ManNewDim} from roughly the same era.
Anyone interested in mathematics and physics, and especially their relationship, will delight in~\cite{ManMP}.

Beyond his writing, however, we have also benefited throughout our careers from the attitudes and atmosphere that his persona modeled and spread among mathematicians:
a generosity with mathematical ideas but also the sharing of poetry and literature, 
an exquisite care and sensitivity to organizing and explaining mathematics (and much else),
and a kindness in conversation, where he brought his whole attention to you.
While a graduate student at Northwestern, 
OG had the unbelievable luck to overlap with Manin and thus to take several courses and have many conversations with him, 
getting to witness his personification of the Platonic idea of mathematician as intellectual in the broadest sense.
We also both spent extensive time, for example, at the Max Planck Institute for Mathematics,
which is suffused with his spirit.
We wish we had the chance to discuss the ideas from this survey with Manin;
we offer it as a small gift in gratitude for all he shared with us.

\subsection{Acknowledgements}

This survey attempts to summarize the efforts and insights of a large body of researchers.
So much of our understanding is shaped by conversations and interactions with them,
and we thank them for their generosity with ideas and examples.
We owe particular thanks to Kasia Budzik, Dylan Butson, Kevin Costello, Richard Eager, Chris Elliott, Davide Gaiotto, Nik Garner, Vasily Gorbounov, Zhengping Gui, Benjamin Hennion, Mikhail Kapranov, Justin Kulp, Si Li,  Natalie Paquette, Surya Raghavendran, Nick Rozenblyum,  Ingmar Saberi,  Pavel Safronov,  Matt Szczesny, Ed Segal,  Johannes Walcher, Minghao Wang, Jingxiang Wu, Philsang Yoo, and Jie Zhou.

The prompt to write this survey came from Andrey Lazarev, and we thank him for his patience and encouragement during the long period of putting it together.

The National Science Foundation supported O.G. through DMS Grant 2042052 during the writing of this paper.
The Max Planck Institute for Mathematics has hosted both O.G. and B.W. throughout their careers, and it provided a wonderful environment for O.G. during the last stages of writing.

\section{Holomorphic field theories}
\label{sec: HFT}

We begin by explaining the basic framework of classical field theory, with a few examples.
For a mathematician, the key point might be that physics suggests certain systems of equations are particularly important and studies their spaces of solutions, 
often unearthing aspects that might not catch a mathematician's eye.
In other words, physics can be a source of interesting moduli spaces.

We then offer a quick overview of derived geometry, aiming to motivate and not to give a systematic treatment.
For this survey, comfort with basic homological algebra and sheaf cohomology will suffice, 
so we will explain how our terminology relates to those notions.
In Section~\ref{sec: pert dfn} we offer a sharp definition well-suited to perturbative field theory 
(i.e., studying formal deformations of a fixed solution to the equations of motion),
which we will use to produce examples of higher algebras.

\subsection{A first pass}

In physics, a classical field theory is usually described in terms of a few pieces of data:
\begin{itemize}
\item a manifold $M$ that plays the role of ``spacetime,''
\item a space of ``fields,'' typically the sections of some fiber bundle $F \to M$ or connections on a principal bundle on $M$, and
\item a system of partial differential equations (PDE) for the fields, known as the ``equations of motion.''
\end{itemize}
These equations of motion are variational: 
the solutions to these PDE are the critical points of some function on the space of fields,
called an ``action functional,''
and such PDE are often called Euler-Lagrange equations.
(One must take a little care in interpreting the claim of taking critical points.)
In practice an action functional is a formal expression, 
written as an integral over $M$ of some polydifferential operator on the fields that takes values in densities of~$M$.

A few examples will illuminate this sketch.
We begin by giving two examples that depend on Riemannian geometry, 
then give an example of a topological field theory,
and finally give an example of a holomorphic field theory.
All the equations that show up are well-known to geometers.

\begin{eg}[Massless scalar field theory]
\label{eg: massless scalar}
Let $M$ be a smooth manifold and equip it with a Riemannian metric~$g$.
Let $\Delta_g$ denote the associated Laplace-Beltrami operator, 
and let $\d{\rm vol}_g$ denote the associated volume form.
In the case $M = \RR^n$ with the Euclidean inner product, the Laplace-Beltrami operator is the differential operator
\[
\Delta = \frac{\partial^2}{\partial x_1^2} + \cdots + \frac{\partial^2}{\partial x_n^2}.
\]
Let the space of fields be the smooth functions $C^\infty(M)$.
A smooth function $\phi$ is called, in physics, a {\em scalar} field.
The action functional
\[
S(\phi) = \int_M \phi \Delta_g \phi \, \d{\rm vol}_g
\]
yields the equation of motion $\Delta_g \phi = 0$. 
Thus this scalar field theory has the moduli space of {\em harmonic functions} as its space of solutions.
\end{eg}

The preceding example exhibits how a moduli space of central interest in differential geometry --- the harmonic functions --- appears naturally as a basic example in physics (the ``free massless scalar field'').
More complicated moduli spaces show up too.

\begin{eg}[$\phi^4$ theory]
As a small variant of the example just given,
consider the action functional
\[
S(\phi) = \int_M \phi \Delta_g \phi + \frac{1}{4!} \phi^4\, \d{\rm vol}_g.
\]
Its equation of motion is the nonlinear PDE 
\[
\Delta_g \phi + \frac{1}{3!} \phi^3 = 0.
\] 
Thus this scalar field theory has a more complicated, nonlinear moduli space of solutions.
\end{eg}

We now describe an example of a theory whose space of solutions is a moduli space of interest in geometric topology: character varieties.

\begin{eg}[Chern-Simons theory]
Pick a compact Lie group $G$ and equip it with an invariant metric,
which induces a $G$-invariant inner product $\kappa$ on its Lie algebra~$\fg$.
(The standard choice is the Killing form.)
Take an oriented smooth 3-manifold $M$, 
and consider the space of connections on the trivial principal $G$-bundle on $M$.
This space of fields is isomorphic to $\Omega^1(M, \fg)$, the $\fg$-valued one-forms, 
by viewing $\d + A$ as a connection where $A \in \Omega^1(M, \fg)$.
The equation of motion is the zero-curvature equation
\[
F_A = 0,
\]
which can be expressed equivalently as the Maurer-Cartan equation
\[
\d A + \frac{1}{2} [A,A] = 0.
\]
These equations arise as the Euler-Lagrange equation for the action functional
\[
S(A) = \int_M CS(A)
\]
where $CS$ denotes the Chern-Simons form
\beqn
\label{eqn:CS}
CS(A) = \frac12 \kappa(A, \d A) + \frac16 \kappa(A, [A,A]).
\eeqn
Hence this theory is known as {\em Chern-Simons theory}.\footnote{Note that $M$ needs to be an oriented 3-manifold for this action functional to make sense.
The equation of motion exists for a principal bundle on any manifold of any dimension,
but it is only variational for oriented 3-manifolds.}
These equations do not depend on a metric on $M$ or any other geometric structure,
and hence define a ``topological'' field theory in the sense that they are independent of geometry.
This theory is interesting because the moduli space of solutions consists of the flat $G$-bundles on $M$,
a topic of deep interest in geometry and topology.
For physicists, it is interesting because it does not involve geometry, 
by contrast to earlier examples of field theories,
and yet exhibits rich phenomena (notably applications to the quantum Hall effect and 3-dimensional gravity).
\end{eg}

\begin{rmk}
In the Chern-Simons example, there is an important aspect we did not yet acknowledge:
we need to identify different solutions related by gauge transformations.
In other words, the correct moduli space is the quotient stack
\[
\left[ \{ A \,:\, F_A = 0\}/\Map(M,G) \right].
\]
The usual mode of description of a field theory, by specifying fields and an action functional (or equations of motion), 
must be supplemented by specifying whether to quotient by some group of symmetries of the space of solutions.
More generally, a field theory ought to present some kind of derived moduli stack.
\end{rmk}

We are now in a good position to offer an example of a {\em holomorphic field theory}:
the holomorphic analog of Chern-Simons theory.

\begin{eg}[Holomorphic Chern-Simons theory]
\label{eg: holCS}
Pick a {\em complex} Lie group $G$, 
and equip its Lie algebra $\fg$ with a complex-linear, $G$-invariant nondegenerate symmetric bilinear form~$\kappa$.
Let $M$ be a three-dimensional \textit{complex} manifold
equipped with a holomorphic volume form $\Omega \in \Omega^{3,0}(M)$; 
in other words, $M$ is a Calabi--Yau manifold.
(This condition is analogous to being an \textit{oriented} 3-manifold.)
Given a smooth principal $G$-bundle $P \to M$,
there is a moduli space of holomorphic structures on~$P$.
Suppose there exists at least one holomorphic structure $\cP \to M$,
which has an associated Dolbeault operator (or $\dbar$-connection) $\dbar_{\cP}$.
Any other holomorphic bundle structure determines a Dolbeault operator,
but our reference point $\dbar_{\cP}$ provides an isomorphism
\[
{\rm Conn}^\dbar_G(P) \cong \Omega^{0,1}(M, \fg_P),
\]
i.e., we express another $\dbar$-connection as $\dbar_{\cP} + A$ where $A$ is a $(0,1)$-form valued in the adjoint bundle.
The equation of motion is again a Maurer--Cartan equation
\beqn
\dbar_{\cP} A + \frac12 [A,A] = 0 ,
\eeqn
or equivalently the $(0,2)$-term of the curvature vanishes.
This PDE is the Euler--Lagrange equation for the action
\beqn
S_{hCS} (A) = \int_M \Omega \wedge CS(A) 
\eeqn
where $CS(A)$ is exactly as in~\eqref{eqn:CS}.\footnote{Notice that due to the presence of the holomorphic volume form $\Omega$, only the $\dbar A$ component of $\d A$ contributes to the integrand of the action functional, because $\Omega \wedge \partial A$ vanishes.} 
As in the case of Chern-Simons theory, one must quotient by the appropriate group of gauge transformations to obtain the correct moduli space, 
which is here the moduli stack of holomorphic principal $G$-bundles whose underlying smooth structure is given by $P \to M$.
\end{eg}

For a systematic introduction to classical field theory, 
in a style amenable to a complex geometer,
we suggest \cite{DeligneFreed}.
In the appendix (see Section~\ref{sec: ex of theories}) we describe many more examples of holomorphic field theories and describe their spaces of solutions using the language of derived geometry, which we discuss next.

\begin{rmk}
Recently there is a substantial effort to revisit field theory with insights from derived geometry.
See \cite{AlfonsiYoung} for motivation and overview, 
but for a detailed development, see \cite{Steffens23, Steffens24}.
\end{rmk}

\subsection{A dose of derived geometry}
\label{sec: DAG}

Derived geometry is a sophisticated subject but the motivations and first steps are surprisingly simple,
as we now attempt to explain.
In this paper we only need the reader to understand the rough outline of derived geometry to follow our discussion,
but there are several, much more substantial expositions now available.
(For algebraic geometers consider \cite{Toen,PridEug}; for symplectic or Poisson geometers consider \cite{PTVV, Safronov}; for differential geometers consider \cite{Joyce}; and for representation theorists consider \cite{GaiRoz}.)

Before saying anything technical, we indicate how derived geometry generalizes more familiar bits of mathematics.
The starting place is the theory of vector spaces, which is very concrete and provides a rich source of both geometric intuitions and algebraic tools.
This linear setting generalizes in two directions:
\begin{itemize}
\item the geometry of nonlinear spaces, whether in the style of algebraic geometry or in the style of differential geometry, and
\item homological algebra, or differential graded vector spaces (and modules, sheaves, and so on).
\end{itemize} 
Both of these generalizations pass back to linear algebra in natural ways:
\begin{itemize}
\item at a point $x$ on a manifold $M$ (or variety), the tangent space $T_xM$ offers the best linear approximation to the space at that point (where it makes sense), and
\item given a cochain complex $(V,d)$, a cohomology group $H^0(V)$ is a vector space encoding essential information about the complex.
\end{itemize} 
As a first pass at derived geometry, one might say that it seeks a common generalization of nonlinear geometry and homological algebra.
With slightly more precision, one might ask for a notion of derived space $M$ such that
\begin{itemize} 
\item each point $x$ in $M$ has a tangent {\em complex} $\TT_x M$, i.e., the linearization of a point yields a cochain complex, and
\item $M$ should have a truncation $\tau(M)$ to an ordinary space.
\end{itemize}
Note that there will be different versions of derived geometry, depending on the kind of nonlinear geometry concerned: 
there is derived algebraic geometry, derived differential geometry, derived symplectic geometry, and so on.

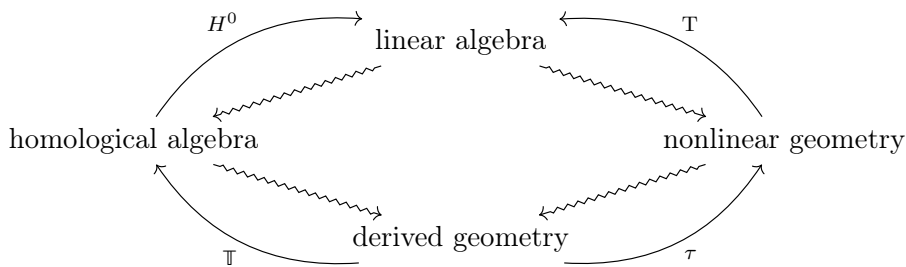
\begin{center}
\begin{tikzcd}
&\text{linear algebra} \arrow[squiggly, dr] \arrow[squiggly, dl] &\\
\text{homological algebra} \arrow[squiggly, dr] \arrow[ur, bend left, "H^0"] &&\text{nonlinear geometry}\arrow[squiggly, dl] \arrow[ul, bend right, "\rm{T}"']\\
&\text{derived geometry} \arrow[ul, bend left, "\TT"] \arrow[ur, bend right, "\tau"']&
\end{tikzcd}
\end{center}

A related desideratum of derived geometry is that 
\begin{itemize}
\item an ordinary scheme or manifold (or whatever class of underived space) should provide an example of a derived space, and
\item any cochain complex should provide an example of a derived space.
\end{itemize}
That is, derived geometry contains both the domains it seeks to generalize.

There is an interesting subtlety here: quasi-isomorphic cochain complexes should present the same derived space, just as they present the ``same'' object in homological algebra.
This requirement means that the true home for homological algebra is the derived category,
obtained by localizing cochain complexes at quasi-isomorphisms.
Similarly, the true home of derived geometry is also typically presented somewhat indirectly,
with different descriptions provided of the same derived space.
To do this carefully requires the theory of $\infty$-categories (or related tools like model categories).
We will not dwell on this aspect of the subject,
but we hope the comparison with the derived category gives a sense of the issues involved.
(It also explains the terminology: it is the same use of ``derived'' for category and for geometry.)

To present a derived space, two approaches are used and they mimic familiar methods for manifolds or schemes:
\begin{itemize}
\item a topological space equipped with a structure sheaf, here a sheaf of dg commutative algebras, or
\item a functor of points, i.e., as a moduli space.
\end{itemize}
The devil is in the details of pinning down, for example, what kind of sheaves appear as structure sheaves or what the source and target categories appear for the functor of points,
but let us postpone dealing with the devil and discuss a few examples used in this paper.
We will use complex algebraic geometry as our model for nonlinear spaces, for simplicity,
and we will emphasize the first style of describing a space.

\subsubsection{Derived baby steps}

Recall that the $n$-dimensional vector space $V = \CC^n$ determines a complex variety $\AA^n$ called an affine space.
All the essential information is encoded in its algebra of functions
\[
\Gamma(\AA^n, \cO) = \CC[x_1,\ldots, x_n] = \Sym(V^*)
\]
a polynomial ring in $n$ generators, where the generators form a basis for the dual vector space $V^*$ to~$V$.
We can view $\AA^n$ as a ringed space whose underlying topological space is $\CC^n$ with the Zariski topology
and whose structure sheaf is obtained by localizing the polynomial ring accordingly.
In derived geometry, a $\ZZ$-graded vector space 
\[
V = ( \ldots, V^{-1}, V^0, V^1, \ldots)
\]
should likewise present a derived space.
(We assume $V$ has finite total dimension, for simplicity.) 
The essential information is encoded in its graded algebra of functions
\[
\Sym(V^*) = \bigoplus_{n \geq 0} \Sym^n(V^*)
\]
where $V^*$ is the dual graded vector space with $(V^*)^k = (V^{-k})^*$ and where the symmetric powers involve the Koszul sign rule.
That is, if $k$ is odd, then the $n$th symmetric power $\Sym^n(V^k)$ is concentrated in degree~$nk$ 
but it is isomorphic to the $n$th exterior power of the {\em un}\/graded vector space~$V^k$.
In other words, $\Sym(V^*)$ is a graded-commutative version of a polynomial ring.
One can view this derived ``affine space'' as a ringed space whose underlying topological space is $V^0$ with the Zariski topology 
and whose structure sheaf is obtained by localizing $\Sym(V^*)$ accordingly.

Let's complicate this example a little.

Let $X$ be a smooth complex variety, and let $V \to X$ be an algebraic vector bundle.
It can be viewed as a family of complex vector spaces parametrized by~$X$.
Similarly, we could pick an integer $k$ and consider the graded vector bundle $V[k] \to X$ whose fiber $V[k]_x$ at a point $x \in X$ is the shift by $k$ by the fiber $V_x$.
(For us, $V[k]$ is placed in degree $-k$ so that $V[k]^m = V^{k+m}$. 
In particular, $V[1]$ is in degree~$-1$.)
But we can also talk about the total space of a vector bundle $V \to X$ as a complex variety in its own right.
In derived geometry, we can view this graded vector bundle $V[k] \to X$ as a derived space.
Its graded algebra of functions is (by definition)
\[
\Gamma(X,\Sym_{\cO_X}(\cV^*[-k])) = \bigoplus_{n \geq 0} \Gamma(X, \Sym^n_{\cO_X}(\cV^*[-k]))
\] 
where $V^*$ is the dual vector bundle, $\cV^*$ denotes the sheaf of sections of this vector bundle, and the symmetric powers again involve the Koszul sign rule.
Notice that when $k = 0$, we recover the usual ring of functions for the total space of~$V \to X$.

This class of examples appears often throughout this paper.
In particular, we often work with $T[k]X$, the total space of the shifted tangent bundle $TX[k] \to X$ of a complex variety~$X$,
or with $T^*[k]X$, the total space of the shifted cotangent bundle of a complex variety~$X$.

These shifted cotangent bundles admit natural {\em shifted} symplectic structures,
just as the unshifted cotangent bundle $T^*X$ has a natural symplectic structure.
The symplectic form on $T^*[k] X$ is given by exactly the same formula, locally
\[
\omega = \sum_j \d p_j \wedge \d x_j,
\]
where the $x_j$ are coordinates on $X$ and the $p_j = \partial/\partial x_j$ are the associated coordinates for cotangent direction.
The difference with the standard case is that now the $p_j$ carry cohomological degree $k$, so $\omega$ is a 2-form (in the de Rham complex of $T^*[k]X$) of cohomological degree~$k$.
(Note that in a derived setting cohomological degree need not coincide with form degree.)
One says that $T^* [k] X$ is a $k$-shifted symplectic space.
Implicit in this discussion is the fact that one can construct a de Rham complex with dg commutative algebra\footnote{The 1-forms are built by taking the derived functor of K\"ahler differentials, and the rest follows the usual process (e.g., take exterior powers over the dg commutative algebra, using the derived tensor product). The de Rham complex will be bigraded, with a form grading and a cohomological grading.}
and that many maneuvers and definitions carry over {\it mutatis mutandis}.

So far we have only worked with graded-commutative algebras,
as a stepping stone in the derived direction.
Let's now see something fully derived.

\subsubsection{Derived intersections and stacky quotients}

Algebraic geometry has already experienced a vast extension of the notion of space via {\it stacks}.
These appeared, in essence, to provide a more pliable and informative approach to taking a quotient of a scheme (or algebraic space) by a group action,
and this framework is the appropriate place for studying important moduli spaces, like moduli of algebraic curves or vector bundles.
From a rather abstract point of view, we replace the category of schemes by a bigger category where the colimits (a quotient is an example of a colimit) remember more information.
Derived geometry arises by asking for a (yet) bigger category where the {\it limits} remember more information.
The key geometric example of a limit is the intersection of two subvarieties in an ambient variety (i.e., the fiber product);
derived geometry changes what the intersection is.

The essential idea is very simple.
In algebraic geometry, if we have two schemes mapping to another scheme
\[
\Spec(A) \to \Spec(C) \leftarrow \Spec(B),
\]
then the fiber product, or intersection, is
\[
\Spec(A) \times_{\Spec(C)} \Spec(B) = \Spec(A \otimes_C B).
\]
Derived algebraic geometry says that the derived fiber product should be encoded by the derived tensor product~$A \otimes_C^{\LL} B$ instead.
Thus the underlying underived intersection is encoded in 
\[
H^0(A \otimes_C^{\LL} B) = A \otimes_C B
\]
but there are sometimes ``derived corrections'' in negative cohomological degrees.

For example, given a {\em non}\/regular sequence of irreducible polynomials $\mathbf{f} = (f_1,\ldots,f_k) \in \CC[x_1,\ldots,x_n]$, 
then the underived intersection of the hypersurfaces $V_j = \{f_j = 0\}$ is some quite singular subvariety $X$ of $\AA^n$.
A resolution of its algebra $\CC[X]$ over the polynomials is very complicated.
By contrast the derived intersection is encoded by the Koszul complex determined by the sequence $\mathbf{f}$,
which is a very simple cochain complex (in fact, a semi-free dg commutative algebra). 
It remembers information about the hypersurfaces (and hence $\mathbf{f}$), whereas the underived intersection only knows about the ideal generated by~$\mathbf{f}$.
In this sense it is a ``smarter'' or less ``forgetful'' limit.

Notice that such derived tensor products will naturally produce dg commutative algebras concentrated in nonpositive degrees,
as they have the nature of resolutions.
We thus generalize commutative algebras $\CAlg$ by nonpositively graded dg commutative algebras~$\CAlg^{\dg}_{\leq 0}$, rather than the bigger category of all dg commutative algebras~$\CAlg^{\dg}$.
The category of derived affine schemes is defined, then, as the category of nonpositively graded dg commutative algebras, up to quasi-isomorphism.\footnote{In fact, one works with the $\infty$-category, which is not just the localization at quasi-isomorphisms, but remembers how quasi-isomorphisms relate up to homotopy, etc.}
Each derived affine scheme $\dSpec(A)$ has a truncation to an affine scheme, which is $\Spec(H^0(A))$.
A derived scheme can be understood as ringed space patched together from derived affine schemes.\footnote{Up to some subtleties, one can think of derived manifolds in terms of manifolds equipped with structure sheaves that are nonpositively graded dg commutative algebras built from $C^\infty$-modules.}

Derived stacks then arise from derived schemes by intelligently taking quotients and, more generally, colimits.
Recall that a scheme in algebraic geometry has a functor of points valued in $\Sets$ while a stack has a functor of points valued in $\Gpds$, 
the category of groupoids.\footnote{Every equivalence relation determines a groupoid where the relation gives the morphisms. This groupoid knows a lot more than the underlying quotient set, which is known as the ``connected components'' of this groupoid. In this way, groupoids offer a refinement of traditional quotients.}
More accurately we work with groupoids localized at Morita equivalence ({\it cf.} quasi-isomorphism of cochain complexes and the derived category).
In derived geometry we push farther and work with topological spaces up to weak homotopy equivalences (or simplicial sets up to weak homotopy equivalence),
which remembers more;
every topological space has an underlying fundamental groupoid.
Thus the nature of derived geometry is indicated by the (co)domains of the functors of points:
\[
\begin{tikzcd}[column sep=large]
\CAlg^{\dg}_{\leq 0} \arrow[dd, bend right=40, "H^0"'] \arrow[rr, "\text{derived stack}"] & & \Spaces \arrow[d, bend left=40, "\Pi_{\leq 1}"]\\
&& \Gpds \arrow[d, bend left=40,  "\pi_0"] \arrow[u, hook]\\
\CAlg \arrow[uu, hook] \arrow[rr, "\text{scheme}" near end] \arrow[urr, "\text{stack}"] && \Sets \arrow[u, hook]
\end{tikzcd}
\]
(This diagram is purloined from Vezzosi's wonderful (and short!) ``What is \dots a derived stack?''~\cite{VezAMS}.)
A derived stack admits a truncation to an ordinary stack by restricting the domain category and mapping down to groupoids as the codomain category.

To ``derive'' a notion from traditional geometry, one typically needs to recast it using the functor of points and then ponder how its features might generalize to derived affine schemes and~$\Spaces$.
By now many such notions have received such a treatment, including analytic geometry~\cite{Porta} and Poisson geometry~\cite{CPTVV}.

\subsubsection{On deriving classical field theory}

The basic ingredients of classical field theory are a ``spacetime'' manifold $M$, the ``fields'' given by sections of some fiber bundle $E \to M$, and a system of PDE on the fields known as the ``equations of motion.''
The essential goal is to describe the space of solutions to these equations of motion.
All these features admit derived versions: one can replace the manifold $M$ by a derived space $\MM$ of some kind and the fiber bundle by a map of derived spaces $\EE \to \MM$ and the fields by sections thereof.
The theory of differential equations admits a generalization too (perhaps most cleanly by using the perspective of $D$-modules and its outgrowths).
In practice, as we'll see in this paper, many of the moduli spaces are variants of spaces already studied by geometers.
In fact, they are often ``derived enhancements'' of well-known moduli spaces: we work with a derived stack whose truncation is a stack already known.\footnote{In \cite{STV} a number of important moduli spaces, like moduli of stable maps of curves or moduli of perfect complexes, receive derived enhancements.}
For example, holomorphic Chern-Simons theory is best viewed as studying the derived enhancement of the moduli space of $G$-bundles on a Calabi-Yau 3-fold~$M$, where $G$ is a reductive complex algebraic group.

There is an important feature that we have not discussed so far:
the equations of motion arise from ``variational calculus.''
Loosely speaking, these equations identify the critical points of some function on the space of fields.
When approached from a derived perspective, 
the space of solutions typically obtains a natural $-1$-shifted symplectic structure.\footnote{The derived critical locus of a function $f$ on an ordinary smooth variety $X$ is the derived intersection inside $T^*X$ of the graph of $\d f$ and the zero section. The structure sheaf can be modeled by functions $\cO(T^*[-1]X)$ equipped with a differential determined by~$f$. The shifted symplectic structure on $T^*[-1]X$ carries over to this derived critical locus. See \cite{VezdCrit} for more discussion.}
Hence derived field theories lead to $-1$-shifted symplectic derived stacks,
and we have examples like $T^*[-1]\Bun_G(M)$ or~$T^*[-1] \RR\Map(M,X)$.

\subsection{Perturbative field theory as derived deformation theory}
\label{sec: pert dfn}

We now turn to giving a mathematical definition that incorporates insights from derived geometry.
We will stick with studying perturbative field theory: 
we fix a solution $\phi_0$ to the equations of motion and study formal deformations of it.
In other words, we want to study the formal derived geometry of the moduli of solutions to the equations of motion.

Here we can exploit the heuristic principle that every derived deformation problem (in characteristic zero) is described by a dg Lie algebra.\footnote{Or, sometimes more conveniently, by an $L_\infty$ algebra.}
That is, if $\cM$ is a derived stack and $p: \Spec(\CC) \to \cM$ is a point, 
then the formal neighborhood $\cM^\wedge_p$ of $p$ is modeled by the Maurer-Cartan locus of some dg Lie algebra $\fg_p$.\footnote{More technically, the shifted tangent complex $\TT_p \cM[-1]$ comes equipped with a homotopy-coherent Lie algebra structure, 
and this determines a functor from small (i.e., artinian local) nonnegatively-graded dg algebras to the $\infty$-category of spaces.}
This idea goes back to Deligne, Drinfeld, and Schlessinger-Stasheff, and perhaps farther, 
but Hinich, Lurie, and Pridham have given a precise articulation of the {\em fundamental theorem of derived deformation theory}:
there is an equivalence of $\infty$-categories between formal moduli spaces and dg Lie algebras.

Applying this perspective to classical field theory was prominently and successfully advocated by Kontsevich, 
and it is an essential component of the formalism developed by Costello.
The idea is that if want to study the perturbation theory around some solution~$\phi_0$ on the whole spacetime manifold~$X$,
then on each open set $U$ in~$X$,
there is a formal moduli space describing deformations of the restriction~$\phi_0|_U$.
In other words, there is a sheaf on $M$ with values in formal moduli spaces,
and this sheaf describes the perturbative field theory.

By the fundamental theorem, we can describe this sheaf of formal moduli spaces by supplying a sheaf of dg Lie algebras, which can often be quite convenient.
We follow this approach here.

\begin{dfn} 
A {\em holomorphic local Lie algebra} on a complex manifold $X$ is
\begin{itemize}
\item a graded holomorphic vector bundle $V$ on $X$ whose sheaf of holomorphic functions is $\cV^{hol}$, and
\item a holomorphic differential operator $Q^{hol} \colon \cV^{hol} \to \cV^{hol}$ of cohomological degree one and a holomorphic bidifferential operator $[\cdot,\cdot] \colon \cV^{hol} \times \cV^{hol} \to \cV^{hol}$. 
\end{itemize}
This data is required to endow $(\cV^{hol}, Q^{hol}, [\cdot,\cdot])$ with the structure of a sheaf of dg Lie algebras.
\end{dfn}

In quantum field theory, it is essential that the operator $Q^{hol}$ be elliptic.
We will not explicitly mention any consequences of ellipticity in this note, so we will not emphasize this condition.

Here is a quick example, as motivation.
Let $X$ be a Calabi-Yau 3-fold equipped with a holomorphic principal $G$-bundle $P \to X$.
There is an associated adjoint bundle $\fg_P \to X$, and hence a coherent sheaf $\cG_P$ of holomorphic sections of the adjoint bundle.
Note that $\cG_P$ is a holomorphic local Lie algebra, using the fiberwise Lie bracket,
and it describes the formal neighborhood of $P$ inside the space $\Bun_G(X)$.
In other words, it encodes the perturbative part of holomorphic Chern-Simons theory of Example~\ref{eg: holCS}.

We now rephrase this data in the style of complex differential geometry and then further into more physical terminology.

Any holomorphic vector bundle has an associated Dolbeault complex $\Omega^{0,\bu}(X,V)$. 
If $V$ is a holomorphic local Lie algebra, 
then $\Omega^{0,\bu}(X, V)$ is a dg Lie algebra. 
The differential is $\dbar + Q^{hol}$, 
and the bracket is given by extending the bracket on $\cV^{hol}$ with the wedge product of differential forms. 
Notice that this dg Lie algebra has the appealing feature that it is the $C^\infty$-sections of a smooth vector bundle on~$X$.

In the case of our running example, we obtain the dg Lie algebra
\[
\Omega^{0,\bu}(X,\fg_P)
\]
from the adjoint bundle of the holomorphic principal $G$-bundle on our Calabi-Yau 3-fold~$X$.
A Maurer-Cartan element in this dg Lie algebra is a $(0,1)$-form $A$ such that
\[
\dbar_P A + \frac{1}{2}[A,A] = 0,
\]
meaning that $\dbar_P + [A,-]$ defines another $\dbar$-connection, and hence holomorphic structure, on the smooth bundle $P \to~X$.

The take-away message is that the Dolbeault complex of a holomorphic Lie algebra leads to a PDE, by taking the Maurer-Cartan equation,
and this PDE should be seen as the equation of motion.

In this example the PDE arises from the variations of the action functional
\[
\int_X \Omega \wedge \left( \frac{1}{2}\langle A, \dbar A\rangle + \frac{1}{3!}\langle A, [A,A]\rangle \right).
\]
The PDE does not involve the choice of a holomorphic volume form $\Omega$ or a pairing on the Lie algebra $\fg$, 
but the action functional does.

In general, for a holomorphic local Lie algebra to arise from a field theory,
it is necessary to provide more data so as to be able to write down an action functional.

\begin{dfn}
A {\em local cyclic structure} on a holomorphic local Lie algebra is a degree $-3$ map
\[
\langle -,-\rangle \colon V \otimes V \to K_X
\]
of holomorphic vector bundles\footnote{Here $\otimes$ denotes the Whitney tensor product of vector bundles on $X$.} such that 
\begin{itemize}
\item if $f, g$ are holomorphic sections of $V$, then 
\[
\langle Q^{hol}f, g \rangle + (-1)^{|f|} \langle f, Q^{hol} g \rangle = 0
\]
so that the pairing is a cochain map, and
\item if $f, g, h$ are holomorphic sections of $V$, then 
\[
\langle f, [g,h] \rangle = (-1)^{|f|} \langle [g,f],h \rangle
\]
so that the pairing is invariant for the bracket.
\end{itemize}
\end{dfn}

Given a local cyclic structure, there is an action functional
\[
\int_X \frac{1}{2} \langle f, \dbar f + Q^{hol} f \rangle + \frac{1}{3!} \langle f, [f,f] \rangle 
\]
whose equation of motion is the Maurer-Cartan equation
\[
\dbar f + Q^{hol}f + \frac{1}{2} [f,f] = 0.
\]
Thus a holomorphic local Lie algebra with a local cyclic structure models a perturbative holomorphic field theory.

The reader might notice that such Maurer-Cartan equations only ever involve quadratic and linear terms in the field.
To incorporate cubic and higher terms in the equations of motion,
one should define holomorphic local $L_\infty$ algebras and local cyclic structures,
which are a modest and straightforward generalization, if one knows what $L_\infty$ algebra means.
Such definitions can be found in~\cite{CG2}.

%


\subsection{From supersymmetric field theories to holomorphic field theories}


This section offers a bridge from theories of physical interest to holomorphic field theories,
and hence may be of particular interest to physicists.
Those without a strong interest in physical examples can skip this section without issue.

In his work on Donaldson invariants of four-manifolds, Witten introduced the idea of a topological twist of a supersymmetric field theory \cite{WittenTwist}.
The topological twist of a supersymmetric field theory (when it exists) can be attacked mathematically using the path integral.
A key feature, however, of topological field theories is that they are amenable to a functorial paradigm.
This firmly establishes their meaning and importance within pure mathematics.
Slightly after Witten's work, the first shadows of what would later be called \textit{holomorphic twists} of supersymmetric field theories appeared, see \cite{Johansen:1994aw,Johansen:1994ud,NekThesis,Losev:1996up,Losev:1995cr,Johansen:2003hw,Eager:2018oww}. 
Originally, this work did not garner the same attention as the parallel works on topological twists; especially in mathematics.
Much later, Costello introduced the precise notion of a holomorphic twist, and focused on examples in dimensions two and four \cite{CosSUSY}.
The characterization of all holomorphic, topological, and mixed holomorphic-topological theories of Yang--Mills type in dimensions $2 \leq d \leq 10$ appears in \cite{ESW}.
Holomorphic twists of more exotic quantum field theories, such as theories of supergravity and higher form gauge fields, have recently been studied in \cite{CLsugra,SWtensor,SWpure,Raghavendran:2021qbh,Eager:2021ufo,Garner:2023wrc}.

\subsubsection{Holomorphicity as a symmetry property}

When physicists describe a classical field theory,
they typically work with the ``spacetime'' manifold $M = \RR^n$,
which is genuinely a spacetime when equipped with the Minkowski (pseudo)metric and is a space (with no time) when equipped with the Euclidean metric.\footnote{A theory on Minkowski spacetime is called {\it Lorentzian} and on Euclidean space is called {\it Euclidean},
although people also use those terms for theories on nonlinear manifolds with Lorentzian pseudometrics or Riemannian metrics, respectively.}
There is, in either case, a group of isometries (the Poincar\'e group or Euclidean group, respectively),
and it is useful to work with theories that are equivariant for this group.
That means the fields are sections of equivariant bundles,
and the equations of motion must be equivariant as well.
In fact, this ``covariance'' is typically viewed as a condition required of a good theory.

It is possible to identify holomorphic theories on $\RR^{2n}$ in terms of this covariance.
In fact, it suffices to focus on how the subgroup of translations acts, as follows.

Equip $\RR^{2n}$ with a complex structure.
Then the complexified Lie algebra of translations decomposes into holomorphic and anti-holomorphic derivatives:
\[
\CC \otimes_{\RR} \RR^{2n} 
= \CC\left\{\frac{\partial}{\partial z_1}, \ldots, \frac{\partial}{\partial z_n}\right\} 
\oplus \CC\left\{\frac{\partial}{\partial \zbar_1}, \ldots, \frac{\partial}{\partial \zbar_1}\right\}.
\]
We would like to say a translation-equivariant theory is {\em holomorphic} if the anti-holomorphic vector fields $\del / \del {\zbar_i}$ act trivially.
More carefully, we can say there is Lie algebra map from real translations to vector fields on on the space of solutions to the equations of motion (i.e., to derivations on the algebra of observables),
and when complexified, every anti-holomorphic translation $\partial/\partial \zbar_j$ acts by zero.
In the derived world this condition is encoded by the data of chain homotopy trivialization of the action by the anti-holomorphic translations. 

Note that for the notion introduced in Section~\ref{sec: pert dfn},
it is straightforward to provide such a homotopical trivialization: 
the fields are encoded by a Dolbeault complex on $\CC^n$, 
a real translation $\partial/\partial x_j$ acts by the Lie derivative $L_{\partial/\partial x_j}$,
and so we trivialize the action $L_{\partial/\partial \zbar_j}$ by using the contraction $\iota_{\partial/\partial \zbar_j}$,
thanks to Cartan's formula.

We now turn to using this point of view in the setting of supersymmetric theories.

\subsubsection{Supersymmetry and twisting}

The goal of this section is to show that any even dimensional supersymmetric field theory produces a holomorphic field theory through a process called {\em twisting}.
This process provides us with a wealth of examples that are intimately connected to theories of genuine interest in physics.\footnote{Strictly speaking, twisting requires the existence of of more than two chiral supercharges.Thus, two-dimensional theories with $\cN = (n,m)$ supersymmetry, and $n,m \leq 1$, do not admit holomorphic twists.}
As we will explain later, many interesting phenomena for supersymmetric theories have analogs in holomorphic field theory.
 
By definition, a supersymmetric field theory on $\RR^d$ is a theory that is acted upon by a super Poincar\'{e} algebra. 
We focus on Euclidean field theories and work in Riemannian signature,
so for us the ``super Poincar\'{e} algebra'' is a super Lie algebra of the form
$\mathfrak{so}(d) \ltimes \ft$
where $\ft$ is the super Lie algebra of {\em supertranslations} whose even part $\ft^0 = \RR^d$ is the Lie algebra of ordinary translations and whose odd part $\ft^1$ is a sum of spin representations. 
While the Lie algebra $\RR^d$ of ordinary translations is abelian, the Lie algebra $\ft$ carries a nontrivial (super) Lie bracket,
and this Lie bracket is defined in terms of a $\mathfrak{so}(d)$-equivariant non-degenerate symmetric pairing
\beqn
\label{e:Gamma}
\Gamma \colon {\rm Sym}^2(\ft^{1}) \to {\rm Sym}(\ft^0) \cong \RR^d 
\eeqn
by the formula $[\cQ, \cQ'] = \Gamma(\cQ, \cQ')$. 

By a {\em supercharge}, one means an odd supertranslation $\cQ \in \ft$, and
a {\em twist} is a (nonzero) supercharge $\cQ$ such that $[\cQ,\cQ] = \cQ^2 = 0$. 
The classification of twisting supercharges is a completely algebraic question:
see~\autocite{EllSaf} for a complete classification in dimensions from $1$ to~$10$.

In even dimensions $d = 2n$, basic properties of the supersymmetry algebra guarantee that there exists a twist that provides a complex structure to $\RR^{2n}$.\footnote{This claim fails in exactly one example, namely for $\cN=(1,0)$ supersymmetry in two dimensions, as the dimension and amount of supersymmetry are too small.}
Indeed, so long as a twisted supercharge $\cQ$ is nonzero, 
one finds that the image of the map $[\cQ,-] = \Gamma(\cQ,-)$ from \eqref{e:Gamma} is at least $n$-dimensional.
A {\em minimal} or {\em holomorphic} supercharge is one for which the dimension is exactly~$n$,
and we declare the image of $\cQ$ as spanned by the anti-holomorphic vector fields~$\{\del / \del \zbar_i\}$. 

If a classical field theory $\cT$ has $\ft$ as a symmetry, 
then a choice of twist $\cQ$ determines a deformation of the theory that we call a {\em twisted supersymmetric theory} $\cT^\cQ$.
If $S$ denotes the action functional of $\cT$, then the action functional $S^\cQ$ of $\cT^\cQ$ has the form
\[
S^\cQ(\varphi) = S(\varphi) + \int_{\RR^d} \varphi  \left(\cQ \cdot \varphi \right) + \cdots 
\]
where the deformation arises from how $\cQ$ acts on the theory 
(the $\cdots$ leaves room for terms non-linear in $\cQ$).
For a systematic discussion of how the twisted theory is defined, 
we refer to~\autocite{CosSUSY, ESW}. 

Building upon the physical literature, Elliott, Safronov, and the second author have given a complete characterization of the all twisted supersymmetric Yang--Mills theories in dimensions $2 \leq d \leq 10$ \autocite{ESW}. 
These twisted theories are BF theorires or Chern--Simons theories, possibly along with other fields in associated bundle; 
depending on the twist, that might mean a purely topological or purely holomorphic BF or Chern-Simons theory.
In even dimensions $d = 2n$, a minimal twist renders the theory holomorphic.
$\cQ$ itself determines the homotopies $\{\eta_i\}$ such that 
\beqn
[\cQ,\eta_i] = \frac{\del}{\del \zbar_i} ,
\eeqn 
recovering the characterization of holomorphic field theories mentioned at the beginning of the subsection. 

As a representative example of twisting, 
consider four-dimensional space $\RR^4$ with $\cN=1$ supersymmetry.
The odd part of the supersymmetry algebra $\ft^1$ is four-dimensional, 
and every nonzero twist is minimal and gives us a copy of $\CC^2$.\footnote{And in a natural sense, they are all equivalent.}
That is, with this amount of supersymmetry, there does {\it not} exist a topological twist.
consider the fundamental example of pure four-dimensional supersymmetric Yang--Mills theory with $\cN=1$ supersymmetry. 
Its holomorphic twist is equivalent to holomorphic BF theory on~$\CC^2$ described in Section~\ref{sec: hol BF dfn}.

There are many other appearances of holomorphic field theories in the context of supersymmetry.
Here is a (very) non-exhaustive list:
\begin{itemize}
\item Topological string theories yield rich classes of holomorphic field theories, as target space field theories. 
For instance, the topological B-model on a Calabi--Yau threefold $X$ leads to Kodaira--Spencer theory on $X$, as proposed in \cite{BCOV} using the framework of closed string field theory. 
The perspective of Kodaira--Spencer theory as a holomorphic field theory has been explored with great success by Kevin Costello and Si Li \cite{CLbcov, CLtypeI}.
\item
It is possible to make sense of Kodaira--Spencer theory outside of three complex dimensions.
In \cite{CLsugra}, Costello and Li have argued for conjectural descriptions of twists of supergravity and superstring theories in terms of this more general version of Kodaira--Spencer theory.
Holomorphic field theory has also appeared in applications to (twisted) holography through the lens of Koszul duality, see \cite{CLsugra, CGholography, CostelloPaquette}.
\item The pure spinor formalism of Berkovits and Cederwall  \cite{Cederwall:2001dx, Berkovits:2005bt}, see also \cite{Eager:2018dsx,perspectives}, starts with a field theory or string theory defined on the (sometimes singular) nilpotence variety of the super Poincar\'e algebra and then uses the action of supersymmetry to produce ordinary (non-twisted) supersymmetric field theories. 
Often, the theory on the nilpotence variety is holomorphic \cite{ESW,SWpure}. 
\item Similar in spirit to the pure spinor formalism is a systematic relationship between holomorphic field theories on (super) twistor space and ordinary (Riemannian) field theories.
Such work can be viewed as extending (and enhancing) the original Penrose correspondence which relates holomorphic BF theory on (super) twistor space associated to $\RR^4$ to the self-dual limit of (super) Yang--Mills theory \cite{Penrose:1969ae, Penrose}---see \cite{adamo2018lecturestwistortheory} for a modern review.
See \cite{Mason, Mason_2009, CostelloTwistor,Bittleston_2023,garner2024scatteringtwistoriallinedefects} for recent work on the twistor correspondence in QFT.
\end{itemize}

We recommend \cite{bah2024panoramaphysicalmathematicsc} for more discussion and references.

\section{Algebra and holomorphic field theories}
\label{sec: algebra}

In this section we describe  associative algebras and their ``homotopical'' generalizations, like $A_\infty$ algebras, related to holomorphic field theories,
following the relationship of associative algebras to quantum mechanics and of vertex algebras to chiral conformal field theory.
A particular focus is upon higher Kac-Moody algebras \cite{FHK},
as these offer a tantalizing direction to explore in search of analogs of the rich connections between representation theory, algebraic geometry, and physics familiar to those who have worked with loop groups.
In a later section --- and it is a central point of this survey --- we will explain how factorization algebras provide a direct conduit from holomorphic field theories to these algebraic constructions.
Throughout this section, we will be a bit cavalier with certain subtleties (e.g., about infinite-dimensional spaces and functions on them), 
emphasizing concepts and motivations over mathematical precision.
We will be a bit long-winded in the first two subsections, to provide motivations,
so that in the final subsection, we can rapidly sketch these higher Kac-Moody and Weyl algebras.

\subsection{Algebras in mechanics}

A key feature of quantum mechanics is that the observables (or operators) live in an associative algebra.
In many cases this associative algebra is a deformation of a commutative algebra, typically arising as functions on a manifold or variety.
The quintessential example is the Weyl algebra
\[
\CC \langle x, p \rangle /(xp - px = i \hbar),
\] 
which is generated by observables $x$ (``position'') and $p$ (``momentum'') for a quantum particle moving along a line.
There is a parameter $\hbar$, which if sent to zero, recovers a commutative algebra $\CC [ x, p ]$ of complex-valued polynomial functions on the cotangent bundle $T^* \RR$ of the real line~$\RR$.
This example will be a model for much of what we discuss in this paper.
So far we have ignored a lot of features of quantum mechanics (e.g., $\ast$-structures, Hilbert spaces, unitarity) that are important in physics,
and we will continue to do so.
Note as well that there are many more observables and we have only discussed subalgebras of all observables (at the classical level, just polynomials in $p$, $q$);
we will often focus on such tractable subalgebras. 

There are some features of this example that we would like to foreground. 
First, the commutative algebra arises as functions on a {\em symplectic} space,
which is the usual mathematical setting for classical mechanics (aside from more subtle situations that require Poisson geometry).
Second, the symplectic form $\omega = \d x \wedge \d p$ equips this algebra 
with a Poisson bracket where $\{ x, p \} = 1$,
which controls the deformation to the Weyl algebra: following Dirac, we promote the Poisson relation to a commutator relation.
These two features motivate the deformation quantization problem: 
given a symplectic (or Poisson) manifold,
describe deformations of its commutative algebra of functions to an associative algebra with the requirement that, to first order,  the commutator recovers the Poisson bracket.
Thanks to Kontsevich \cite{KonDQ}, there is a beautiful answer to this question, which has spawned a mountain of fascinating mathematics.

Our view on field theory is motivated by this perspective on quantum theory (we discuss it further in Section~\ref{s:factinqft}), 
and it might help the reader to bear in mind a variant of this question:
given a classical holomorphic field theory, what are the natural deformations of its algebra of functions?
Many of the results we discuss later can be seen as generalizations of phenomena associated with deformation quantization.

In particular, consider how symmetries appear in mechanics.
Let $M$ be a symplectic manifold encoding the ``phase space'' of a classical mechanical system,
and let $\{-,-\}$ denote the Poisson bracket on $C^\infty(X)$.
If a Lie group $G$ acts on $X$ by symplectomorphisms (i.e., respects the symplectic structure),
there is a map of Lie algebras $\rho \colon\fg \to {\rm SympVect}(M)$ so that an ``infinitesimal symmetry'' $x \in \fg$ acts by a symplectic vector field.
In many cases, such infinitesimal symmetries are realized as observables: 
there is a map of Lie algebras $H_\rho \colon \fg \to C^\infty(M)$ such that $\{ H_\rho(x), - \} = \rho(x)$.
One says that each element $x$ has a Hamiltonian function $H_\rho(x)$ whose associated Hamiltonian vector field is $\rho(x)$.
Classic examples include the momenta for $T^* \RR^n$ arising from translation and rotational symmetry.

How does this set-up fit into the deformation quantization view?
First, note that the free commutative algebra $\Sym(\fg)$ has a canonical Poisson bracket by defining $ \{x, y \}   = [x, y]$ for generators $x, y \in \fg$ and extending by the Leibniz rule.
Hence the Hamiltonian map $H_\rho$ extends to a Poisson map $H_\rho: \Sym(\fg) \to C^\infty(M)$.
We can then ask: does this map quantize? That is, can we compatibly deform the domain and range of the maps to interesting associative algebras as well as deforming to a map of associative algebras?
Here it might be useful to recognize that the enveloping algebra $U\fg$ is a natural deformation quantization of $\Sym(\fg)$, thanks to the Poincar\'e-Birkhoff-Witt theorem.
Thus our question might be reformulated to asking whether there is an associative algebra map $H^q_\rho: U\fg \to C_\hbar^\infty(M)$, where $C_\hbar^\infty(X)$ is a deformation quantization, such that $H^q_\rho$ recovers $H_\rho$ in the classical limit.
Such a map is sometimes called a {\em quantum moment map},
as $H_\rho$ is the pullback of (polynomial) functions along a moment map~$\mu_\rho \colon M \to \fg^*$.

This class of questions appears throughout mathematics,
and it has played a key role in geometric representation theory and the theory of $D$-modules, where $D$ denotes a ring of differential operators.
After all, for any manifold $M$, the algebra of differential operators $D_X$ is a natural deformation quantization of functions on $T^* M$.
Hence, for any $G$-manifold $M$, it is natural to ask whether there is a representation $\fg \to D_M$ deforming the canonical Poisson representation.
In this spirit, we will search for analogues of the kind of mathematics that has grown out of results by Beilinson-Bernstein \cite{BB}, Brylinski-Kashiwara \cite{BryKash}, Kostant \cite{Kostant} and others.

\begin{rmk}
We have only spoken of symplectic manifolds, and emphasized vector spaces,
but fermions play a crucial role in quantum mechanics too.
Here the symmetric algebra on a symplectic vector space (as classical observables) 
is replaced by an exterior algebra on an inner product space,
which can be seen as a symmetric algebra on an {\em odd} vector space
(i.e., a super vector space with purely odd component).
The deformation quantization of the exterior algebra is then a Clifford algebra.
We will set fermions aside for the remainder of this section,
but the discussion below extends easily to them and fermions appear naturally in the study of holomorphic field theories.
\end{rmk}

%

\subsection{A loopy version: the algebra of modes}

What could provide the {\em holomorphic} analog of the story above?
Let us indicate one possible answer before we motivate it from the point of view of holomorphic field theory.
Roughly, the idea is that associated to any finite-dimensional symplectic space is another, infinite-dimensional, symplectic space of ``loops'' into the original symplectic space.


Let $X$ be an algebraic variety defined over $\CC$, equipped with a symplectic form $\omega_X$.
Then we posit for the relevant phase space, the ``symplectic space'' $LX$---the formal loop space of $X$.
This is the space of maps $\Hat{D}^\times \to X$ from the formal punctured disk to $X$.
%
(There are clearly variations on this idea, such as the holomorphic loop space of all holomorphic maps from $\CC^\times$ to~$X$, when $X$ is a complex manifold.)

When $X=V$ is a complex vector space equipped with a linear symplectic form $\omega_V$, we can identify $LX$ with the vector space $V((z))$ of $V$-valued Laurent series in a single variable.
The symplectic form on this vector space is given in terms of the residue
\beqn
\omega_{LV} (f,g) = \underset{z=0}{\text{Res}} \, \omega_V(f,g) \d z .
\eeqn
One can get rid of the appearance of the holomorphic line element $\d z$ if we `twist' Laurent series by the square-root element $\d z^{1/2}$---in other words, we take sections of the line bundle $K^{1/2}_{D^\times}$.
(Again, there are clearly alternative algebras to consider.)
When we quantize, this algebra will be related to a well-known vertex algebra.

There is a feature we would like to point out, before we turn to explaining our interest in this space and algebra.
The inclusion $\Hat{D}^\times \hookrightarrow \Hat{D}$ means that any algebraic map  $\phi: \Hat{D} \to V$ restricts to an algebraic map $\phi \colon \Hat{D}^\times \to V$.
In our notation, this restriction map is the inclusion $V[[z]] \hookrightarrow V((z))$ of $V$-valued power series in $V$-valued Laurent series.
When we quantize, this relationship will produce the underlying vector space (or Fock space, or vacuum module) of the vertex algebra.

Now we will discuss why we might focus on these constructions,
and in what sense they are holomorphic versions of usual mechanics.

Let us start by offering a motivation for the ``answers'' (e.g., $\RR[q,p]$ as classical observables) in the setting of standard classical mechanics.
The model problem in mechanics is to describe a point particle moving in $\RR^n$,
subject to some forces that specify the differential equations governing the particle's motion.
(We call these the ``equations of motion.'')
Newton's law says that this equation is second order, so that a solution (i.e., trajectory) is specified by giving the position and velocity of a particle at one instant in time.

There is another view, which generalizes more naturally to field theories and which is dubbed the Lagrangian formalism.
Here one notes that there is a space of all imaginable trajectories, namely the path space ${\rm Map}(\RR, \RR^n)$, where we view the source $\RR$ as ``time,''
and there is a subspace of trajectories realized by the particle, namely solutions to the equations of motion.
This subspace is, in some sense, the critical set of an action functional; 
the variational calculus provides the relevant mathematical framework.
The space of solutions is naturally isomorphic to $T \RR^n$,
where an isomorphism is given by fixing an instant $t_0$ in time $\RR$ and then sending a solution $\phi$ to the pair $(\phi(t_0), \dot{\phi}(t_0))$, the position and velocity of the solution at~$t_0$.
With a little care, one finds that the variational calculus equips $T \RR^n$ with a natural symplectic form and a symplectomorphism $T\RR^n \cong T^* \RR^n$, with its canonical symplectic form. 
Putting everything together, we have a natural symplectomorphism of $T^* \RR^n$ with the space of solutions.
The observables of the classical system are the algebra of functions on the space of solutions,
and so this isomorphism tells us that functions on $T^* \RR^n$ provide the observables.
The polynomial functions on $T^* \RR^n$ are a subalgebra of all observables, 
and they suffice to distinguish distinct solutions.

This view generalizes nicely to the holomorphic setting,
and it amounts to using the Lagrangian formalism to study the holomorphic field theories we have already introduced.
In this model case, we replace the path space ${\rm Map}(\RR, \RR^n)$ by ${\rm Map}(\CC^\times, \CC^n)$ or, if one wishes, ${\rm Map}(S, \CC^n)$ with $S$ a Riemann surface.
We start with $\CC^\times$ since it has a ``time'' direction given by the radial coordinate ($t$ becomes $r = e^t)$).
Letting us view $\CC^\times$ as a cylinder $\RR \times S^1$, we can view the mapping space as  ${\rm Map}(\RR, \Map(S^1,\CC^n))$, namely mechanics into the loop space of $\CC^n$.
We replace Newton's equations of motion with a holomorphic version
so that the space of solutions is given by holomorphic maps from $\CC^\times$ into $\CC^n$.
If we wish to focus on a more algebraic version, as we did at the beginning,
we could restrict to the algebraic maps from $\CC^\times$ into~$\CC^n$,
which sit inside the holomorphic maps.

A confession is necessary here, because this replacement is a bit misleading.
The holomorphic version of the model case, known as the free $\beta\gamma$ system,
actually has holomorphic maps into $T^* \CC^n$ as the space of solutions.
The algebraic maps are precisely $T^* \CC^n$-valued Laurent polynomials.\footnote{One can extract the algebraic model from the holomorphic model by taking formal power series expansions at the level of observables.}

We now turn to quantization, and our phase space will be, for simplicity, the formal loop space $V((z))$ of a symplectic complex vector space $V$.
Following the case of mechanics, we might ask for a deformation of the Poisson algebra
\beqn
\cO\left( V((z)) \right) .
\eeqn
Like the symplectic form on $V((z))$, the Poisson bracket is given in terms of the residue pairing.
The standard Heisenberg commutation relations determine a Weyl algebra $W[V]$ for this formal loop space.
As a vector space $W[V]$ is isomorphic to $\cO(V((z)))$, 
just as the simplest Weyl algebra is isomorphic to $\cO(V)$ as a vector space,
but the product structure is modified by quantization.
We expect, by analogy with the simplest case, that any quantization of a bigger algebra (e.g., for observables of the holomorphic loop space) contains this quantization as its algebraic skeleton.

The subspace of ``contractible formal loops" $V[[z]]$ inside $V((z))$ produces a module for $\cO(V((z)))$ given by $\cO(V[[z]])$. 
(This module is in fact a quotient algebra.)
One can quantize this module to a module ${\rm Vac}[V]$ for $W[V]$,
where as a vector space, it is still $\cO(V[[z]])$.
This module structure means there is a map of algebras $W[V] \to \End({\rm Vac}[V])$,
so we have a natural ``Hilbert space'' or Fock module for~$W[V]$.

It is a remarkable fact, suggested by physicists, that there is also a map
\[
Y \colon {\rm Vac}[V] \to \End({\rm Vac}[V])[[z, z^{-1}]],
\]
known as the {\em vertex operator} or {\em state-field correspondence},
which is closely related to the $W[V]$-action.
The intuition behind $Y$ is that for any map $\phi: \CC^\times \to V$ and for any point $w \in \CC^\times$, 
there is a restriction of $\phi$ to a little disc around $w$.
That is, there is a restriction map $r_w: V((z)) \to V[[t]]$, where $t$ is the local coordinate $t = z-w$.
Hence there is a $w$-dependent map $\cO(V[[t]]) \to \cO(V((z)))$ by pulling back a function along $r_w$;
this map determines a $w$-dependent action of $\cO(V[[t]])$ on $\cO(V[[z]])$.
One can, in essence, quantize this map, by trying to extend the formulas of the $W[V]$-action,
and this leads to the map~$Y$.
If one axiomatizes the behavior of $Y$, one is led to the notion of a vertex algebra.
The algebra $W[V]$ is known as the \textit{algebra of modes} associated to this vertex algebra.

Let us briefly comment on how symmetries extend to the holomorphic setting.
In other words, we wish to explain how the affine Lie algebras $\widehat{\fg}$ and free field realizations arise in analogy with our discussion of quantum moment maps.

Suppose now that $G$ is a complex Lie group and it acts on the symplectic complex vector space $V$ by holomorphic symplectomorphisms.
Then $G$ also acts on the formal loop space $V((z))$ as a ``global'' symmetry:
the action is independent of $z$.
On the other hand, if we consider the infinitesimal action $\rho \colon \fg \to {\rm SympVect}(V)$,
there is, in fact, a natural extension to a ``local'' symmetry by the ``loop algebra'' 
\beqn
\rho^{\rm loc} \colon L \fg = \fg ((z)) \to {\rm SympVect}(V((z))) .
\eeqn
Geometrically, a $\fg$-valued function $x(z)$ acts pointwise (with respect to the formal disc) on the formal loop space.
This action factors through a map $H_\rho^{\rm loc}: \fg((z)) \to \cO(V((z)))$,
and hence is Hamiltonian.
Thus the loop algebra $\fg((z))$ arises naturally as symmetries of a classical holomorphic field theory.

When we try to quantize, we run into an interesting phenomenon:
the Hamiltonian action $H_\rho^{\rm loc}$ does not extend as a Lie algebra map into the Weyl algebra $W[V]$,
although it makes sense as a linear map.
The failure to respect the brackets determines, however, a cocycle on the loop algebra
by
\[
\alpha_V( x \otimes f(z) , y \otimes g(z)) = \Tr_V(xy) \Res(f \, \partial_z g).
\]
There is thus a central extension of the loop algebra 
\[
\CC K\to \widehat{\fg} \to \fg((z))
\]
so that 
\[
[x \otimes f, y \otimes g] = [x,y] \otimes fg + \alpha_V( x \otimes f , y \otimes g) c,
\]
and there is a Lie algebra map
\beqn\label{eqn:freefield1}
\widehat{H}^{\rm loc}_\rho \colon \widehat{\fg} \to W[V] .
\eeqn
This is an example of a {\em free field realization} for the affine Lie algebra~$\widehat{\fg}$.
This relationship --- and its extension to a map between the vertex algebras associated to $\widehat{\fg}$ and $V$ --- is small part of the rich dialogue that has developed between complex analysis, representation theory, and the physics of chiral conformal field theory.
(To get Riemann surfaces into the game, we will need a richer framework than ordinary algebra.)

\subsection{Reduction vs. chiralization}
\label{sec: red vs chiral}

{\it Dimensional reduction}, an idea from physics, offers a view on how the Weyl algebra on generators $p, q$ relates to its loopy version.
As we already discussed, mechanics (the Weyl algebra of $V$) has to do with a one-dimensional system, i.e., time-line, 
while its loopy version ($W[V]$) has to do with a two-dimensional system, such as the punctured plane~$\CC^\times$.
There is a familiar identification (e.g., polar coordinates)
\[
\CC^\times \cong \RR \times S^1
\]
and so
\[
\Map(\CC^\times,V) \cong \Map(\RR, \Map(S^1,V)) = \Map(\RR, LV).
\]
Dimensional reduction here means studying the two-dimensional physics as a one-dimensional system with a more complicated target, namely the loop space.\footnote{One can play this game whenever the spacetime manifold is product $X \times Y$ or a fiber bundle.
Sometimes the term {\it compactification} is used instead, with reduction reserved for fixing the fiberwise fields to be constant---for a discussion of compactification within the setting of factorization algebras see Section~\ref{s:compact}.}
There is, in this case, a kind of dual process that we call {\it chiralization},
that turns a one-dimensional physical system (mechanics) into a holomorphic field theory.
We will now describe this relationship in the language of fields and Lagrangians.

The ordinary Weyl algebra on a symplectic vector space $V$ arises from the quantization of ordinary mechanics on the symplectic vector space $V$.
A Hamiltonian $H \in \cO(V)$ acts as a derivation of the Weyl algebra by the commutator $[H,-]$.
The fields of this system are smooth maps $\phi = \phi(t)$ valued in $V$, 
and the action functional is 
\[
S_{1d} (\phi) = \int \omega(\phi , \d \phi) + H(\phi) \d t ,
\]
where $\d \colon C^\infty(\RR) \otimes V \to \Omega^1(\RR) \otimes V$ is the de Rham differential. 

The chiralization of this setup goes as follows.
Let a field $\phi = \phi(z,\zbar)$ be a smooth function valued in $V$, 
where $(z,\zbar)$ are (anti)holomorphic coordinates on the manifold $\CC^\times$.
The action functional is 
\[
S_{2d}(\phi) = \int \omega(\phi, \dbar \phi) \d z ,
\]
where $\dbar \colon C^\infty(\CC^\times)\otimes V \to \Omega^{0,1}(\CC^\times) \otimes V$ is the Dolbeault differential.

{\it Reduction} goes from this holomorphic field theory to a mechanical system. 
Take as this one-dimensional space the radial direction in $\CC^\times$.
Concretely, any smooth function on $\CC^\times$ can be decomposed in polar coordinates as
\[
\phi(r,\theta) = \sum_{k \in \ZZ} \phi_k(r) e^{2 \pi i k \theta} ,
\]
where $\phi_k(r)$ is a smooth function depending only on the radius $r \in \RR_+$.
Alternatively, by taking the logarithm we can view $\phi_k(t)$ as a function of the ``time'' parameter $t = \log r \in \RR$. 
We interpret $\sum_{k \in \ZZ} \phi_k(r)$ as smooth map from $\RR$ to the formal loop space $LV$. 

At the level of action functionals, 
we can integrate $S_{2d}(\phi)$ over the angular coordinate to obtain the following formal sum of one-dimensional action functionals
\beqn
S_{red}\left( \ldots, \phi_{-1}, \phi_0 , \phi_1, \ldots \right) = \sum_{k\in \ZZ} \int  \phi_{-k} \,\d \phi_{k} + k\, \phi_k \, \phi_{-k} \, \d t .
\eeqn 
Notice that the $k=0$ summand is just the Lagrangian of ordinary mechanics into $V$ with zero Hamiltonian. 
In general, we find mechanics into the symplectic vector space $LV$, where the symplectic structure uses the residue pairing in addition to the symplectic form on $V$, with a nontrivial Hamiltonian that encodes the $U(1)$-symmetry arising from rotating $LV$ in the natural way.


\subsection{Algebras in higher dimensions }
\label{s:alghigh}

The story we have just told admits a natural extension to higher complex dimensions,
but it requires a first step in the direction of derived geometry.

Let $V$ denote a complex vector space (no longer required to be symplectic).
The naive idea is to replace the algebraic loop space $\Map(\CC^\times,  V)$, or its formal version $\Map(\Hat{D}^\times, V)$, with a higher-dimensional analog $\Map(\CC^d -\{0\}, V)$.
If one considers ordinary algebraic or holomorphic maps for $d > 1$, 
then every such map extends across the origin, by Hartogs' lemma,
so this seems uninteresting.
One should notice, however, that punctured affine space $\mathring{D}^d = \Hat{D}^n -\{0\}$ is, as a scheme, not affine, 
and so the {\em derived} global sections of the structure sheaf $\cO_{\rm alg}$ are interesting:
\[
H^\bu(\RR \Gamma(\mathring{D}^d, \cO)) = \begin{cases} 
0, & \bu \neq 0, d-1 \\ 
\CC[[z_1,\ldots,z_d]], & \bu = 0 \\ 
\CC[z_1^{-1},\ldots,z_d^{-1}] \frac{1}{z_1 \cdots z_d}, & \bu = d-1 
\end{cases} .
\]
We are providing a natural basis for the cohomology 
that aims to emphasize analogies with the $d=1$ case.
When $d=1$, note that this recovers the Laurent series,
so when $d > 1$, 
we view the cohomology in degree $d-1$ as providing the derived replacement of the polar part of the Laurent polynomials.
Checking this computation is a straightforward exercise in algebraic geometry;
for instance, use the cover by the affine opens that are the complements of the coordinate hyperplanes $\{z_i =0\}$.
The full cohomology has the structure of a commutative algebra. 
In degree zero, the algebra structure is the ordinary one on functions, or polynomials.
The module structure in the basis above is given by the apparent multiplication: $z_i^k \cdot z_i^{-l-1} = z_i^{k-l-1}$ if $k < l$ and zero otherwise.

There is also an algebraic version where we replace the formal disk by the affine space $\AA^d$ and its punctured version $\mathring{\AA}^d$.
A similar result holds in analytic geometry, of course,
so that we have natural maps
\[
\RR \Gamma(\mathring{D}^d, \cO) \leftarrow \RR \Gamma(\mathring{\AA}^d, \cO) \to \RR \Gamma(\mathring{\CC}^d, \cO_{\rm an}) 
\]
given by the embedding of polynomials into power series and holomorphic functions, respectively.

It is important to have an explicit dg commutative algebra that models the derived global sections for functions on punctured formal disk,
and not just the cohomology groups.
There is a beautiful, explicit dg model $\sfA^\bu_d$ for the formal version derived global sections due to Faonte, Hennion, and Kapranov \cite{FHK} and based on the Jouanolou method for resolving singularities. 
While we won't go into details about this model, we remark on a few important points all of which we take from \cite{FHK}:
\begin{itemize}
\item The complex $\sfA^\bu_d$ is concentrated in degrees $0,\ldots,d-1$ and its cohomology is the same as the cohomology of $\RR \Gamma(\mathring{\DD}^d, \cO)$ described above.
\item There is a polynomial version $\sfA_{d,poly}^\bu$ and explicit embeddings of commutative dg algebras 
\beqn
\sfA_{d}^\bu \hookleftarrow \sfA_{d,poly}^\bu \hookrightarrow \Omega^{0,\bu}(\CC^d - 0) .
\eeqn
The left embedding replaces the inclusion of Laurent polynomials into Laurent series.
The right embedding replaces the inclusion of Laurent polynomials into the Dolbeault complex on $\CC^\times$.
\item There is a higher residue map $\text{Res} \colon \sfA_d^{\bu} \to \CC[-d+1]$ that is compatible with integration over the $2d-1$ sphere
\beqn
\label{eq:res}
\Res(\alpha) = \oint_{S^{2d-1}} \alpha \wedge \d^d z,
\eeqn
under the above embeddings.
Here $\d^d z = \d z_1 \wedge \cdots \wedge \d z_d$ is the standard holomorphic volume form.
\end{itemize}



The ingredients are in place to mimic the construction of the Weyl algebra for the formal loop space, of the affine Lie algebras, and of free field realization, but now with $\mathring{\CC}^d$ in place of $\mathring{\CC}^1 = \CC^\times$.
We must be careful about parity, however. 
In dimension $d = 1$, the Laurent polynomials valued in a symplectic vector space $V((z))$ again has the structure of a symplectic vector space via the residue. 
Because of the appearance of the cohomological shift in the higher dimensional residue it is no longer the case that 
\beqn
V \otimes \RR\Gamma(\mathring{D}^d, \cO) \simeq V \otimes \sfA_d
\eeqn
is symplectic if $V$ is.

The residue map $\Res$ in equation \eqref{eq:res} is of cohomological degree $1-d$. 
Thus, to get a symplectic pairing on $V \otimes \sfA_d$ via the residue we can ask that $V$ is equipped with a {\em $(d-1)$-shifted symplectic structure}. 
A typical example of such a space is $V = T^*[d-1] L = L \oplus L^*[d-1]$, where $L$ is any vector space.
\footnote{In the generality of shifted symplectic geometry, we expect that the derived mapping space $\text{Map}(\mathring{D}^n, X)$ is symplectic provided that $X$ is a $(d-1)$-shifted symplectic space.}

Let $\omega_V$ be a $(d-1)$-shifted symplectic pairing on the graded vector space $V$. 
Then $\omega_V$ extends to a map $\til \omega_V \colon (V \otimes \sfA_2)^{\otimes 2} \to \sfA_d$ where we use the algebra structure on $\sfA_d$. 
We then obtain a ($0$-shifted) symplectic pairing on $V \otimes \sfA_d$ by the formula
\beqn
\label{eqn:omega}
\omega(a,b) = \Res ( \til \omega_V (a,b) \wedge \d^d z) .
\eeqn
As in the construction of the ordinary Weyl algebra, we can think about $\omega$ as defining two-cocycle on the (abelian) dg Lie algebra $V \otimes \sfA_d$.
In particular, $\omega$ determines a central extension of dg Lie algebras
\beqn
\CC \to \lie{h}_d[V] \to V \otimes \sfA_d .
\eeqn
From here, we can define the higher-dimensional analog of the Weyl algebra. 

\begin{dfn}
Let $V$ be a $(d-1)$-symplectic graded vector space. 
Then, the $d$-dimensional Weyl algebra $W_d[V]$ is the enveloping algebra of the dg Lie algebra $\lie{h}_d [V]$ where we set the central term equal to the unit.
\end{dfn}

There is a close relationship of the higher dimensional Weyl algebra to the higher dimensional $\beta\gamma$ system on $\CC^d - 0$ just as the loopy Weyl algebra is related to the $\beta\gamma$ system on~$\CC^\times$.
In this sense, $W_d[V]$ is an example of a `mode algebra', where the modes are now parametrized by the cohomology of punctured $d$-dimensional affine space.

The trick of turning providing a higher dimensional version of the Weyl algebra can be easily modified to obtain higher dimensional version of current, or loop, algebras. 
This is one of the starting points of \cite{FHK} in their construction of higher dimensional Kac--Moody algebras. 

\begin{dfn}[\cite{FHK}]
\label{dfn:sphere}
For a Lie algebra $\fg$, the {\em sphere algebra} in complex dimension $d$ is the dg Lie algebra~$\fg \otimes \sfA_d$.
We denote it by~$\fg^\bullet_d$.
\end{dfn}

There are natural central extensions of this sphere algebra 
parallel to the affine Lie algebras appearing as central extensions of the loop algebra.
They are given by the data of a $\fg$-invariant degree $d+1$ polynomial $\theta \in \Sym^{d+1}(\fg^*)^G$.
The cocycle that defines the $L_\infty$ central extension is
\[
\label{fhk cocycle}
\begin{array}{cccc}
\theta_{\rm FHK} : & (\sfA_d \tensor \fg)^{\tensor (d+1)} & \to & \CC\\ 
& a_0 \otimes \cdots \otimes a_d & \mapsto & \Res \left(\theta(a_0,\partial a_1,\ldots,\partial a_d)\right)
\end{array}.
\]
This cocycle has total cohomological degree $+2$ as an element in the Lie algebra cohomology of $\fg^\bu_{d}$ and so determines a central extension of dg Lie algebras
\[
\CC c \to \Hat{\fg}^\bullet_{d, \theta} \to \fg^\bullet_d .
\]
There is a particularly nice model for this dg Lie algebra as an $L_\infty$ \textit{algebra} that we use in \cite{GWkm}.
An $L_\infty$ algebra is a generalization of a dg Lie algebra;
it is a graded vector space $L$ equipped with a countable collection of ``higher brackets'' $[-]_k \colon L^{\otimes k} \to L$ that satisfy a system of equations generalizing the Jacobi identity.
For an overview (and an actual definition!) we recommend \cite{LV}.
In our setting the usual bracket (the 2-bracket) of $\fg^\bullet_d$ is unchanged but there is now a nontrivial $d$-bracket defined by~$\theta_{\rm FHK}$.

The algebras $\Hat{\fg}^\bullet_{d, \theta}$ are not just introduced on a whim.
In \cite{FHK} they are shown to appear naturally in various contexts,
notably in studying the derived moduli space of $G$-bundles on a complex varieties of dimension~$d$.

With the higher dimensional Weyl and Kac--Moody algebras in place, a natural question to ask is whether there are higher dimensional versions of the free field realization as in \eqref{eqn:freefield1}. 
If $V$ is a representation of a Lie algebra $\fg$,
there is a local symmetry 
\[
\fg^\bullet_d \to \cO(V \otimes \sfA^{\bu}_d),
\]
just as $\fg((z))$ acts on $\cO(V((z)))$. 
One can ask if this action lifts to a quantum symmetry.
The obstruction, see \cite{GWkm}, is one of the cocycles we discussed above:
\[
(x_0 \otimes f_0, \ldots, x_{d} \otimes f_{d}) \mapsto \Tr_V( x_0 \cdots x_{d}) \,\Res(f_0 \, \partial f_1 \wedge \cdots \wedge \partial f_{d}).
\]
That is, the invariant degree $(d+1)$ polynomial on $\fg$ that determines the $L_\infty$ extension of the sphere algebra $\fg_d^\bu$ is given by ${\rm Tr}_V (x_0 \cdots x_{d})$.
In this sense, the higher Kac-Moody algebras admit a higher free field realization,
but a richer story is possible once we have a higher version of vertex algebras.
We will return to this topic once we have the rich language of factorization algebras.

\section{Algebraic structures parametrized by configuration spaces}
\label{sec: config}

Associative algebras arise in quantum mechanics, we have argued, because the spacetime of a mechanical system is a line (``time'').\footnote{In the context of holomorphic field theory on $\C^n$, we have treated the radial direction as the ``time''.}
If we organize the observables (or operators) by when they happen in time,
it corresponds to writing them left-to-right as products in an algebra.
A natural question to ask then is what kind of higher algebraic structure captures how observables behave for theories in higher dimensions?
The shortest answer is ``factorization algebras,'' and we will describe and explore this notion in Section~\ref{sec: fact alg}
but there are nice answers for two special classes of theories that we will develop now:
\begin{itemize}
\item For a topological field theory on $\RR^d$ (which we'll define in a moment),
the observables form an algebra over the {\em little $d$-discs operad}~$\Disc_d$ (sometimes also known as the $E_d$ operad). 
Equivalently, one says the observables are an $E_d$-algebra. 
This notion was introduced by algebraic topologists in the 1960s, 
so there is a rich body of work to draw upon.
\item For a holomorphic field theory on $\CC^d$,
the observables form an algebra over the $d$-polydiscs operad~$\PDisc_d$.
This notion was introduced only recently in~\cite{CG1}.
\end{itemize}
We will make the assertions above precise (in part by adding some natural but necessary hypotheses on a theory)
and develop them in parallel.
As a particularly illuminating example, we will explore the $d=1$ case rather closely and explain in what sense a vertex algebra is a kind of ``holomorphic associative algebra.''
We will also discuss the polydiscs algebras for which the higher Kac-Moody algebras are a Lie-algebraic shadow.

\begin{rmk}
In the $d=1$ case there is a lovely exploration of how the highly developed theory of vertex algebras relates to this perspective in the pair of papers by Bruegmann \cite{Bruegmann1,Bruegmann2}.
For a discussion of how conformal nets --- an alternative formalism for chiral CFT --- relates to our perspective, see~\cite{Henriques, BPS}.
\end{rmk}

Before jumping into the main story, we mention a motivation for talking about $E_d$- and $\PDisc_d$-algebras in parallel.
By dealing with both, we are suggesting, of course, that there is a natural analogy and hence a useful transfer of intuition between the settings of topological and holomorphic field theories.
But there is a stronger, technical relationship that we'll develop via a version of dimensional reduction.
From each holomorphic field theory in complex dimension $d$, 
we'll extract an $E_d$-algebra.
This construction generalizes the relationship between the Kac-Moody vertex algebra of a Lie algebra $\fg$ at level $k$ and the central extension of the loop algebra~$L\fg$.
We expect this correspondence to bear rich fruit.

There is an explicit model for an algebra over the cohomology of the $E_d$-operad called a $P_d$-algebra.
Essentially, these are graded commutative algebras with a Poisson bracket of cohomological degree $-d+1$.
In \cite{CosScheim} a parallel model for the cohomology of the holomorphic $\PDisc_d$ colored operad was proposed, 
though axioms were not written down.
In recent work, similar models to that one have appeared,
using the nomenclature of $\lambda$-brackets \cite{BGKWY,GKW},
following Kac's approach in the $d=1$ setting. 
The multivariable $\lambda$-brackets are similar in spirit to Lie conformal algebras,
which use the notion of single-variable $\lambda$-bracket to encode the ``singular part'' of the data of a vertex algebra.
These papers use this algebraic structure to give a beautiful proposal for computing Feynman diagrams in holomorphic (and mixed holomorphic-topological) field theories (see Section~\ref{sec: renorm} for more on renormalization within the context of holomorphic QFT).

\subsection{Generalizing associative algebras for topological field theories}

We want to describe the algebraic structures that appear in $d$-dimensional topological field theories,
as a model for holomorphic field theories.
It will be helpful to use the language of operads, which are an efficient and powerful tool for defining and studying algebraic structures,
so we will gloss the key ideas of operads briefly and then turn to describing the operads relevant to TFTs.
A more extensive discussion can be found in \cite{CG1}, with \cite{EllSaf} studying the topological setting in an elegant way.

An associative algebra is, by definition, a vector space $A$ with a linear map $\mu: A\otimes A \to A$ satisfying that the linear map $\mu \circ (\id \otimes \mu): A^{\otimes 3} \to A$ equals the linear map $\mu \circ (\mu \otimes \id)$.
(Let's ignore the other data and conditions for the moment.)
Nearly all flavors of algebras --- Lie, Poisson, and so on --- have a similar definition:
one has a collection of multilinear maps that satisfy some relations.
An operad is a way of specifying such data, so there is an associative operad $Ass$,
a Lie operad $Lie$, and so on.

Before sketching the definition, we give a concrete example.
Let $V$ be a vector space.
Then every $n$-fold multilinear operator on $V$ is an element of $\Hom(V^{\otimes n}, V)$;
that vector space parametrizes all possible $n$-ary operations on $V$.
One can also compose such operations to produce new ones.
For instance, given $\mu : V^{\otimes m} \to V$ and $\nu: V^{\otimes m} \to V$, we can put the output of $\mu$ into the $i$th input of $\nu$ to produce a new operation
\[
\begin{array}{cccc}
\nu \circ_i \mu: & V^{\otimes m+n-1} & \to & V\\
& v_1 \otimes \cdots \otimes v_{m+n-1} & \mapsto & \nu(v_1 \, \ldots, \mu(v_i, \ldots, v_m), v_{m+1}, \ldots, v_{m+n-1})
\end{array}
\]
where $i$ can run over any index between 1 and $n$.
We call this collection of data the {\em endomorphism operad of $V$}, denote $\End_V$.
It consists of 
\begin{itemize}
\item a list of vector spaces labelling $n$-ary operations, with 
\[
\End_V(n)= \Hom(V^{\otimes n}, V),
\]
with $n \in \NN$, and
\item a list of linear maps
\[
\circ_i : \End_V(n) \otimes \End_V(m) \to \End_V(n+m-1)
\]
for every pair $n$ and $m$ and where $1 \leq i \leq n$.
\end{itemize}
There are certain natural relations when one does multiple compositions, like $\mu \circ_i (\nu \circ_j \pi)$, arising from how multilinear maps compose.

Just as the model case of an associative algebra is $\End(V) = \Hom(V,V)$ (the ring of square matrices),
the endomorphism operad is the model case of operads.
A linear operad $P$ consists of a list of vectors $\left(P(n) \right)_{n \in \NN}$, where each $P(n)$ is a vector space describing $n$-ary operations, and a list of linear maps $\left(\circ_i\right)$.
We will not describe the relations, but resemble those of an endomorphism operad.

\begin{rmk}
The reader might ask why the adjective {\em linear} was appended before {\em operad}.
The reason is that we discussed algebraic operations on vector spaces,
but one can, of course, discuss analogs in other settings.
For instance, a monoid is a set $M$ with a multiplication map satisfying an associativity constraint.
One can talk about operads in sets (``an operad in sets $P$ has a set $P(n)$ of $n$-ary operations for every $n \in \NN$ \dots'') or in topological spaces or in dg vector spaces.
In fact, the theory works in any symmetric monoidal category. 
For the most part, we will work with vector spaces below.
\end{rmk}

Let us now raise another useful analogy.
To define a module for an associative algebra $A$, one specifies an algebra map $\rho: A \to \End(V)$.
This map encodes how $A$ acts on $V$.
Similarly, to give an algebra for an operad $P$, one specifies an operad map $\rho: P \to \End_V$.
This map picks out specific multilinear operations on $V$ that will satisfy the relations we want.
For instance, if $Pois$ is the Poisson operad, then $Pois(2) \cong \CC^2$ and this vector space is spanned by an element $\cdot$, the commutative product, and $\{\}$, the Poisson bracket.
Given a $P$-algebra $V$, the map $\rho(2): Pois(2) \to \End_V(2)$ thus tells us both the commutative and Poisson bracket put upon~$V$.
(The reader might enjoy trying to work out $Pois(3)$ and how to get elements there by composing binary operations.)
Alternatively, by applying the hom-tensor adjunction to $\rho(2)$, we have a map
\[
Pois(2) \otimes V^{\otimes 2} \to V
\]
that evaluates a binary operation with a pair of inputs.

By working at the level of operads, one can obtain universal results that hold for all algebras over a given operad.
For example, there is an operad map $b: Lie \to Pois$ that, at the level of the binary operations, includes the bracket, 
and this map induces the functor $b^*: \Alg_{Pois} \to \Alg_{Lie}$ where a Poisson algebra is ``forgotten'' to a Lie algebra.

So far we have only discussed putting algebraic structures on a single vector space $V$,
but it is possible to allow collections of vector spaces.
For instance, the pair of an algebra and a module provides an example of ``two-colored operad'' as the collection always consists of two vector spaces.
When one allows multiple objects, one works with {\em colored operads} or, as an alternative terminology, {\em multicategories}.

\begin{rmk}
For motivational introductions, we recommend \cite{StashAMS} and \cite{Val12}.
For systematic development of the theory, one might start with \cite{LodVal,Fresse,LurieHA}.
Operads now have an extensive literature, and there is a body of deep results (such as Koszul duality) with far-reaching consequences throughout mathematics.
\end{rmk}

With this language to hand, we now turn to generalizing associative algebras for the setting of topological field theories. 
Our goal is to formulate a colored operad that captures the behavior of the observables of a topological field theory.

The heuristic idea is straightforward.
Given a field theory $\cT$ on $\RR^d$, there is a vector space (or possibly dg vector space) $\Obs_\cT(D_r(x))$ of measurements one can make in the open disc $D_r(x)$ of radius $r$ around the point $x$.
Any measurement in a small disc should provide a measurement on a bigger disc,
so there is a map $\Obs_\cT(D_r(x)) \to \Obs_\cT(D_R(y))$ whenever $D_r(x) \subset D_R(y)$.
Thus we have a collection of vector spaces labeled by the discs in $\RR^d$,
and we have identified a class of unary operations determined by inclusion.

There should also be natural $n$-ary operations encoding how to combine observables from several distinct discs to give an observable on a big disc.
Start with $n$ discs $D_{r_1}(x_1)$, \dots , $D_{r_n}(x_n)$ whose closures are pairwise disjoint, and let $D_R(y)$ be a big disc containing all the smaller discs.
Then we want a map
\[
\Obs_\cT(D_{r_1}(x_1)) \otimes \cdots \otimes \Obs_\cT(D_{r_n}(x_n)) \to \Obs_\cT(D_R(y)).
\]
In the limit when the big disc is all of $\RR^d$ (i.e., $R = \infty$) and the input radii $r_i$ are very small, 
this map should encode the $n$-point functions beloved of physicists.

So far we have described a structure that observables should have for any field theory on~$\RR^d$.
For a topological field theory, we expect that the size of a disc does not matter:
for any inclusion $D_r(x) \subset D_R(y)$, the map $\Obs_\cT(D_r(x)) \to \Obs_\cT(D_R(y))$ is an equivalence.
This condition implies that we can always shrink the input radius as small as we want and we can always take the output radius to be a fixed value, say $R =1$.
In other words, the unary operations are parametrized by a point in the unit disc;
we have a local system of isomorphisms.

For $n$-ary operations we can likewise shrink the input radii as small as we want.
Thus the $n$-ary operations should be parametrized by the configuration space of $n$ points inside the unit disc.
(Rather, we can recover all $n$-ary operations knowing this information.)
We expect that these $n$-ary operations should also define a local system.

Our list of expectations is motivated by the examples of topological field theories that physicists typically give, rather than the Atiyah-Segal approach in terms of functors out of bordism categories.
Loosely speaking, a classical field theory (defined by a Lagrangian density, also known as an action functional) is ``topological'' if solutions to the equations of motion do not depend on the geometry of the spacetime manifold but only on its smooth structure (or other differential-topological features, like orientations).
A model example is Chern-Simons theory on oriented 3-manifolds for compact Lie group~$G$:
the space of solutions to the equations of motion on $M$ is the moduli space of flat $G$-bundles on $M$.
The observables for a classical field theory on $M$ means the commutative algebra of functions on solutions to the equations of motion on $M$.
In particular, for $M$ a disc, the size of the disc will not affect the algebra of observables,
and inclusions of discs leads to isomorphisms of observables.
Several of these features extend to the quantum field theories that arise by quantizing such topological classical theories;
the observables continue to form a vector space but are not commutative algebras anymore.
(See \cite{CG1, CG2} for a precise mathematical formulation of this idea, which is a natural extension of deformation quantization to higher dimensions.)

\begin{rmk}
There are examples of TFTs from physics where the solutions do not satisfy this condition,
such as the $A$-model of mirror symmetry.
In such cases, however, the observables are typically the de Rham complex (or some other model for cohomology) of the space of solutions,
and these observables do satisfy the conditions we've specified.
\end{rmk}

We are now in a position to offer a nice operad that should describe the observables of a topological field theory.
This operad was introduced by algebraic topologists for wholly distinct purposes (understanding $d$-fold based loop spaces),
so a remarkable dialogue has opened between topologists and physicists.
We will start, following the topologists, in the category $Top$ of topological spaces
where we use the cartesian product for the symmetric monoidal structure.
An operad in this category will be called a {\em topological operad}.

\begin{dfn}
The {\em little $d$-discs operad} $\Disc_d$ is the topological operad whose space of $n$-ary operations is
\begin{align*}
\Disc_d(n) = \{ (x_1, \ldots, x_n; r_1, &\ldots, r_n)  \in D_1(0)^n \times (0,1)^n \;|\; \\
&\text{the closure of the discs $\overline{D_{r_i}(x_i)}$ is pairwise disjoint}\}.
\end{align*}
The composition $\mu \circ_i \nu$ is given by rescaling the unit disc that labels the output of $\nu$ to $r_i$ and then embedding it into the disc~$D_{r_i}(x_i)$.
\end{dfn}

To illustrate the idea, we take $d=2$ and draw a point in the configuration space~$\Disc_2(2)$:
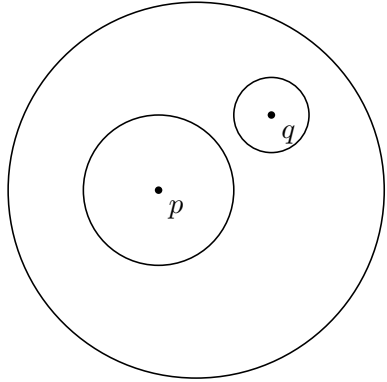
\begin{center}
  \begin{tikzpicture}
   \draw[semithick] (0,0) circle (2.5);
   \draw[semithick] (-0.5,0) circle (1);
   \draw[semithick] (1,1) circle (0.5);
   \fill (-0.5,0) circle (0.05) node [below right] {$p$};
   \fill (1,1) circle (0.05) node [below right] {$q$};
\end{tikzpicture}
\end{center}
where the outermost disc is the unit disc,
the input disc centered at the point~$p$ has larger radius, 
and the input disc centered at the point~$q$ has smaller radius.
By varying the locations of $p$ and $q$ and by varying the radii of the inner circles,
we move through~$\Disc_2(2)$.

This picture clearly resembles the ``pair of pants" product from two-dimensional topological field theories, but pressed flat.
To a physicist, this picture encodes a way of multiplying an operator with support near 0 by an operator with support near $p$ to get an operator with support near the origin (although the support is larger);
they use the term {\em operator product expansion} when they express this multiplication in formulas.

This operad was invented by algebraic topologists to capture the rich algebraic structure of the $d$-fold based loop space $\Omega^d X$ of a pointed space $(X,x)$.
(The usual homotopy group $\pi_d(X)$ is a shadow of this richer structure.)
Thus there is a deep reservoir of tools and results for studying algebras over the little discs operads in algebraic topology.

Observables, however, are supposed to live in vector spaces or dg vector spaces,
so we need to explain how to get algebras in that setting.
Thankfully, it is straightforward to get a linear version of this operad a dg operad (i.e., operad with values in dg vector spaces and using the tensor product as the symmetric monoidal product).
One simply takes the singular chain complex on each space of operations:
for each natural number $n$, take $C_\bu(\Disc_d(n); \CC)$ and take the induced chain map 
\[
C_\bu(\circ_i): C_\bu(\Disc_d(n); \CC) \otimes C_\bu(\Disc_d(m); \CC) \to C_\bu(\Disc_d(n+m-1); \CC)
\]
to produce a dg operad we denote~$C_\bu(\Disc_d)$.
(Note that as we always want to work with cochain complexes, we regrade the usual singular chains $C_\bu(X, \CC)$ to sit in nonpositive cohomological degrees. 
Likewise the composition $C_\bu(\circ_i)$ is then a cochain map.)

\begin{dfn}
An {\em algebra over the little $d$-discs operad} (or little $d$-discs algebra or $E_d$-algebra) is a cochain complex $A$ that is an algebra over the dg operad~$C_\bu(\Disc_d)$.
\end{dfn}

Examples appear very naturally from topology: for any space $X$ with a basepoint $x$, take the singular chains $C_\bu(\Omega^d X)$ of the $d$-fold based loop space $\Omega^d_x X$.
A much more subtle source of examples was conjectured by Deligne and subsequently proven in many ways: 
for any associative algebra $A$, the Hochschild cochain complex ${\rm Hoch}^\bu(A,A)$ can be equipped naturally with an $E_2$-algebra structure.

We can state now encapsulate our expectations from above:
\begin{quote}
the observables of a topological field theory on $\RR^d$ are a little $d$-discs algebra.
\end{quote}
There is a setting where a precise version of this statement was proven, by Elliott and Safronov, for a mathematically well-defined class of Lagrangian field theories.

\begin{thm}[\cite{EllSaf}]
\label{thm: EllSaf}
Let $\cT$ be a classical field theory on $\RR^d$
and let its observables $\Obs_\cT$ be constructed following~\cite{CG2}.

If the theory $\cT$ satisfies 
\begin{enumerate}
\item[(i)] it is equivariant under the action of the translation group,
\item[(ii)] the action of the translation Lie algebra is homotopically trivialized,\footnote{View the translation algebra $\RR^d$ as spanned by the basis $\{\partial_j\}$ of constant vector fields, and let $\rho: \RR^d \to \Der(\Obs_\cT)$ denote the translation action on observables.
By a homotopical trivialization, we mean here that there is a linear map $\eta: \RR^d \to \Der(\Obs_\cT)$ of degree $-1$ such that $[\d_{\cT}, \eta(\partial_j)] = \rho(\partial_j)$.
In other words, $\eta$ provides a chain homotopy trivialization of the action~$\rho$.} 
and
\item[(iii)] the extension of observables along an inclusion of discs is a quasi-isomorphism,
\end{enumerate}
then the observables form a little $d$-discs algebra.
\end{thm}

This theorem applies to many well-known examples from physics:
\begin{itemize}
\item The topological $B$-model of maps from an oriented surface into a Calabi-Yau manifold $X$ yields an $E_2$ algebra. It is quasi-isomorphic to the well-known Gerstenhaber algebra of polyvector fields on~$X$. (The quantization was constructed by Li and Li \cite{LiLi} in a setting where \cite{EllSaf} applies.)
\item The Chern-Simons gauge theory of flat $G$-bundles on an oriented 3-manifold yields an $E_3$ algebra, which encodes the quantum group $U_\hbar \fg$ under a version of Koszul duality. 
(The quantization was constructed, essentially, by Axelrod-Singer \cite{AxeSing1,AxeSing2} and Kontsevich \cite{KonECM} but articulated in this setting by Costello~\cite{CosBook}.) 
\item The Rozansky-Witten theory of maps from an oriented 3-manifold into a complex manifold yields an $E_3$ algebra. (The quantization was constructed by Chan, Leung, and Li~\cite{ChanLeungLi}.)
\item The Kapustin-Witten gauge theory for oriented 4-manifolds yields an $E_4$ algebra, related to deforming functions on the coadjoint quotient stack~$\fg^*/G$. (The quantization was constructed in~\cite{EGW}.)
\end{itemize}
It applies to other examples of topological AKSZ theories, such as BF theories and the Poisson $\sigma$-model.
In all these examples, the underlying elliptic complex appearing in the field theory is a version of the de Rham complex,
so that the chain homotopy $\eta$ is easily constructed by contracting a vector field with differential forms.

This theorem builds an explicit bridge between topology and physics. 
Thanks to the theorem of Elliott and Safronov, we know that most topological field theories in the physicist's sense provide $E_n$ algebras (with the examples above as evidence),
but they can thus be analyzed using the powerful machinery of modern homotopy theory. 
Conversely, intuition from physics about the behavior of topological field theories suggests novel constructions and examples to topology. 
Much of the amazing resonance between operads, homological algebra, and string-theoretic physics over the last few decades can be understood from this perspective.
(See Section~\ref{sec: reln to functorial TFT} for the relation with the Atiyah-Segal approach to TFT.)

Let us comment on some aspects of the theorem and its proof.
The notion from \cite{CosBook} articulates Lagrangian field theories whose linearized equations of motion are elliptic in nature;
it uses (and mathematically articulates) the Batalin-Vilkovisky formalism,
which is inherently homological.
In particular, there is an underlying graded vector space of ``fields,''
and the observables are a dg commutative algebra given a symmetric algebra on the linear dual to the fields and equipped with a differential $d_\cT$ that is determined by the Lagrangian density of the theory.
Item (ii) above then means that each vector field $\partial_j$ acts by a degree 0 cochain map $D_j$ on the observables and there is a degree one cochain homotopy $\eta$ on the observables such that $[d_\cT, \eta(\partial_j)]~=~D_j$.

For physicists, we note that the proof uses a version of Witten descent.
Indeed, descent can be used to produce interesting operations in these disc algebras,
so that many familiar constructions from the physical approach to TFTs fit cleanly and easily into this framework.

\subsection{Generalizing algebras for holomorphic field theories}

We now turn to the problem of describing observables of a holomorphic field theory on~$\CC^d$,
and we follow an approach similar to the topological case just discussed.

In the setting of a holomorphic field theory, we expect the observables to have a nontrivial dependence on the radius of a disc.
Consider, for example, that 
when $D_r(x) \subset D_R(y) \subset \CC$ are open discs in the plane,
the restriction map of holomorphic functions
\[
\cO(D_R(y))\to \cO(D_r(x))
\]
is not an isomorphism,
and the sheaf $\cO$ provides a fundamental model for how the solutions of the equations of motion should behave for a holomorphic field theory.
We will assume, however, that a holomorphic theory is equivariant under translations of $\CC^d$ and that the observables in a disc $D_r(X)$ are equivalent to the observables in a translate~$D_r(x+z)$.
Thus, we will have a colored operad whose colors are labeled by radii~$r \in (0,\infty)$,
as we merely have to specify the observables on a disc~$D_r(0)$.

We now introduce a small tweak.
Instead of working with discs, we work with polydiscs: for $w = (w_1, \ldots, w_d) \in \CC^d$, the polydisc of radius $r$ around $w$~is
\[
PD_r(w) = \{ z \in \CC^d \;|\; |z_i - w_i| < r \text{ for all $i$}\}.
\]
We work with polydiscs because it is easier to describe holomorphic functions on a polydisc than on a disc;
one can borrow results from single-variable complex analysis.
With this notion in hand, we introduce our main character.

\begin{dfn}
The {$d$-polydiscs operad} $\PDisc_d$ is the colored operad in the category of complex manifolds where each positive real number $r \in (0,\infty)$ is a color 
and where there is a complex manifold of $n$-ary operations
\begin{align*}
\PDisc_d(r_1, \ldots, r_n| R) = \{ (w_1, \ldots,& w_n) \in PD_R(0)^n \;|\; 
\text{the closure of the polydiscs $\overline{PD_{r_i}(w_i)}$}\\ 
&\text{ is pairwise disjoint and each is contained in $PD_R(0)$}\}
\end{align*}
for any list of radii $r_1,\ldots, r_n, R$.
This manifold is always an open subset of~$PD_R(0)^n$. 
(Note that the manifold is empty for many values.)
There is a composition
\[
\circ_i: \PDisc_d(r_1, \ldots, r_n| R) \times \PDisc_d(s_1, \ldots, s_m| S) \to \PDisc_d(r_1, \ldots, s_1, \ldots, s_m, r_{i+1}, \ldots, r_n| R)
\]
only when $r_i = S$, 
and in that case, the composition $\mu \circ_i \nu$ is given by embedding an output polydisc into an input polydisc of the same radius.
(Composition is not defined if the radii do not match.)
\end{dfn}

We are again thinking about configuration spaces --- here of polydiscs inside larger polydiscs --- but we are remembering the complex structure explicitly.
To describe observables, we want to produce a dg colored operad that encodes the holomorphic nature of this situation.
We do {\em not} want to take singular chains, which would only remember the homotopy type;
instead, we choose to take the Dolbeault complex of each configuration space.
Note, however, that taking the Dolbeault complex is a contravariant functor,
so we get a {\em co}\/operad rather than an operad.
This notion bears the same relationship to operad as a coalgebra does to an algebra;
one simply reverses all the arrows and asks for an appropriate variant of the relations.
For instance, there is a composition
\begin{multline}
\label{pdcirc}
\circ_i: \Omega^{0,\bu}(\PDisc_d(r_1, \ldots, s_1, \ldots, s_m, r_{i+1}, \ldots, r_n| R)) \to \\
\Omega^{0,\bu}(\PDisc_d(r_1, \ldots, r_n| R)) \otimes \Omega^{0,\bu}(\PDisc_d(s_1, \ldots, s_m| S))
\end{multline}
when $r_i = S$, and in that case, the composition $\mu \circ_i \nu$ is given by embedding an output polydisc into an input polydisc of the same radius.
(Composition is not defined if the radii do not match.)

An important technical point arises here: 
we need to be careful about functional analysis.
The wedge product map
\[
\Omega^{0,\bu}(X) \otimes \Omega^{0,\bu}(Y) \to \Omega^{0,\bu}(X \times Y)
\]
is not a quasi-isomorphism in general (unlike with de Rham complexes),
but the completed projective tensor product as topological vector spaces 
\[
\Omega^{0,\bu}(X) \,\widehat{\otimes}_\pi\, \Omega^{0,\bu}(Y) \to \Omega^{0,\bu}(X \times Y)
\]
is an isomorphism.
(We use here the standard Fr\'echet topology on smooth functions, extended to sections of vector bundles.)
Hence we should work in some category of cochain complexes of well-behaved topological vector spaces, such as convenient vector spaces, or something related.
This requirement makes sense when geometry matters (e.g., using the Dolbeault complex) and not just topology.
For an extensive discussion of these issues, see the appendices of~\cite{CG1}.
From hereon we will ignore function-analytic aspects and write~$\otimes$,
with an appropriate choice of symmetric monoidal category left implicit.

\begin{dfn}
\label{dfn pdisc}
The {\em $d$-polydiscs dg cooperad} is the colored operad in the category of dg vector spaces where each positive real number $r \in (0,\infty)$ is a color 
and where there is a dg vector space of $n$-ary cooperations
\[
\Omega^{0,\bu}(\PDisc_d(r_1, \ldots, r_n| R)) 
\]
for any list of radii $r_1,\ldots, r_n, R$.
There is a composition
\begin{multline}
\circ_i \colon \Omega^{0,\bu}(\PDisc_d(r_1, \ldots, s_1, \ldots, s_m, r_{i+1}, \ldots, r_n| R)) \to \\ \Omega^{0,\bu}(\PDisc_d(r_1, \ldots, r_n| R)) \otimes \Omega^{0,\bu}(\PDisc_d(s_1, \ldots, s_m| S))
\end{multline}
only when $r_i = S$, 
and in that case, the composition $\mu \circ_i \nu$ is given by embedding an output polydisc into an input polydisc of the same radius.
(Composition is not defined if the radii do not match.)
\end{dfn}

An algebra $A$ over this cooperad means that for each radius~$r$, 
there is a dg vector space $A(r)$, and for any list of radii $r_1,\ldots, r_n, R$,
there is a cochain map
\[
\alpha(r_1, \ldots, r_n| R): A(r_1) \otimes \cdots \otimes A(r_n) \to A(R) \otimes \Omega^{0,\bu}(\PDisc_d(r_1, \ldots, r_n| R)) .
\]
To invest the map $\alpha$ with more meaning, 
pick a real, closed submanifold 
\[
T \subset \PDisc_d(r_1, \ldots, r_n| R)
\]
along which one can integrate $(0,k)$-forms.
For instance, for any inputs $a_i \in A(r_i)$, we get an element in $A(R)$~by
\[
\int_T \alpha(a_1 \otimes \cdots \otimes a_n),
\]
where we have suppressed distracting notation from~$\alpha$.

We will now examine the case $d=1$ in some detail and relate such algebras to vertex algebras.
Afterward, we explain how holomorphic field theories produce such algebras,
by a result analogous to Elliott and Safronov's.
In particular, a holomorphic version of Witten descent gives a useful method for producing interesting observables for holomorphic field theories.

\subsubsection{Vertex algebras from the $d=1$ case}

In the case $d=1$, it is possible to relate polydiscs algebras to vertex algebras.
More precisely, there is a functor from a class of 1-polydiscs algebras to vertex algebras,
developed in Chapter 5 of~\cite{CG1}.
In essence, a polydisc algebra is in this class if it has a nice action of the rotation group~$U(1)$,
so that one can extract ``modes'' (or Fourier components) from it.
Examples like the free $\beta\gamma$ system and the Kac-Moody vertex algebras are analyzed in~\cite{CG1}.

Here we will analyze the cohomology of the cooperad and what kind of information it encodes, 
aimed at relating it to vertex algebras.
The takeaway from the discussion below is that one can read off the modes $a_{(n)}$ that assemble into the vertex operator $Y(a,z)$ using a $1$-polydiscs algebra.

It is easiest to work with extremal values of the radii.
We take the output radius to be $R = \infty$,
so that the outputs are observables on all of~$\CC$.
We take the input radii to be $r_i = 0$,
so that the inputs are observables supported at the points $w_i$.
(A physicist might call these ``local operators.'')
This input value is not allowed in the definition, but it arises as a natural limit.
In any case, it is natural to think about the complex manifold ${\rm Conf}_n(\CC)$ of configurations of $n$ distinct points in~$\CC$:
it is
\[
{\rm Conf}_n(\CC) = \CC^n \setminus \{\text{union of the diagonals $\{ w_i = w_j\}$}\}.
\]
Note that for any finite and nonzero values of radii, there is an inclusion
\[
\PDisc_d(r_1, \ldots, r_n| R) \hookrightarrow {\rm Conf}_n(\CC).
\]
We now turn to examining the Dolbeault cohomology.

It is easy to see
\[
{\rm Conf}_n(\CC) \cong \CC \times (\CC^\times)^{n-1}
\]
by sending a configuration $(w_1,\ldots, w_n)$ to $(w_1, w_2 - w_1, \ldots, w_n-w_{n-1})$.
(We could replace $w_1$ by any $w_j$, or even the sum of them all.)
Hence, the Dolbeault cohomology consists of holomorphic functions on this space;
there is no higher cohomology.
Inside the holomorphic functions 
\[
\cO({\rm Conf}_n(\CC)) = H^{0,\bu}({\rm Conf}_n(\CC)),
\] 
as a dense but understandable subalgebra, are the algebraic functions:
\[
\cO_{\rm alg}({\rm Conf}_n(\CC)) = \CC[w_1, (w_1 - w_2)^{\pm 1}, \ldots, (w_n - w_{n-1})^{\pm 1}].
\]
Here $(w_j - w_{j-1})^{\pm 1}$ means that the function $1/(w_j - w_{j-1})$ is a generator, 
as is $(w_j - w_{j-1})$.
Each such algebraic function provides an operation on $H^*(A)$, where $A$ denotes a polydisc algebra.

Let us unravel things in the special case of two points.
For simplicity, we will fix $w_1 = 0$ so that we focus on where place the other point $w_2 \in \CC^\times$.
We will use $z$ to denote this point $w_2$.
Then we have binary co-operations parametrized by the space $\C^\times$.
Let 
\[
V = \lim_{r \to 0} H^*(A(D_r(0)))
\]
denote the cohomological observables ``supported at the origin.''
(It would be better to take the homotopy limit of the $A(D_r(0))$,
but we leave the experienced reader to make such adjustments.)
For any two elements $a, b \in V$ and for any $n \in \ZZ$,
we then obtain an element
\[
a_{(n)} b = \int_{|z| = r} \alpha(a,b)\, z^{-n-1} \d z
\]
in $H^*(A(D_R(0))$ for any radius $R > r$.
(The choice to match $n$ with $z^{-n-1}$ is conventional.)
Any holomorphic function $f \in D_R(0) \setminus \{0\}$ also defines a map
\[
\int_{|z| = r} \alpha(a,b)\, f(z) \d z 
\]
and, if $f$ is meromorphic at the origin, then its Laurent expansion
\[
f(z) = c_{-N} z^{-N} + \cdots + c_0 + c_1 z + \cdots
\]
guarantees that we have
\[
\int_{|z| = r} \alpha(a,b)\, f(z) \d z = \sum_{n=-\infty}^{N-1} c_{-n-1} a_{(n)} b.
\]
Thus we have worked out the modes $a_{(n)}$, as promised.

This construction also gives insight into some challenges with the formalism of vertex algebras,
such as the fact that one cannot simply compose vertex operators but must describe a 3-point product and track the geometry between the insertion points $0$, $z$, and $w$.
For instance, to obtain the modes, we took the outgoing radius to be infinite,
so it cannot be inserted into any finite radius disk.
Nonetheless, the associativity of composition in $\PDisc_1$ casts the shadow of the Cousin property or Borcherds product in vertex algebras.

\subsubsection{Examples related to vertex algebras} 

We have just recalled that in the complex dimension one case, holomorphic  polydiscs algebras (with a compatible $U(1)$ action) correspond to vertex algebras.
Let us briefly point out a few examples of this correspondence.

Given a Lie algebra $\lie{g}$ we obtain a one-dimensional holomorphic disc algebra $A_{\lie{g}}$ by the following construction. 
To a disk $D_r(z)$ of radius $r$ centered at $z \in \CC$ we attach the cochain complex
\beqn
\label{eqn:diskskm}
A^{\lie{g}} (D_r(z)) = C_\bu \left(\lie{g} \otimes \Omega^{0,\bu}_c(D_r(z)) \right) 
\eeqn
Here, $\Omega^{0,\bu}_c(D_r(z))$ denotes the Dolbeault complex of compactly supported $(0,\bu)$ forms on the disc.
Because the space of $(0,\bu)$ forms define a commutative dg algebra, 
the complex $\lie{g} \otimes \Omega^{0,\bu}_c(D_r(z))$ has the natural structure of a dg Lie algebra---the complex on the right hand side is the Chevalley--Eilenberg complex that computes the Lie algebra homology of this dg Lie algebra.
Using the formalism of factorization algebras (which we recall in the next section) it is shown in \cite{CG2} how to equip $A_{\lie{g}}$ with the structure of a holomorphic disc algebra.
The corresponding vertex algebra is the Kac--Moody vertex algebra of $\fg$ at level zero. 
To obtain the vertex algebra associated to a nonzero level $\kappa$ one can deform this holomorphic polydisc algebra $A^{\fg} \rightsquigarrow A^\fg_{\kappa}$.

Let $V$ be a vector space.
Then, we can define a holomorphic disc algebra $A_V$ that corresponds to another familiar vertex algebra: the $\beta\gamma$ system with values in~$V$.
To a disc $D = D_r(z)$ this assigns the following cochain complex
\beqn
\Sym \left(\Omega^{1,\bu}_c(D) \otimes V^* [1] \oplus \Omega^{0,\bu}_c (D) \otimes V [1]\right) 
\eeqn
where the differential is $\dbar + \triangle$ that we will explain momentarily. 
We denote the linear generators of this symmetric algebra by $\gamma \in \Omega^{1,\bu}_c(D) \otimes V^* [1]$ and $\beta \in \Omega^{0,\bu}_c(D) \otimes V^* [1]$.
Now, for the differential, $\dbar$ is the familiar operator (which acts in the standard way on $(0,\bu)$ and $(1,\bu)$ forms) and $\triangle$ can be defined on $\Sym^2$ summand by the formula
\beqn
\triangle (\gamma \cdot \beta) = \int_D \<\gamma, \beta\> . 
\eeqn
The operator can be extended to the entire complex $A_V(D)$ by the rule that it is a BV operator for the $1$-shifted Poisson bracket given by wedge and integration.
The vertex algebra associated to this disk algebra is the familiar $\beta\gamma$-vertex algebra.

\subsubsection{Revisiting the higher Kac-Moody algebras}

It is apparent that the formula in \eqref{eqn:diskskm} makes sense for a (poly)disk of arbitrary complex dimension $d$.
In this way, we obtain functor from Lie algebras to $d$-polydisks algebras $\fg \mapsto A^{\fg}_d$.
(Here $A^{\fg}_1 = A^{\fg}$.)

In this section we have mostly focused on disks.
If we evaluate $A^{\lie{g}}_{d}$ on a punctured disk rather than a disk then we obtain a complex that receives a dense embedding from the complex
\beqn
\Sym \left(\lie{g}^\bu_{d}\right) 
\eeqn
where $\lie{g}^\bu_d$ is the higher sphere algebra from Definition \ref{dfn:sphere}.
A natural question is how the Lie algebra structure appears; to address this we will use the language of factorization algebras in the next section.

Given an invariant polynomial $\theta \in \Sym^{d+1}(\lie{g}^*)^G$ as in \S \ref{s:alghigh} we can obtain a twisted version of this $d$-polydisks algebra $A^{\lie{g}}_{d,\theta}$.
At the level of the punctured disk this corresponds to the higher Kac--Moody algebra that is the central extension of the sphere algebra $\lie{\fg}^\bu_{d}$.
This is a $d$-polydisks avatar of the higher Kac--Moody algebra, in the same way that the Kac--Moody vertex algebra is related to the affine Kac--Moody Lie algebra~\cite{GWonKM}.

\subsubsection{Polydiscs algebras from holomorphic field theories}

Let 
\[
\ft_\CC = \CC \otimes_\RR \ft = \text{span}_{\CC}\{\partial/\partial z_1, \ldots, \partial/\partial z_d, \partial/\partial \zbar_1, \ldots, \partial/\partial \zbar_d  \}
\]
denote the complexification of the Lie algebra~$\ft$ of translations.
Let $\ft_{\zbar}$ denote the subalgebra spanned by all the~$\partial/\partial \zbar_i$.

\begin{thm}[see Chapter 5 of \cite{CG1}]
Let $\cT$ be a classical field theory on $\CC^d$,
and let its observables be constructed following~\cite{CG2}.
If the theory satisfies 
\begin{enumerate}
\item[(i)] $\cT$ is equivariant under the action of the translation group and
\item[(ii)] the action of the Lie algebra $\ft_{\zbar}$ is homotopically trivialized,
\end{enumerate}
then the observables form a $d$-polydiscs algebra.
\end{thm}

This result is analogous to Elliott and Safronov's, although it predated (and inspired) it.
It applies to the theories that are the focus of this survey.
One can mimic the version of Witten descent for the topological case to produce observables from local observables.

A number of examples from Section~\ref{sec: ex of theories} have been developed in this framework, 
analogous to the list of examples after Theorem~\ref{thm: EllSaf}:
\begin{itemize}
\item The curved $\beta\gamma$ system for a Riemann surface into a complex manifold~$X$ was constructed in \cite{CosWG2}, and then \cite{GGW} demonstrates in detail how the vertex algebra of {\it chiral differential operators} of $X$ is recovered.
\item In \cite{Wthesis} these methods were extended to holomorphic $\sigma$-models from $\CC^d$ to a complex manifold $X$, 
yielding $d$-polydics algebras of chiral differential operators on~$X$.
\item Costello and Li quantized BCOV theory, yielding a 3-polydiscs algebra for holomorphic gravity on the Calabi-Yau 3-fold~$\CC^3$~\cite{CLbcovOC}.
\item The holomorphic twists of 4-dimensional supersymmetric Yang-Mills theories are constructed and analyzed in~\cite{EGW}.
\end{itemize}
A number of examples of topological-holomorphic theories have also been constructed, and these yield polydiscs algebras.
A highlight is Costello's discovery of a theory on $\CC \times \RR^2$ that encodes the Yangian, an infinite-dimensional quantum group~\cite{CosYangian}.

\section{Factorization algebras and holomorphic field theories}
\label{sec: fact alg}

A key idea for us is that much of the content in field theories --- both classical and quantum --- is captured by factorization algebras.
We will introduce that concept now, in the style that sheaves are defined, so first we describe {\em pre}\/factorization algebras and then impose a local-to-global condition to characterize factorization algebras.
Their local structure encompasses the behavior we saw in the preceding section, where $E_d$ algebras captured topological field theories on $\RR^d$ and $\PDisc_d$ algebras captured holomorphic field theories on~$\CC^d$.
Factorization algebras allow one to work on smooth or complex manifolds and hence bring in global aspects of geometry.
Their global behavior encompasses constructions like Hochschild homology (for $E_1$ algebras extended as factorization algebras on $S^1$) and conformal blocks (for $\PDisc_1$ algebras extended as factorization algebras to Riemann surfaces).
We will then sketch the main theorem of \cite{CG2}, 
which explains a precise and general relationship between quantum field theories and factorization algebras.
Finally, we describe some special features and interesting applications of factorization algebras of holomorphic field theories.

Kapranov gave a lecture series explaining how factorization algebras fit into algebraic geometry \cite{KapLect}.
A recent survey of factorization algebras, with a focus on their role in mathematical physics, can be found in~\cite{CosGwEMP}.

\begin{rmk}
Beilinson and Drinfeld introduced the notion of factorization algebras in their work on chiral conformal field theory \cite{BD},
and their version has had great success in algebraic geometry and representation theory,
particularly in the setting of the geometric Langlands program.
Francis, Gaitsgory, and Lurie ported the spirit of this notion into the setting of manifolds,
developing topological chiral homology \cite{LurieHA} and factorization homology \cite{AF}
for questions in algebraic topology.
Motivated by these efforts (and benefiting from conversations as those works were in progress), 
Kevin Costello and the first author developed a version that works well for a broad class of quantum field theories and for more differential-geometric settings~\cite{CG1,CG2}.
These approaches share a common spirit but have important technical differences.
A detailed comparison between the versions for locally constant factorization algebras exists~\cite{KSW},
but a detailed comparison between the Beilinson-Drinfeld and Costello-Gwilliam version is currently open.
(See, however, the last section of \cite{HenKap} for important steps in that direction.)
\end{rmk}

\subsection{The essential idea of a prefactorization algebra}

Let $M$ be a topological space and let $\mc{C}^\otimes$ be a symmetric monoidal category. 
In this paper $M$ is always a smooth manifold (typically a complex manifold) and $\mc{C}$ is ${\rm Vect}$ or ${\rm dgVect}$, with the usual tensor product as the symmetric monoidal product.
(As mentioned in the paragraph before Definition~\ref{dfn pdisc}, one needs to deal with issues of functional analysis for applications to QFT, but we will not discuss those aspects here.
In brief, for purposes of holomorphic field theories, we use an $\infty$-category obtained from working with cochain complexes of differentiable vector spaces, with extensive treatment in~\cite{CG1}.)

\begin{dfn}
A {\em prefactorization algebra} $\cF$ on $M$ taking values in cochain complexes is a rule that assigns a cochain complex $\cF(U)$ to each open set $U \subset M$ along with the following maps and compatibilities.
\begin{itemize}
\item  There is a cochain map $m_V^U: \cF(U) \rightarrow \cF(V)$ for each inclusion $U \subset V$.

\item There is a cochain map $m_V^{U_1,\ldots,U_n} : \cF(U_1) \otimes \cdots \otimes \cF(U_n) \rightarrow \cF(V)$ for every finite collection of open sets where each $U_i \subset V$ and the $U_i$ are pairwise disjoint. The following picture represents the situation.
\begin{center}
 \begin{minipage}[c]{3cm}
 \begin{tikzpicture}[scale=0.25]
 \draw (0,0) circle (5);
 \draw (-1.5,2) circle(1.3) node {$U_1$};
 \draw (-2.2,-1.5) circle (1.5) node {$U_2$};
 \draw (0, -3) node {\tiny$\dots$};
 \draw (2.1,-1) circle (1.8) node {$U_n$};
 \draw (2.1, 3) node {$V$};
 \end{tikzpicture}
 \end{minipage}
\hspace{0.7cm} $\rightsquigarrow$ \hspace{0.5cm}
 \begin{minipage}[c]{8cm}
$\cF(U_1)\otimes\dots\otimes\cF(U_n)\xto{m_V^{U_1,\ldots,U_n}}\cF(V),$
 \end{minipage}
\end{center}

\item The maps are compatible in the obvious way, so that if $U_{i,1}\sqcup\cdots\sqcup U_{i,n_i}\subseteq V_i$ and $V_1\sqcup\cdots\sqcup V_k\subseteq W$, the following diagram commutes.
\begin{center}
\begin{tikzcd}[column sep=small]
{\bigotimes}^{k}_{i=1}{\bigotimes}^{n_i}_{j=1}\cF(U_j) \arrow{dr} \arrow{rr} &   &{\bigotimes}^k_{i=1}\cF(V_i) \arrow{dl}\\
&\cF(W)  &
\end{tikzcd}
\end{center}
\end{itemize}
\end{dfn}

For an explicit example of the associativity, consider the following picture.
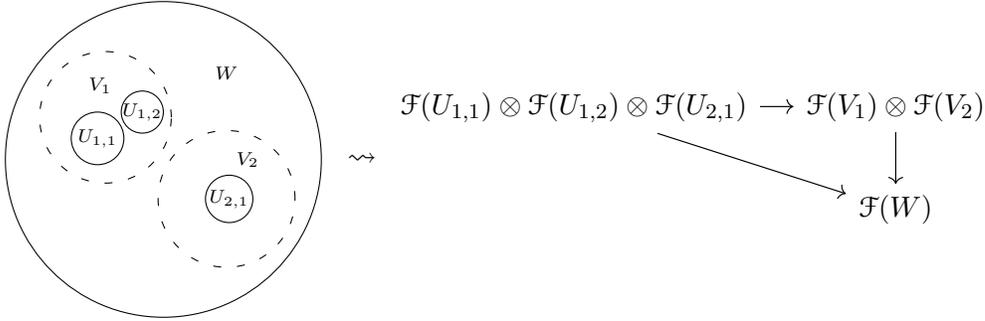
\begin{center}
\hspace{-1cm}
\begin{minipage}{3cm}
{\tiny
\begin{center}
\begin{tikzpicture}[scale=0.35]
\draw (0,0) circle (6);
\draw (2.4, 3.3) node {$W$};
\draw [style=loosely dashed] (-2.2,1.6) circle (2.5);
\draw (-2.4, 2.8) node {$V_1$};
\draw [style=loosely dashed] (2.4,-1.5) circle (2.6);
\draw (3.2, 0) node {$V_2$};
\draw (-2.5, 0.8) circle (1) node {$U_{1,1}$};
\draw (-0.8, 1.8) circle (0.8) node {$U_{1,2}$};
\draw (2.5,-1.5) circle (0.9) node {$U_{2,1}$};
\end{tikzpicture}
\end{center}
}
 \end{minipage}
 \hspace{1.3cm} $\rightsquigarrow$ \hspace{-0.1cm}
\begin{minipage}{8cm}
\begin{tikzcd}[column sep=small]
\cF(U_{1,1}) \otimes \cF(U_{1,2}) \otimes \cF(U_{2,1}) \arrow{dr} \arrow{r}  & \cF(V_1) \otimes \cF(V_2) \arrow{d}\\
&\cF(W) 
\end{tikzcd}
\end{minipage}

The case of $k=n_1=2$, $n_2 = 1$.
\end{center}

This definition bears analogies to familiar objects in mathematics.
On the one hand, $\cF$ resembles a {\em precosheaf}, which is a functor from opens in $M$ to a category like ${\rm Vect}$ or ${\rm dgVect}$.
(A presheaf is a functor out of the opposite category to opens in~$M$.) 
Here, however, $\cF$ also assign values to disjoint unions of opens, and it uses the tensor product rather than direct sum.
This feature leads to the other analogy: $\cF$ resembles an algebra, as
the multilinear maps look like multiplications.
These maps let us multiply elements from disjoint regions to get an element in a larger region.

\begin{figure}
\begin{center}
 \begin{tikzpicture} 
 \begin{scope}[line cap=round,ultra thick]
 \draw (-2.5,1) -- (-1.5,1); 
 \draw (-1,1) -- (-0.5,1); 
 \draw (1,1) -- (2,1); 
 \draw[->,semithick] (0,0.75) -- (0,0.25);
 \draw (-2.75,-0.1) -- (-0.25,-0.1); 
 \draw (0.75,-0.1) -- (2.5,-0.1);  
 \draw[->,semithick] (0,-0.35) -- (0,-0.85);
 \draw (-3,-1.2) -- (3,-1.2); 
 \end{scope}
 
\node at (4,0) {$\rightsquigarrow$ };

\node at (7,0){
\begin{tikzcd}[cramped,column sep=tiny]
a \otimes b \otimes c \arrow[mapsto]{d} & \in & A \otimes A \otimes A \arrow{d} \\
ab \otimes c \arrow[mapsto]{d} & \in & A \otimes A \arrow{d} \\
abc  & \in &A 
\end{tikzcd}
};
\end{tikzpicture}
\end{center}
\caption{The prefactorization algebra $\cF_A$ of an associative algebra $A$}
\label{fig:assasfact}
\end{figure}
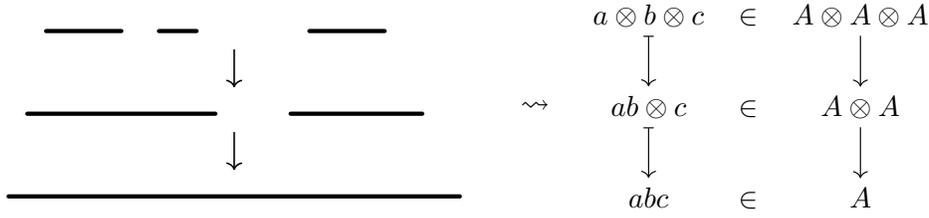

\begin{eg}\label{ex:associativealgebra}
Every associative algebra $A$ defines a prefactorization algebra $\cF_A$ on $\RR$, as follows. 
To each open interval $(a,b)$, we set $\cF_A( (a,b) ) = A$. 
To any open set $U = \coprod_j I_j$, 
where each $I_j$ is an open interval, we set $\cF(U) = \bigotimes_j A$. 
The structure maps simply arise from the multiplication map for $A$. 
Figure ~\ref{fig:assasfact} displays the structure of $\cF_A$. 
Notice the resemblance to the notion of an $E_1$ or $A_\infty$ algebra.~\hfill$\Diamond$
\end{eg}


\begin{eg}\label{ex:SymOfCosheaf}
Another important example for us is the symmetric algebra of a precosheaf. 
Let $F$ be a precosheaf of vector spaces on a space $X$. 
For example, consider $F = C^\infty_c$ the compactly supported smooth functions on a manifold. The functor $\cF  = \Sym F: U \mapsto \Sym(F(U))$ defines a precosheaf of commutative algebras, but it also a prefactorization algebra. For instance, if $U$ and $V$ are disjoint opens, we see that
\begin{align*}
\cF(U \sqcup V) = \Sym F(U \sqcup V) & \cong \Sym (F(U) \oplus F(V)) \\ & \cong \Sym(F(U)) \otimes \Sym(F(V)) = \cF(U) \otimes \cF(V),
\end{align*}
and these isomorphisms provide the structure maps for $\cF$.~\hfill$\Diamond$
\end{eg}

There is another, more sophisticated way to phrase prefactorization algebras, using the framework of operads.
Associated to each topological space $M$, there is a colored operad $\Disj_M$ such that a prefactorization algebra is an algebra over~$\Disj_M$.
For a discussion, see Chapter 3 of~\cite{CG1};
related ideas are developed in~\cite{AF, LurieHA, Ginot}.

\begin{rmk}
The reader fond of higher abstract nonsense will note that we could have prefactorization algebras take values in other target categories (or even better, $\infty$-categories),
not just vector spaces or cochain complexes.
For instance, one might categorify the values: 
$\cF$ might assign an $\infty$-category $\cF(U)$ to each open set~$U$.
This generalization is fruitful.
In the Beilinson-Drinfeld style, it leads to chiral categories, which play an important role in the geometric Langlands program.
In the style advocated above, it would likewise be useful for capturing notions like defects or extended operators from quantum field theory.
\end{rmk}

\subsection{Meet factorization algebras}

We are now in a position to define factorization algebras,
which are prefactorization algebras whose behavior on large open sets is determined by their behavior on small open sets.
To give a precise description, we introduce a Grothendieck topology due to Michael Weiss, 
which is explained nicely and further developed in~\cite{BoavidaWeiss}.

\begin{dfn}
Let $U$ be an open set. A collection of open sets $\mathfrak{U} = \{ U_i \mid i \in I\}$ is a {\em Weiss cover} of $U$ if for any finite collection of points $\{x_1,\ldots,x_k\}$ in $U$, there is an open set $U_i \in \mathfrak{U}$ such that $\{x_1,\ldots,x_k\} \subset U_i$.
\end{dfn}

The Weiss covers define a Grothendieck topology on $\Open(M)$, the poset category of open subsets of a space $M$. 
We call it the {\em Weiss topology} of~$M$. 
Note that Weiss covers are required to ``know'' about every configuration of finitely many points,
and so they often contain many opens.
For a smooth $n$-manifold $M$, a useful Weiss cover is the collection of open sets in $M$ diffeomorphic to a disjoint union of finitely many copies of the open $n$-disc.
(If one has a metric on the manifold, one can work with unions of ``small'' discs.)

Now that we have a notion of cover, we can formulate the local-to-global property for factorization algebras.

\begin{dfn}
A \emph{factorization algebra} is a prefactorization algebra $\cF$ on $M$ with values in a symmetric monoidal category $\cC^\otimes$ that is a cosheaf for the Weiss topology.
\end{dfn}

Being a cosheaf means that for any open set $U \subset M$ and any Weiss cover $\{ U_\alpha\}$ of $U$,
the value of $\cF(U)$ can be reconstructed using the \v{C}ech complex:
\[
\cF(U) \xto{\simeq} \check{C}(\mathfrak{U},\cF).
\]
Replacing cochain complexes with some other category (or $\infty$-category), one needs to take a (homotopy) colimit over the \v{C}ech nerve of the Weiss cover.

As verification that this notion has some real content, let's consider a first, interesting case.
We already explained how to produce a locally constant prefactorization algebra $\cF_A$ on $\RR$ from a dg associative algebra $A$.
We could ask about working on the circle $S^1$ instead of $\RR$.
This construction tells us what to assign to every proper open subset of $S^1$: for instance, to each open interval, we assign $A$.
The \v{C}ech complex should then tell us what to assign to the whole circle.
Remarkably, this construction recovers a well-known invariant: the Hochschild chains ${\rm Hoch}_*(A)$ of~$A$.

\begin{thm}[\cite{LurieHA, AF,KSW}]
The prefactorization algebra $\cF_A$ extends to a factorization algebra on $S^1$, and
\[
{\rm Hoch}_\bu(A) \simeq \check{C}(\mathfrak{U},\cF_A)
\]
for any Weiss cover of~$S^1$.
\end{thm}

This theorem tells us that factorization homology can be seen as a generalization of Hochschild homology,
where one simultaneously generalizes the circle to an arbitrary closed manifold and an algebra to a factorization algebra.

There is a variant on this result that plays a role in Section~\ref{s:compact} below.
Recall that a locally constant sheaf on $S^1$ may have monodromy: 
starting with a section on a small interval around the circle and then continuing it around the circle (by local constancy),
it may not match,
i.e., it may not extend to a global section.
Instead, transporting values around the circle produces an automorphism of the stalk at a point, known as monodromy.
Similarly, a locally constant factorization algebra $\cF$ on $S^1$ restricts to a locally constant factorization algebra on $S^1 - \{p\} \cong \RR$,
so it corresponds to an associative algebra $A$.
But in transporting elements of $A$ around the circle using this factorization algebra on $S^1$,
we may find there is an algebra automorphism $\sigma$ appearing as monodromy.
Then
\[
\int_{S^1} \cF = \cF (S^1) \simeq {\rm Hoch}_\bu(A, A_\sigma)
\]
where $A_\sigma$ denotes $A$ viewed a bimodule where $A$ acts from the left by multiplication and from the right by multiplication twisted by~$\sigma$.

\subsection{How factorization algebras appear in QFT}
\label{s:factinqft}

Factorization algebras on $\RR^n$ and $E_n$ algebras should be closely related,
as can be seen by drawing pictures of the configurations of discs that control operations.

\begin{dfn}
A factorization algebra $\cF$ on an $n$-manifold $M$ is \emph{locally constant} if for each inclusion of open discs $D \subset D'$, then the map $\cF(D) \to \cF(D')$ is a quasi-isomorphism.
\end{dfn}

For $M = \RR$, we have already discussed how locally constant factorization algebras relate to associative algebras, starting with Example~\ref{ex:associativealgebra}.
Lurie has shown the following vast extension of this example. (See section 5.4.5 of \cite{LurieHA}, particularly Theorem~5.4.5.9.)

\begin{thm}\label{thm:locisen}
There is an equivalence of $(\infty,1)$-categories between $E_n$ algebras and locally constant factorization algebras on~$\RR^n$.  
\end{thm}  

\begin{rmk}
We remark that Lurie and Ayala-Francis uses a different gluing axiom than we do. A careful comparison of the different axioms and a proof of their equivalence (for locally constant factorization algebras) can be found in~\cite{KSW}.~\hfill $\Diamond$
\end{rmk}

This perspective suggests that factorization algebras offer a natural generalization to field theories on manifolds, and not just on Euclidean spaces. Indeed \cite{CG1, CG2} provide a systematic relationship between field theories and factorization algebras.
The slogan is 
\begin{quote}
the observables of a field theory living on a manifold $M$ form a factorization algebra on~$M$. 
\end{quote}
To physicists, this claim should not be so surprising:
the Weiss cosheaf condition offers a precise version of the idea that a quantum field theory can be fully encoded by all its $n$-point functions.\footnote{This idea is not seen as accurate when one takes into account some nonperturbative aspects of QFT, notably generalized global symmetries and higher-dimensional defects. We note that it is possible to enlarge the Weiss topology to accommodate those aspects.}

To make a mathematical statement, one needs to be precise about what a field theory is.
Here, we will mean the definitions articulated in \cite{CosBook,CG1, CG2}, which encompass the Euclidean versions of many field theories studied in physics and mathematics.
Those books use the Batalin-Vilkovisky (BV) formalism for quantization, a homological method generalizing the BRST approach and widely regarded as the most powerful and general way to quantize gauge theories. 

\begin{thm}
\label{main}
The observables of a classical field theory on $M$ form a factorization algebra $\Obs^{cl}$ that assigns to every open set, a 1-shifted Poisson algebra (i.e., a commutative dg algebra equipped with a Poisson bracket of cohomological degree one). A BV quantization of this theory yields a factorization algebra $\Obs^q$ that is a flat deformation of $\Obs^{cl}$ over~$\RR[[\hbar]]$.
\end{thm}

This theorem provides an elegant interpretation of BV quantization as a kind of deformation quantization. 
In the setting of mechanics (or one-dimensional field theory, since a particle has a worldline), deformation quantization explains the transition from classical to quantum as deforming the observables from a Poisson algebra to an associative algebra. 
This theorem allows one to interpret quantization of field theories---on manifolds of arbitrary dimension---as a deformation quantization from a 1-shifted Poisson factorization algebra to a plain factorization algebra. 
In fact, in \cite{GLL,GLX} it is shown how BV quantization for one-dimensional quantum field theories connects with the ordinary theory of deformation quantization. 

This theorem is also a key connection between field theories and higher algebras.
Indeed, a the value of a locally constant factorization algebra on a $d$-disk has the structure of a $\op{Disk}_d$-algebra.
The value of a holomorphic factorization algebra on a complex $d$-dimensional polydisk has the structure of a $\op{PDisk}_d$-algebra.

\subsubsection{Relating to the Atiyah-Segal approach to TFTs}
\label{sec: reln to functorial TFT}

Atiyah \cite{AtiTFT} offered a mathematical definition for topological field theories that uses symmetric monoidal functors out of bordism categories.
His approach was inspired by Segal's definition of a conformal field theory \cite{SegCFT},
and these notions have subsequently undergone extensive development.
Notably, Baez and Dolan \cite{BaeDol} suggested and advocated for fully extended TFTs; 
their vision was revisited by Lurie \cite{LurieTFT}, who used modern machinery to provide a roadmap for using fully extended TFTs and for proving the cobordism hypothesis.

An important suggestion from \cite{LurieTFT} is that factorization homology allows one to produce nontrivial examples: 
every $E_d$ algebra provides a fully extended $d$-dimensional framed TFT,
albeit taking values in a higher analogue of the Morita bicategory (which works for the $d=1$ case).
Scheimbauer \cite{Scheim} implemented Lurie's suggestion.

\begin{thm}
For any $E_d$ algebra $A$, the functor
\[
\int_{(-)} A \colon {\rm Bord}^{or}_d \to {\rm Alg}_d
\]
determines a fully extended framed $d$-dimensional field theory with values in a higher Morita $(\infty,d)$-category of $E_d$ algebras.
\end{thm}

In combination with the theorem of Elliott and Safronov, we can connect action functionals --- the physicist's usual description of a field theory --- to Atiyah-Segal functorial field theories.

\begin{cor}
\label{cor: EllSaf}
For any Lagrangian field theory satisfying the hypotheses of Theorem~\ref{thm: EllSaf}, 
there is a fully extended framed $d$-dimensional field theory 
given by $\int_{(-)} \cA^q$ where $\cA^q$ denotes the $E_d$ algebra of quantum observables.
\end{cor}

More generally, given a map of Lie groups $G \to O(d)$ and a compatible $G$-action on the theory, 
a homotopical trivialization of this action on the quantum theory equips $\cA^q$ with the structure of a $G$-framed $E_d$ algebra and, by Lurie's work, a $G$-framed fully extended $d$-dimensional TFT.
It is interesting to ask how these kinds of functorial field theories intertwine with the better-known constructions,
such as the Reshetikhin-Turaev theories, 
where the target category is not ${\rm Alg}_d$ but some other higher-categorical extensions of~$\Vect$.

In this survey we use factorization algebras to capture the rich algebraic structure of holomorphic field theories,
but these connections with functorial TFTs raise the question of producing functorial holomorphic field theories.

%
%

\subsection{Compactification and conformal blocks}
\label{s:compact}

Factorization algebras yield interesting information when their global sections are computed,
much as sheaf cohomology is a rich source of invariants and information.
We focus here on some global aspects of holomorphic factorization algebras,
connecting with important results about vertex algebras and with the physical notion of compactification.

A map $f \colon X \to Y$ induces a pushforward functor between categories of factorization algebras.
If $\cF$ is a factorization algebra on $X$, 
then the pushforward $f_* \cF$ is the factorization algebra on $Y$ that assigns the value $\cF(f^{-1}(U))$ to an open set $U \subset Y$.

Consider a situation where $\cF$ is the factorization algebra of observables of a quantum field theory on $X$.
In the case that $f$ is a smooth map with compact fibers, 
then $f_* \cF$ is the factorization algebra of observables for the \textit{compactification} of the quantum field theory along $f$.
In the case $f \colon Z \times Y \to Y$ is the projection, 
one says that $f_* \cF$ is the compactification of $\cF$ along~$Z$.
This is akin to the famous Kaluza--Klein compactification in physics and should not be confused with the notion of ``dimensional reduction'', which only includes some of the Kaluza--Klein fields.

Another special case is the map $f: X \to \star$.
Then $f_* \cF = \int_X \cF$ is the global sections of the factorization algebra $\cF$ along $X$.
In the case that $\cF$ is a locally constant factorization algebra, 
this computation agrees with the factorization homology of $\cF$ along $X$, as defined by Ayala-Francis and Lurie.
In the case that $\cF$ is a holomorphic factorization algebra on a Riemann surface $\Sigma$, 
the global sections $\int_\Sigma \cF$ is related to a familiar invariant within the theory of vertex/chiral algebras called the ``space of conformal blocks'' \cite[\S 8]{FBZ}.
Precisely, the zeroth homology of $\int_\Sigma \cF$ is linearly dual to the space of (ordinary, non-derived) conformal blocks.

Thus, it is natural to view global sections as a higher-dimensional generalization of the notion of conformal blocks.
That is, when $X$ is any complex manifold and $\cF$ is a holomorphic factorization algebra on $X$, 
then $\int_X \cF$ plays the role of conformal blocks.
One can then ask for analogs in higher dimensions of important results involving Riemann surfaces and vertex algebras.
We sketch here an example of how the character of a (conformal) vertex algebra relates to the conformal blocks over the moduli of elliptic curves.

First, consider a \textit{real} one-dimensional \textit{topological} factorization algebra $\cA$.
If we identify $\cA$ with the associative algebra that it determines, 
then the factorization homology $\int_{S^1} \cA$ along the circle agrees with the Hochschild cohomology of the associative algebra~$\cA$.
In other words, it is the home for all of the ``traces'' defined on the algebra~$\cA$.

For a vertex algebra $\mathbb{V}$, the role of the circle is played by any elliptic curve,
and the role of Hochschild homology is played by conformal blocks on that genus one curve. 
Thus conformal blocks over elliptic curves provides the natural home for the characters of the vertex algebra.
The most basic example of the character of a conformal vertex algebra is its $q$-character,
which is defined in terms of a trace:
\beqn
\op{char}_{\mathbb{V}} (q) = \op{Tr}_{\mathbb{V}} (q^{L_0-c/24}) ,
\eeqn
where $L_0$ denotes the action of the vector field $z \del_z$.
The shift by the central charge guarantees that this expression has nice modularity properties, i.e., we can view $q$ as parametrizing the moduli of elliptic curves.

From the point of view of factorization algebras, we can understand this trace in terms of Hochschild homology of the algebra of modes associated to $\mathbb{V}$.
Indeed, suppose the conformal vertex algebra $\mathbb{V}$ arises from a holomorphic $S^1$-equivariant factorization algebra $\cF_{\mathbb{V}}$ on $\C$.
If $r \colon \C^\times \to \RR_{>0}$ is the radial projection there is a procedure for extracting the associative algebra of modes from the one-dimensional factorization algebra $r_* \cF_{\mathbb{V}}$.
Indeed, while $r_* \cF_{\mathbb{V}}$ may not be a locally constant factorization algebra, it always contains a dense subfactorization algebra which is locally constant.\footnote{To extract it, one assumes that $\cF$ is $S^1$-equivariant.
For each interval $I \subset \RR_{>0}$ one then looks at the subspace of $\cF(\pi^{-1}(I))$ consisting of the eigenspaces for this $S^1$ action.
The assignment of an interval $I$ to this subspace is a locally constant factorization algebra.}
With slight abuse of notation, we will still denote $r_* \cF_{\mathbb{F}}$ the locally constant factorization algebra obtained from this compactification.

If we further put this on the circle, we find Hochschild homology.
Thus, for each modulus~$q$ functoriality of factorization homology gives rise to a map
\beqn
\int_{S^1_{\sigma}} r_* \cF_{\mathbb{V}} = \op{Hoch}(r_* \cF_{\mathbb{V}}, \sigma) \to \int_{E_q} \cF_{\mathbb{V}} .
\eeqn
Here, $\sigma$ is an automorphism defined in terms of the given $S^1$-action on $\cF_{\mathbb{V}}$, and what appears on the left-hand side is the $\sigma$-twisted Hochschild homology.
For free theories, this map is a quasi-isomorphism.
On the other hand, given a conformal block (i.e., linear functional on $\int_{E_q} \cF_{\mathbb{V}}$) we can compose it with that map to define a trace on the modes algebra.
From this point of view, the factorization homology along an elliptic curve encodes characters in the sense of vertex algebras.
For more on this perspective see \cite{GuiLi}, where they carefully develop the notion of a trace for chiral/vertex algebras.

In higher dimensions we can build a similar picture,
replacing the moduli of elliptic curves with the moduli of {\it Hopf manifolds}.
Recall that a Hopf manifold is a quotient of the punctured affine plane
\beqn
H_{\mathbf{q}} = \left(\C^d - \{0\}\right) \slash \left((z_1,\ldots,z_d) \sim (q_1 z_1,\ldots,q_d z_d) \right)
\eeqn
where $\mathbf{q} = (q_1,\ldots,q_d)$ are complex numbers satisfying $0 < |q_i| < 1$.
For generic $\mathbf{q}$ such a complex manifold is diffeomorphic to $S^{2d-1} \times S^1$.
The factorization homology of a holomorphic factorization algebra $\cF$ along a Hopf surface of this kind can be identified with the Hochschild homology of the associative algebra of $S^{2n-1}$-modes:
\beqn
\op{Hoch}(r_* \cF , \sigma) \to \int_{H_{\mathbf{q}}} \cF 
\eeqn
where, again, $r\colon \C^d - \{0\} \to \RR_{>0}$ is the radius and $\sigma$ is an automorphism of $\cF$ induced from the $S^1$-equivariant structure of $\cF$.
Thus, we see that the factorization homology $\int_{H_{\mathbf{q}}} \cF$ naturally encodes traces of the $S^{2n-1}$-mode algebra associated to $\cF$.
In \cite{SWchar}, it is shown how an important invariant called the \textit{supersymmetric index} (or supersymmetric partition function) of a supersymmetric theory along a sphere $S^{2n-1}$ be cast into this framework at the level the holomorphic twist.

\section{Renormalization \& anomalies: constructing examples}
\label{sec: renorm}

In practice, quantum field theories are typically constructed using regularization and renormalization,
by starting with an action functional that defines a classical field theory and then finding a way to extract sensible answers from the yoga of Feynman diagrams (which naively produce divergent integrals).
One can apply similar techniques to produce holomorphic quantum field theories, and the higher algebraic structures associated to them.
A key feature of holomorphic theories is the absence of counterterms (possibly scheme-dependent), 
which can make typical quantum field theories challenging to work with.
In this section we describe the current state of the art for renormalization and anomalies of holomorphic field theories,
after reviewing some background and the similarly pleasant results for topological field theories.

\subsection{Overview and context}

The dream is that a quantum field theory is specified by the path (or functional) integral,
which is supposed to involve integrating a measure over the space of all fields.
This path integral does not (usually) exist using conventional mathematical tools,
so one mimics well-established asymptotic expansions for finite-dimensional oscillating integrals.
Feynman explained how to label the terms in these expansions using graphs;
there is an explicit algorithm that converts the action functional of the classical field theory into integrands (often distributional in nature) on the spacetime manifold $M$ and products $M^k$ of it.
Unfortunately, these formal expressions are often ill-defined as they would involve multiplying distributions (and hence the putative integral would diverge).
Regularization means a way of systematically working around these issues,
often by adjusting the integrand or finding a clever way of extracting a finite answer.
In short hand, one introduces {\it counterterms} that soak up the divergences and hence obtain useful values for the Feynman diagrams.
Renormalization means assembling the data produced by the collection of Feynman diagrams to approximate the ``true'' path integral.

To offer some notation, a graph $\gamma$ (whose exact shape depends on the theory) is expected to define a functional $W_\Gamma$ on the space of fields (i.e., a distribution).
It is described, formally, by an an integral over the product manifold $M^{|V(\gamma)|}$,  
where $V(\Gamma)$ is the vertex set of the graph.
Unfortunately, this integrand, and the distribution it purports to encode, is often not integrable.

These methods are essential --- and often quite convoluted --- for the theories that describe physical phenomena, 
such as quantum electrodynamics,
where Feynman first introduced his diagrammatics.
For topological quantum field theories, like Chern-Simons theory,
life is quite a bit simpler: Axelrod-Singer and Kontsevich developed a a very simple regularization method that removes all divergences.
For these theories, the putative divergences always arise from distributional issues along the diagonals of the product space $M^k$.
On the Fulton-MacPherson compactification of the configuration space of $k$ points in $M$ (i.e., takes the real blow-up along the diagonals),
the integrand admits a {\em smooth} extension and the integrals become well-defined.

These results guarantee the existence of factorization algebras for topological field theories like Chern-Simons theory or the Poisson $\sigma$-model or the topological B-model by the main result (Theorem~\ref{main}) of \cite{CG2}.
Thus these theories produce $E_d$ algebras by Theorem~\ref{thm: EllSaf} and Atiyah-Segal-style TFTs by Corollary~\ref{cor: EllSaf}.

This configuration space method traditionally used to study topological theories does not admit an immediate generalization to holomorphic theories.
Nevertheless, the renormalization for holomorphic field theories is equally well-behaved.

Costello, inspired by the methods of Axelrod-Singer and Kontsevich,
found a framework for perturbative renormalization that works with a large class of elliptic complexes.
When applied to Dolbeault complexes arising in Kodaira-Spencer (or BCOV) theory, 
he and Li found some remarkable simplifications and good behavior \cite{CLbcov}.
In \cite{LiVertex} Li showed, by clever analysis, that holomorphic field theories in complex dimension one (i.e., chiral conformal field theories) are finite to all orders in perturbation theory; 
this result was later extended to arbitrary Riemann surfaces in \cite{LiZhou}.
Following up, the second author laid out systematic mathematical foundations:
\begin{itemize}
\item he established key analytic results about holomorphic renormalization, 
showing that it is highly manageable, with no counterterms to 1-loop; and
\item characterized the one-loop anomalies  (i.e., obstructions to BV quantization),
\end{itemize}
see  \cite{Wthesis,Wrenorm}.
These apply to all holomorphic theories on~$\CC^n$,
but some of the techniques should extend to more general complex manifolds.

In the last few years,
there has a been a buzz of recent activity on the renormalization of holomorphic field theories.
Budzik, Gaiotto, Kulp, Wu, and Yu developed a strategy for proving that in a holomorphic quantum field theory,
certain graphs (called Laman graphs) do not produce divergences~\cite{BGKWY}.
The key idea is to combine the heat kernel analysis with the sort of compactifications that are used in the theory of configuration spaces.

Recently, Minghao Wang has extended this analysis to show that holomorphic field theories on $\C^n$ admit a regularization scheme that is divergence-free to \textit{all} orders in perturbation theory~\cite{Wang}.
In other words, both topological and holomorphic field theories are completely free of the divergences present in typical quantum field theories!
We briefly recall the approach to this regulation scheme.

\subsection{No divergences: Wang's approach}

We will allow graphs with finitely many vertices and edges, 
but also leaves (also known as hair or half-edges).\footnote{That is, a vertex may have, in addition to edges connecting it to other vertices, some incident edges that do not connect to another vertex. 
These will label inputs to an operation associated to the graph.}
We allow multiple edges between vertices, as well as loops, the precise form of which are determined by the interactions of the theory.

For a graph $\gamma$, let $E(\gamma)$ denote the number of edges,
$V(\gamma)$ denote the number of vertices,
and $L(\gamma)$ denote the number of leaves.
We fix a ``cut-off" length scale $L > 0$ that is used to remove the divergences that often come from integrating over a noncompact manifold (a so-called IR divergence), 
which is a different kind of divergence (and much easier to resolve) than those that typically plague Feynman diagrammatics.

\begin{dfn}
The {\em Schwinger space} $\op{Schw}(\gamma)$ associated to a graph $\gamma$ and a length $L$ is the non-compact manifold with corners~$(0,L]^{E(\gamma)}$.
\end{dfn}

The desired Feynman graph weight $W_\gamma$ of a graph $\gamma$ is obtained as the integral over this Scwhinger space. 

\begin{dfn}
For a perturbative field theory with space of fields $\cE$, there is a distribution valued in differential forms on Schwinger space
\[
w_{\gamma} \colon \cE^{L(\gamma)}_c \to \Omega^\bu \left( \op{Schw}(\gamma) \right)
\]
that we call the {\em regularized weight of $\gamma$}.
\end{dfn}

The possibly ill-defined weight $W_\gamma$ is obtained by integration of $w_\gamma$ over Schwinger space.
This family of distributions is constructed as an integral over $\CC^{v(\gamma) d}$,
where each edge $e$ is labeled by a smooth function $P_{t(e)}$ depending on the length of the edge and where each vertex is labeled by a distribution arising from an interaction term of the theory.
For general theories, one cannot smoothly continue the regularized weight from the positive length values to zero lengths.
The putative integrals become divergent, which is reflected by the non-compactness of Schwinger space from this perspective.

For holomorphic field theories, however, something remarkable happens.

\begin{thm}[\cite{Wang}]
Let $\cE$ be the space of fields underlying a holomorphic field theory on $\C^n$.
For every graph~$\gamma$, the differential form $w_{\gamma}$ has a natural smooth extension $\til {w_\gamma}$ to a compactification $\til{\op{Schw}(\gamma)}$ of Schwinger space.
In particular, the composite
\[
\cE^{L(\gamma)}_c \xto{\til{w_\gamma}}\Omega^\bu \left( \til{\op{Schw}(\gamma)} \right) \xto{\int} \C
\]
is a continuous linear functional on $\cE_c^{L(\gamma)}$ and shows the desired weight~$W_\gamma$ is well-defined.
\end{thm}

In other words, divergences do not appear when constructing the perturbative quantization of {\em any} holomorphic field theory on~$\CC^d$.
(Similar results hold for topological-holomorphic theories~\cite{WangWilliams}.)

This result ensures the existence of myriad examples of the higher algebraic structures we've discussed:
it provides the analysis to quantize the kinds of examples from Section~\ref{sec: ex of theories} and hence to produce polydiscs algebras and factorization algebras.

\subsection{Anomalies that obstruct quantization}

The discussion so far concerned the analytic issues that are present when trying to make the path integral precise.
For some quantum field theories, like gauge theories, there is another issue that could prevent the theory from being well-posed:
these are called {\em anomalies}.
An anomaly measures the failure of the path integral to enjoy the same symmetries as the underlying classical field theory.
In gauge theory, for instance, the perturbative renormalization might not be gauge invariant, although the classical action functional is.
Specifically, the gauge variation of a particular Feynman graph integral may be nonzero---this leads to an anomaly in the quantum theory.
In topological theories, it is known that these sort of anomalies are absent.
For holomorphic theories, however, anomalies are present.
(Within the context of Costello's formalism for perturbative QFT,
these anomalies arise as obstructions to solving a Maurer-Cartan equation known as the quantum master equation.)

Since renormalization is so well-behaved, however, such anomalies are amenable to explicit characterization.
For example, at one-loop in perturbation theory, an explicit formula for such anomalies for holomorphic theories on $\C^d$ is given in \cite{Wrenorm}.
It involves a sum over Feynman weights $W_\Gamma$ where $\Gamma$ is a polygon with $(d+1)$~vertices.

As an example, consider holomorphic BF theory on $\C^2$ with Lie algebra $\lie{g}$
and couple it to the $\beta\gamma$ system on $\C^2$ valued in a representation $(V, \rho)$ of $\lie{g}$.
The action functional is
\[
S(A,B, \beta, \gamma) = \int \langle B, \dbar A\rangle + \frac{1}{2} \langle B, [A,A]\rangle + \langle \beta, \dbar \gamma\rangle + \langle \beta, \rho(A)\gamma\rangle,
\]
where the final term encodes the ``coupling'' ({\it cf}. the discussion of these theories in Section~\ref{sec: ex of theories}).
The second and fourth terms are cubic functions of the fields,
and they provide trivalent vertices to build graphs with.
The first and third terms are quadratic, and they control what labels the edges of graphs.

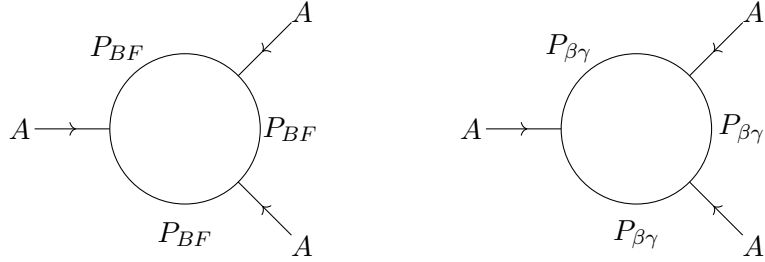
\begin{figure}
\begin{center}
\begin{tikzpicture}
	\begin{scope}


		\coordinate (O) at (0,0);
		\node[draw] (ci) at (O) [circle through=(right:1)] {};
		\draw[fermion](45:2) -- (45:1);
		\draw[fermion](180:2) -- (180:1);
		\draw[fermion](315:2) -- (315:1);
		\node at (45:2.2) {$A$};
		\node at (180:2.2) {$A$};
		\node at (315:2.2) {$A$};
		\node at (130:1.4) {$P_{BF}$};
		\node at (270:1.4) {$P_{BF}$};
		\node at (0:1.4) {$P_{BF}$};
	    	\clip (0,0) circle (1cm);
		\end{scope}

		\begin{scope}[shift={(6,0)}]
				\coordinate (O) at (0,0);
		\node[draw] (ci) at (O) [circle through=(right:1)] {};
		\draw[fermion](45:2) -- (45:1);
		\draw[fermion](180:2) -- (180:1);
		\draw[fermion](315:2) -- (315:1);
		\node at (45:2.2) {$A$};
		\node at (180:2.2) {$A$};
		\node at (315:2.2) {$A$};
		\node at (130:1.4) {$P_{\beta\gamma}$};
		\node at (270:1.4) {$P_{\beta\gamma}$};
		\node at (0:1.4) {$P_{\beta\gamma}$};
	    	\clip (0,0) circle (1cm);
 \end{scope}
\end{tikzpicture}
\caption{Diagrams whose gauge variation contribute to the anomaly in holomorphic gauge theory on~$\C^2$}
\label{fig:localization}
\end{center}
\end{figure}

The gauge variation of the two diagrams which contribute to the gauge anomaly are shown in Figure~\ref{fig:localization}.
The difference between these two diagrams is how we label the internal edges: in the first diagram they are labeled by $A-B$ fields and in the second they are labeled by $\beta-\gamma$ fields.
An anomaly can be represented by a (local) functional on the space of fields.
In each diagram at hand, this local functional depends only on the holomorphic gauge field~$A$ and has the form
\beqn
c \int \op{Tr}_V (A \wedge \del A \wedge \del A)
\eeqn
where $c$ is some nonzero constant and $V$ is some representation of $\lie{g}$.\footnote{In this example, the value of $c$ is the same for both diagrams.}
Specifically, we have:
\begin{itemize}
\item In the first diagram the representation $V$ is the adjoint representation. 
Thus, this expression vanishes when $\lie{g}$ is simple, e.g., when $\lie{g} = \lie{sl}(n)$.
\item In the second diagram the representation is that used in defining the $\beta\gamma$ system.
\end{itemize}
The sum 
\[
c \int \left(\op{Tr}_M (A \del A \del A) + \op{Tr}_V(A \del A \del A)\right)
\]
is the total anomaly for this field theory.
These anomalies of a holomorphic gauge theory on $\CC^2$ are holomorphic avatars of the Adler--Bell--Jackiw, or chiral, anomaly in ordinary gauge theory on~$\RR^4$.

In the case that $\lie{g}$ is simple (so that the first diagram does not contribute to the anomaly), 
there is an easy way to produce a non-anomalous holomorphic QFT: 
take $V$ to be a symplectic representation of the gauge Lie algebra $V = Q \oplus Q^*$, where $Q$ is any $\lie{g}$-representation.

\section{Vista: Seiberg duality and its consequences}
\label{sec: seiberg}

A major endeavor in theoretical physics is to develop a deep and effective understanding of theories like quantum chromodynamics (QCD), which governs the behavior of quarks and gluons.
One natural approach is to study theories that are close cousins but that are easier to analyze,
and a notable way is to look for theories with the same ingredients but with larger groups of symmetries.
Supersymmetric Yang--Mills theories are quite appealing in this regard,
and there have been remarkable successes in this direction.
A well-known high point is the work of Seiberg and Witten on 4-dimensional $\cN =2$ supersymmetric Yang--Mills theory,
where they explained how confinement could be realized by a version of the Mandelstam-'t Hooft mechanism (i.e., a kind of electromagnetic duality identifies confinement with a Meissner effect).
Their work offered powerful new tools that opened up much subsequent progress,
and even impacted mathematics through the Seiberg-Witten equations (which are related to the topologically twisted versions of $\cN=2$ theory).

The class of supersymmetric theories closest to the Standard Model, however, has only $\cN=1$ supersymmetry.
There too Seiberg offered inspiring and powerful insights about SQCD ($\cN=1$ supersymmetric QCD),
and one of his leading contributions is referred to as {\em Seiberg duality}.
It has had a large impact in theoretical physics,
but the impact in mathematics is much more muted,
possibly due to the lack of a topological twist for any $\cN=1$ theory.\footnote{Although citations offer only coarse insight into impact, Seiberg's paper introducing the conjecture \cite{Seiberg} has over 1800 citations listed on Inspire but less than 60 citations listed on MathSciNet.}
Such theories do admit {\em holomorphic} twists, however, 
and so we explore here what holomorphic Seiberg duality might mean and what its consequences could be in mathematics.
In our opinion this direction could be quite fruitful and we welcome others to join us in exploring it.

We begin this section by describing the holomorphic twist of $\cN=1$ supersymmetric Yang--Mills theory,
and then we offer a holomorphic version of Seiberg duality, 
initially conjectured by Richard Eager.
We end by sketching how this duality, when compactified along Riemann surfaces,
should produce equivalences between different kinds of 2-dimensional B-model theories;
these reductions are related to work of Hori and Hori--Tong and its mathematical development, notably by Segal and collaborators~\cite{Hori, HoriTong, HMSegRen, ADS}.

\subsection{Holomorphic chromodynamics and steps towards holomorphic QCD}

A 4-dimensional $\cN=1$ supersymmetric Yang--Mills theory depends on a choice of compact Lie group $G$, with Lie algebra $\fg$, 
and it lives on  the spacetime $M=\RR^4$.
The field content of the theory~is
\begin{itemize}
\item a {\it vector multiplet} that consists of a gauge field $A \in \Omega^1(M) \otimes \fg $ together with some fermions $\lambda \in C^\infty(M) \otimes \fg$, and
\item a {\it matter multiplet} (often called the chiral multiplet), depending on a choice of a representation $V$ of $\fg$, that consists of a scalar $\phi \in C^\infty(M) \otimes V$ together with some fermions $\psi \in C^\infty(M) \otimes V$. 
Additionally, there is a {\it superpotential} $W \in \cO(V)^\fg$, a $G$-invariant polynomial on~$V$. 
\end{itemize}
We always work with super vector spaces, so a fermion means an element of an odd vector space while a boson means an even element.
Thus, we view a fermion $\lambda$ as living in $\Pi (C^\infty(M) \otimes \fg)$,
where $\Pi$ denotes odd parity.

From hereon we fix a nonzero, square zero supercharge $\cQ$ of the 4-dimensional $\cN=1$ supersymmetry algebra.
This element determines a complex structure on $\RR^4$.
We also use it to twist the supersymmetric theory to produce a holomorphic field theory on $\CC^2$.
This twist also has a ``vector'' part and a ``matter'' part.
This twist of 4-dimensional $\cN=1$ supersymmetric Yang--Mills theory has been analyzed and described in \cite{CosYangian, ESW, SWchar}, within that framework for (primarily perturbative) field theory.

This twist has a moduli-theoretic description that we now sketch;
we call a point in such a moduli space a {\it solution to the equations of motion}. 
(We give later a description in terms of fields and action functionals.)
Let $G_\CC$ denote the complex Lie group associated to $G$,
and let $\fg_\CC$ denotes its complex Lie algebra.
We assume that $\fg_\CC$ is equipped with a non-degenerate symmetric invariant pairing and will freely identify $\fg_\CC$ with its linear dual using this data.

The twisted vector multiplet defines a field theory on any complex surface $X$ often called {\em holomorphic BF theory}. 
It will be convenient for us to assume that $X$ is equipped with a holomorphic volume form that we denote by $\Omega_X$. 
A solution encodes a holomorphic $G_\CC$-bundle $\cP \to X$ (by picking out a $\dbar$-connection on the underlying smooth bundle) and a holomorphic section of the adjoint bundle $\mathrm{ad}(\cP) \to X$.
The moduli space can be expressed as the mapping space ${\rm Map}(X, T^*[1]BG_{\CC})$, as a derived stack.
(Note that $T^*[1]BG_{\CC}$ is the coadjoint quotient stack $[\fg_\CC^*/G_{\CC}]$.)

\def\Sslash{\!\sslash\!}

A twisted matter multiplet also defines a field theory on any complex surface $X$. 
The superpotential determines a derived affine scheme: the derived critical locus ${\rm Crit}^d(W)$, which can be modeled as the spectrum of the dg algebra given by polyvector fields $PV(V)$ on $V$ with differential $\iota_{\d W}$.
(It is the Koszul complex associated to the Jacobi ring of $W$.)
This description makes manifest an odd symplectic structure on ${\rm Crit}^d(W)$,
as polyvector fields are an odd Poisson algebra via the Schouten bracket.
By itself (i.e., before coupling to a gauge theory) the twisted matter multiplet corresponds to the mapping space ${\rm Map}(X, {\rm Crit}^d(W))$.\footnote{Notice that AKSZ formalism equips this mapping space with a $1$-shifted symplectic structure as opposed to the desired $(-1)$-shifted symplectic structure.
Hereon we totalize all gradings to $\ZZ/2$.}

The holomorphic twist of a theory with both vector and matter multiplets can be described, globally, as a derived mapping space from $X$ to the target
\beqn
[( \fg_\CC^*\times {\rm Crit}^d(W))/ G_\CC] .
\eeqn
This is a derived quotient stack whose first factor arises from the BF theory and whose second factor arises from the matter multiplet.

\textit{Supersymmetric quantum chromodynamics} (SQCD) refers to a specific sort of four-dimensional $\cN=1$ gauge theories.
In the special, yet interesting, case that the gauge group is $G = SU(N)$ (and so $G_\CC = SL(N)$) the  matter fields consist of two types: quarks, valued in some number of copies of the fundamental (defining) representation, and antiquarks, valued in the same number of copies of the anti-fundamental (dual to the defining) representation.
The number of fundamental representations that appear is called the number of \textit{flavors}, and will be denoted $F$.
The number $N$ is called the number of \textit{colors}.
The following theorem describes the holomorphic twist of SQCD (for $G = SU(N)$). 
Henceforth, we will refer to the resulting holomorphic field theory as \textit{holomorphic QCD}.


\begin{thm}[\cite{ESW,SWchar}]
The holomorphic twist of $G = SU(N)$ SQCD with $F$ flavors and super potential $W$ on $\C^2 = \RR^4$ is equivalent to the following holomorphic field theory:
\begin{itemize}
\item the gauge fields are those of holomorphic BF theory valued in $\sl(N)$, whose fields are
\[
A \in \Omega^{0,\bu}(\CC^2 , \sl(N))[1] \quad\text{and}\quad B \in \Omega^{2,\bu}(\CC^2, \sl(N)^*);
\]
\item the ``quarks'' are fields of a holomorphic $\beta\gamma$ system valued in 
\[
V_Q = \Hom(\CC^F, \underline{N}) \cong \underline{N}^{\oplus F} \cong \CC^{NF},
\] 
whose fields are
\[
Q \in \Omega^{0,\bu}(\CC^2 , V_Q) \quad\text{and}\quad R \in \Omega^{2,\bu}(\CC^2 , V_Q^*)[1];
\]
\item the ``antiquarks'' are fields of a holomorphic $\beta\gamma$ system valued in the dual representation 
\[
V^*_Q \cong \Hom(\underline{N}, \CC^F) \cong (\underline{N}^*)^{\oplus F} \cong \CC^{NF}, 
\]
whose fields are
\[
Q^\dag \in \Omega^{0,\bu}(\CC^2 , V^*_Q) \quad\text{and}\quad R^\dag \in \Omega^{2,\bu}(\CC^2 , V_Q) [1].
\]
\end{itemize}
(The superscript $\dag$ does not mean we apply an operation to $Q$ or $R$; 
it is simply to remind us that these fields are the ``antiparticles''.)

The action functional is
\[
S_{QCD} = \int B\, F^{0,2}_A + \int R \,\dbar_A Q + \int R^\dag \,\dbar_A Q^\dag + \int\d^2 z \, W(Q) ,
\]
where $\dbar_A = \dbar + A$ is the covariant $\dbar$-operator.\footnote{Although this perturbative theory only depends on the Lie algebra of $G$, the group structure will be relevant when we discuss local operators and moduli of vacua.}
\end{thm}

\begin{rmk}
Notice that to write down the superpotential term, we are utilizing the holomorphic volume form~$\d^2 z$. 
Relatedly, unless $W = 0$, this theory is only $\ZZ/2$-graded, so all fields are considered merely as even or odd.
For instance the component of the gauge field $A^{0,1}$ is even while the component of the matter field $Q^{0,1}$ is odd. 
When $W$ is homogenous one can lift the theory to a $\ZZ$-grading by making additional choices.
\end{rmk}

To formulate Seiberg duality, we will need to allow auxiliary matter fields that we call ``mesons.''
As above, we fix $G_{\CC} = SL(N,\CC)$.

\begin{dfn}
A holomorphic QCD theory {\em with mesons}
is a holomorphic QCD with $N$ colors and $F$ flavors along with {\em meson fields}, which form a $\beta\gamma$ system valued in the vector space $V_M$, whose fields we denote by
\[
\mu \in \Omega^{0,\bu}(\CC^2 , V_M) \quad\text{and}\quad \nu \in \Omega^{2,\bu}(\CC^2 , V_M^*) [1] .
\]
The action functional is
\[
S_{QCDwM} = S_{QCD} + \int \nu \,\dbar \mu + \int \d^2 z \; \til W(Q, Q^\dag, \mu) .
\] 
\end{dfn}

Note that the fields $\mu$ are ``uncharged,'' i.e., $SL(N)$ acts trivially on $V_M$.
Hence we use $\dbar$ rather than $\dbar_A$.
The only coupling that the meson $\mu$ has to the other fields is through the new superpotential~$\til W$. 

\subsection{A first pass at Seiberg duality}

Our long-term goal is to show that two special cases of holomorphic QCD are {\em dual}.
The notion of duality can refer to several kinds of relationships,
but loosely speaking it means that some aspects of the two theories are equivalent under some correspondence.
We will postpone saying anything precise till we have some terminology in place.


The proposed duality involves two holomorphic theories on $\CC^2$. 
We call one side {\em electric} and the other {\em magnetic}.
These turn out to be the twists of the original electric and magnetic supersymmetric theories of the famous Seiberg duality, 
which motivates our terminology.

\begin{dfn}
Fix a positive integer $F$, for flavor, and fix integers $N$ and $\check{N}$ such that $F = N + \check{N}$ and $F \geq N, \check{N} \geq 2$.
\begin{itemize}
\item
The {\em electric theory} is the holomorphic QCD with $N$ colors and $F$ flavors and {\em without} mesons.
We denote it by~$\cT_{\sfE}(F, N)$.
We use $A, B, Q, R, Q^\dag$, and $R^\dag$ for the fields.
For simplicity, we don't consider any superpotential.
\item
The {\em magnetic theory} $\cT_{\sfM}(F, \check{N})$ is the holomorphic QCD with $\check{N}$ colors and $F$ flavors, but with mesons $\check{\mu}$ that take values in $$V_{\check{M}} = \End(\CC^F) \cong \CC^{\oplus F^2}$$
(and $\check{\nu}$ for the partner field to $\check{\mu}$).
We use $\check{A}, \check{B}, \check{Q}, \check{R}, \check{Q}^\dag$, and $\check{R}^\dag$ for the other fields.
The superpotential is 
\beqn\label{eqn:checkW}
\til W =  \Tr(\check{Q}^\dag \check{\mu}  \check{Q})
\eeqn
where we compose the linear maps before tracing over~$V_{\check{N}}$.
\end{itemize}
\end{dfn}

Note that the lower bound on $N$ and $\check{N}$ is simply to ensure that we have nontrivial gauge groups.


Seiberg duality states how certain ``electric'' and ``magnetic'' versions of SQCD should be equivalent.
The electric and magnetic holomorphic QCDs are precisely the holomorphic twists of those that appear in Seiberg duality,
and so Richard Eager suggested that Seiberg duality should survive after taking the holomorphic twist.

\begin{conj*}[\cite{Eager}]
Let $F = N + \check{N}$ and $F \geq N, \check{N} \geq 2$.
Construct a quantization of the theories $\cT_{\sfE}(F,N)$ and $\cT_{\sfM}(F,\check{N})$.
There is an equivalence of 2-dimensional holomorphic factorization algebras
\[
\Obs_{\sfE} \simeq \Obs_{\sfM}  
\]
where these are the observables of the electric and magnetic theories, respectively.
\end{conj*}

\begin{rmk}
The statement is an equivalence of factorization algebras valued in $\ZZ/2$-graded cochain complexes since the field theories are only $\ZZ/2$-graded.
\end{rmk}

This claim should seem rather surprising.
In many ways these theories look, on their face, like very different theories:
the gauge groups $SL(N)$ and $SL(\check{N})$ are (usually) different, and
one side has mesons and the other does not.
Of course, relating quite different theories also makes the conjecture powerful.

One appealing aspect of Eager's conjecture is that it is mathematically precise:
duality would mean giving an equivalence of factorization algebras.
It is analogous to claiming that two vertex algebras are equivalent,
such as in the well-known boson-fermion correspondence or in the correspondence between the chiral rings of the A- and B-models of mirror symmetry.
In this sense it is a target at which a mathematician can aim,
rather than the more amorphous statements of duality that one often finds in the physics literature.
In that literature, however, there is guidance on what map of factorization algebras provides the equivalence,
which we describe in the next section.

In general, Eager's suggestion opens up a fascinating vein of inquiry:
\begin{quote}
how do physical results about supersymmetric gauge theories behave under twisting to holomorphic or topological theories?
\end{quote}
We have found, so far, that these results often admit clean mathematical formulations and that the holomorphic (or topological) twists involve constructions of natural interest in geometry and representation theory.
In \cite{EGW} we and Chris Elliott explored twists of $\cN=4$ supersymmetric Yang--Mills theory,
provided their quantizations, and began analyzing their factorization algebras.
One might can also ask about the twisted versions of S-duality and how to provide those conjectured equivalences.
The answer to this question has been studied by Raghavendran and Yoo \cite{RYsduality} using Costello and Li's description of spacetime twists of superstring theories in terms of topological string theory \cite{CLsugra,CLtypeI}. 
A further localization of $S$-duality has also been studied by Costello and Gaiotto~\cite{CGholography}.

\begin{rmk}
Our discussion here is ahistorical.
Witten initiated this vein of research for twists that are topological \cite{WittenTwist},
and Costello advocated holomorphic twists \cite{CosSUSY}.  
Costello and Li initiated the use of holomorphic twists to study dualities particularly in the setting of the topological B-model (see \cite{CLbcov} and follow ups), and further extended these ideas to non-topological strings \cite{CLsugra}. 
These ideas are extensively used in the program of twisted holography, which is growing rapidly at the moment.
(The interested reader might track the papers that cite those we have mentioned.)
Eager's suggestion was a natural variation on this theme but it (seems to) involve physical ideas and methods that have not (yet) played a large role in Costello's program.
\end{rmk}

\subsection{A remark relating to homological projective duality}

A standard maneuver used in physics is compactification:
given a field theory $\cT$ in $n$ dimensions, one takes compact ($n$-$k$)-manifold $N$ and produces a field theory $\cT_N$ in $k$ dimensions by studying $\cT$ on manifolds of the form $N \times \Sigma$, where $\Sigma$ is $k$-dimensional.
One then considers the push forward of the theory along the projection map $N \times \Sigma \to \Sigma$, to obtain a theory along $\Sigma$.
(A close cousin is dimensional reduction, where one restricts to solutions of the classical field theory that are constant along~$N$.)

We now apply this maneuver to holomorphic Seiberg duality:
by compactifying along a Riemann surface, we should generate a conjectural equivalence between two-dimensional field theories.
As we will see, this new conjecture resembles the Hori-Tong dualities that have led to Hori-mological mirror symmetry~\cite{Hori, HoriTong, HMSegRen}.

Consider a holomorphic theory defined on $\C^\times \times \C^\times$ 
and push forward along the map $\rho: \CC^\times \times \CC^\times \to \RR_{>0} \times \RR_{>0}$ 
that sends $(z_1,z_2)$ to $(|z_1|, |z_2|)$.
In other words, take the torus compactification of the holomorphic theory.
For a holomorphic factorization algebra $\cF$, 
the pushforward $\rho_* \cF$ contains a locally constant factorization algebra on $\RR_{> 0}^2 \cong \RR^2$ (in essence, by taking the Fourier modes on the torus fibers).
This locally constant factorization algebra corresponds to an $E_2$-algebra.\footnote{Recall Section~\ref{sec: red vs chiral}, where we compactified on a circle and then broke things up using Fourier modes. This construction extracted an associative ($E_1)$ algebra that is a Weyl algebra of a formal loop space from a holomorphic field theory on~$C^\times$.}

Let's start with holomorphic QCD and see what Seiberg duality would imply.
On the electric side, 
we have a two-dimensional holomorphic BF theory coupled to a complex two-dimensional $\beta\gamma$ system.
By itself, the $\beta\gamma$ system with target $\C^{NF}$ yields, by torus compactification, 
the topological B-model theory whose target is the formal double loop space~$L^2\C^{NF}$. 
The algebra of functions of this target space is $\cO(\C^{NF}) ((z))((w))$, 
where $z,w$ are formal loop parameters that we imagine as describing the winding of the two circles we have compactified. 
The dimensional reduction then picks out the constant modes along each circle, 
so we find a topological B-model with target $\C^{NF}$.
Similarly, holomorphic BF theory dimensionally-reduces to topological BF theory,
which can be viewed as the B-model with target the classifying space~$BGL(N)$.
In sum, the dimensional reduction of electric holomorphic QCD is the topological B-model with target the stack $\C^{NF} / GL(N)$.
As a derived mapping space, the moduli of solutions is given by
\beqn
{\rm Map}\left(\RR^2, T^*[1] \C^{NF} / GL(N)\right) .
\eeqn
Dimensional reduction of the factorization algebra of observables produces the $E_2$ algebra of (holomorphic) polyvector fields on the stack~$\C^{NF} / GL(N)$. 

On the magnetic side the reduction is similar, except this time we must take into account the superpotential. 
What results is a Landau--Ginzburg B-model whose target is the stack
\beqn
\left(\C^{(F-N)F} / GL(F-N) \right) \times \C^{F^2} .
\eeqn
The first factor comes from the reduction of the dual quarks coupled to the dual gauge theory, and the second factor comes from the reduction of the mesons. 
The potential for the Landau--Ginzburg model is precisely the potential $\check{W}$ seen in~\eqref{eqn:checkW}.
This procedure produces an $E_2$ algebra by deforming the differential on the $E_2$ algebra of polyvector fields on the stack $\C^{(F-N)F} / GL(F-N) \times \C^{F^2}$, 
namely use the differential $\{\check{W},-\}$ where $\{-,-\}$ is the Schouten bracket.

Holomorphic Seiberg duality thus predicts an equivalence of two-dimensional topological field theories,
or at least of $E_2$-algebras.

\begin{conj*}
There is an equivalence of $E_2$ algebras between
\begin{itemize} 
\item the {\em electric} $E_2$ algebra of polyvector fields on the stack $\C^{NF} / GL(N)$, and 
\item the {\em magnetic} $E_2$ algebra of polyvector fields on the stack $(\C^{(F-N)F} / GL(F-N)) \times \C^{F^2}$ with differential~$\{\check{W},-\}$. 
\end{itemize}
\end{conj*}


This conjectured equivalence resembles Hori-Tong dualities.
Note that it is often said that Hori-Tong dualities are like Seiberg dualities!
We hope that insights from the study of these dualities --- such as the deep work on homological projective duality \cite{KuzHPD, ThomasHPD, Perry} --- might apply here, and also feed back towards holomorphic Seiberg duality.

\section{Appendix: A zoo of holomorphic field theories}
\label{sec: ex of theories}

In this section we describe several examples of holomorphic field theories,
in the style we used for holomorphic Chern-Simons theory.
Our aim here is to give a sense of how take physical ideas and find holomorphic analogs.

In traditional physics there are a few basic types of field theories, 
and they can be assembled into combinations,
of which the Standard Model offers a nice example.
The basic types are\footnote{For the experts, we remark that we focus above on the field content that motivates or characterizes, and mostly ignore (just for the moment!) issues like gauge transformations and the associated ghost fields.}
\begin{itemize}
\item {\em Scalar} field theories where the fields are maps from $M$ to a vector space $V$, or natural generalizations like taking sections of a vector bundle $V \to M$. The Higgs boson appears as part of a scalar field theory.
\item {\em Fermionic} field theories where the fields are maps from $M$ to an {\em odd} vector space (in the sense of supermathematics), typically a section of a spinor bundle. Such theories are often used to describe matter, such as electron fields.
\item $\sigma$-{\em models} where the fields are maps from $M$ into a manifold~$X$. (One can view these as nonlinear generalizations of scalar field theories.)
\item {\em Gauge} field theories where the fields are connections on a principal $G$-bundle over $M$, with $G$ a Lie group. The Yang--Mills theories for the electroweak and strong forces are examples.
\item {\em Gravity} theories where the fields are metrics on $M$, or something that similarly controls the geometry of $M$. Einstein's theory is a paradigmatic example.
\end{itemize}
The Standard Model involves the first three types,
and the different fields interact in intricate ways.\footnote{For instance, the scalar and fermionic fields are sections of associated bundles to the principal bundle appearing in the gauge theory.}

There are holomorphic versions of all these basic types, 
and we now spell out some simple examples.\footnote{See \cite{DonThom} for a pioneering and illuminating discussion of how to find such analogs, by mathematicians. See Nekrasov's thesis \cite{NekThesis} for a master physicist's view on holomorphic field theories, before he  reconfigured our understanding of much of gauge theory.}
To give the readers --- possibly a mixed audience of mathematicians and physicists --- some feel for the examples, 
we describe them in several ways: 
in terms of PDE and action functionals but also as moduli spaces.
We do not offer a thorough justification of the dictionary, for which see \cite{CosSUSY} for a wonderful and motivating exposition,
and see \cite{AlfonsiYoung} and \cite{Steffens23} for recent work on this topic.

For an example that mixes several types, see Section~\ref{sec: seiberg}
where we discuss a holomorphic analog of quantum chromodynamics (QCD), the theory that governs the strong force and quarks.

A motivation for us to explore holomorphic field theories is that physical insights and conjectures, such as Seiberg duality, translate into surprising predictions in the holomorphic setting,
which is typically more accessible for mathematics.
(One might also hope that progress in the holomorphic setting then feeds back useful insights to physics.)

\subsubsection{Holomorphic scalar field theories}

We met the standard example of a scalar field theory in Example~\ref{eg: massless scalar},
where the field $\phi$ a smooth real-valued function on a manifold $M$
and the equation of motion is $\Delta \phi = 0$,
so that the moduli space of solutions is the space of harmonic functions on~$M$.

A natural analogue is to have the field be a complex-valued function $\gamma$ on a complex manifold $M$ with an equation of motion $\dbar \gamma = 0$,
so that the moduli space of solutions is the space of holomorphic functions on~$M$.

Notice that we are shifting from a second-order differential operator $\Delta$ to first-order differential operators $\dbar$.
Reworking a theory to express its equations of motion using first-order operators is sometimes called the first-order formalism.

We will start with $M$ being a Riemann surface for simplicity.
Consider the action functional
\beqn
\label{eqn: 1d betagamma}
\int \beta \wedge \dbar \gamma
\eeqn
where $\beta$ is a $(1,0)$-form on $M$.
Note that this action functional also imposes the equation of motion $\dbar \beta = 0$,
so that $\beta$ should be a holomorphic $1$-form.

It should seem reasonable to view two solutions $\beta_1$ and $\beta_2$ as equivalent if they differ by a $\dbar$-exact term.
If one imposes this equivalence relation (a kind of ``gauge symmetry''),
then one is on the path towards setting up the theory in the Batalin-Vilkovisky formalism.

Derived geometry is the mathematical version of BV formalism.
From this point of view, the moduli space of solutions for this $\beta\gamma$-system is the shifted cotangent bundle
\[
T^*[-1] \Map_\dbar(M,\CC) = \RR\Gamma(M,\cO) \times \RR\Gamma(M,\Omega^1_{hol}),
\]
where $\Omega^1_{hol}$ is the sheaf of holomorphic one-forms.

\begin{rmk}
It is possible to directly relate the usual massless scalar field theory to this $\beta\gamma$-system as follows.

When $M$ is two-dimensional, oriented and Riemannian,
we can view $M$ as equipped with a complex structure.
Every complex-valued harmonic function then decomposes as a sum of a holomorphic and anti-holomorphic function. 
That is, if we complexify the scalar field $\phi$, 
we can view it as a pair of fields $\gamma$ and $\Bar{\gamma}$,
where the equation of motion of $\gamma$ is $\dbar \gamma = 0$ and the equation of motion is $\partial \Bar{\gamma} = 0$.
These components are sometimes called the chiral and anti-chiral ``sectors.'' 

The $\beta\gamma$-system also plays a role in superstring theory as the free field theory underlying bosonic ghosts corresponding to worldsheet supersymmetry~\cite[Chapter 3]{Polchinski}.
There is likewise a complex-conjugate $\Bar{\beta}\Bar{\gamma}$-system encoding the antichiral sector of the free boson.
To recover the full boson, one must couple the chiral and anti-chiral sectors.
See \cite{Kapustin:2005pt,GGW} for more details.~\hfill$\Diamond$
\end{rmk}

This example admits a natural generalization to an arbitrary complex manifold.
The input data is a complex manifold $X$ of arbitrary complex dimension $n$, 
together with a holomorphic vector bundle $V$ on~$X$.
(Above we took this bundle to be trivial.)
Each field of the theory consists of a pair $(\gamma,\beta)$ where 
\[
\gamma \in \Gamma(X,V) \quad\text{and}\quad \beta \in \Omega^{n,n-1}(X,V^*) .
\]
The action functional is 
\[
\int_X \beta \wedge \dbar \gamma
\]
and the equations of motion are
\[
\dbar \gamma = 0 = \dbar \beta.
\]
Solutions consist of a holomorphic function $\gamma$ and a holomorphic $(n,n-1)$-form.
We again suggest to identify solutions to the second equation $\dbar \beta = 0$ that differ by a $\dbar$-exact form.

\subsubsection{Holomorphic fermions}

It is straightforward to modify the example above to allow holomorphic {\it super} vector bundles,
i.e., holomorphic $\ZZ/2$-graded vector bundles.

For example, given a holomorphic vector bundle $V$ on a complex $n$-fold~$X$,
let $\Pi V$ denote the parity-reversed vector bundle,
i.e., we now give $V$ an odd grading.
There is a $bc$ system that is the fermionic analog of the $\beta\gamma$ system.
The fields of the theory consist of a pair $(c,b)$ where 
\[
c \in \Gamma(X,\Pi V) \quad\text{and}\quad b \in \Omega^{n,n-1}(X,\Pi V^*) .
\]
Thus, $b$ and $c$ are fermionic fields.
The action functional is 
\[
\int_X b \wedge \dbar c,
\]
and the equations of motion are
\[
\dbar c = 0 = \dbar c.
\]
We again suggest to identify solutions to the equation $\dbar c = 0$ that differ by a $\dbar$-exact form.

There is another class of fermionic theories that is quite interesting as they are more fundamental than the previous example.
Let $X$ be a complex $n$-fold and $L \to X$ a line bundle on $X$ such that
\begin{itemize}
\item $n \equiv 1 \mod 4$ (i.e., $n=4k+1$ for some $k$) and
\item $L$ is a square root of the canonical bundle $L^{\otimes 2} \cong K_X$.
\end{itemize}
A field is a section
\[
\psi \in \Omega^{0,2k}(X,\Pi L),
\]
and the action functional is
\[ 
\int_X \psi \wedge \dbar \psi .
\]
The equation of motion for this system is simply $\dbar \psi = 0$.
To study the quantization of this system, it is necessary to identify solutions that differ by $\dbar$-exact terms (and identifying identifications that differ by $\dbar$-exact terms, etc.).

From the point of view of derived geometry, the cochain complex
\[
\Omega^{0,\bullet}(X, \Pi L)[2k]
\]
presents a derived super stack that should be viewed as the moduli of solutions for this theory.

\subsubsection{Holomorphic $\sigma$-models}

A $\sigma$-model is a field theory whose space of fields contains a mapping space $\Map(M,N)$.
When $N = \RR$, this mapping space simply becomes functions and so the $\sigma$-model reduces to a scalar field theory.

There is a natural generalization of the $\beta\gamma$ system to a $\sigma$-model.
Let $M$ be a complex $d$-manifold and let $N$ be a complex $n$-manifold.
The fields are a pair $(\gamma, \beta)$ where
\[
\gamma \in \Map(M,N) \quad\text{and}\quad \beta \in \Omega^{n,n-1}(M,\gamma^* T^*_N) .
\]
Note that $\dbar \gamma$ lives in $\Omega^{0,1}(M, \gamma^* T_N)$, 
a section of the pullback of the tangent bundle of the target.
The action functional is 
\[
\int_X \beta \wedge \dbar \gamma,
\]
where we use the canonical evaluation pairing between the pullback of the tangent and cotangent bundles of the target.
The equations of motion are
\[
\dbar \gamma = 0 = \dbar \beta,
\]
so that solutions consist of a holomorphic map $\gamma: M \to N$ and a holomorphic section $\beta$ of the pullback tangent bundle.
As usual, we suggest to identify solutions to the equation $\dbar \beta = 0$ that differ by a $\dbar$-exact form.

From derived geometry point of view, the moduli space of solutions for this $\beta\gamma$-system is the shifted cotangent bundle
\[
T^*[-1] \Map_\dbar(M,N) ,
\]
where $\Map_\dbar(M,N)$ is the derived stack of holomorphic maps $M \to N$.
This is a natural space to study in complex geometry.

\subsubsection{Holomorphic gauge theories}
\label{sec: hol BF dfn}

We have already met holomorphic Chern--Simons theory, which is defined on Calabi--Yau three-folds.
There is a another class of holomorphic gauge theories 
defined on any complex manifold, with no restriction on dimension 
and without the need for a Calabi--Yau structure. 

Let $X$ be a complex manifold of complex dimension~$n$.
Fix a complex Lie group $G$ and an invariant symmetric bilinear form $\langle-,-\rangle$ on its Lie algebra~$\fg$.
Fix a holomorphic principal $G$-bundle $P \to X$.
The fields consist of pairs $(A,B)$ where 
\[
A \in \Omega^{0,1}(X, \fg_P)
\] 
where $\fg_P$ is the adjoint bundle of $P$ and
\beqn
B \in \Omega^{n,n-2}(X, \fg_P^*) ,
\eeqn
where $\fg_P^*$ is the coadjoint bundle of $P$. 
The action functional is 
\beqn
S(A,B) = \int_X \langle B \wedge F_A \rangle
\eeqn
where $F_A$ is the curvature of the connection $A$.
The equations of motion are
\begin{align*}
F_A^{0,2} = \dbar A + \frac12 [A,A] & = 0 \\
 \dbar B + [A,B] & = 0.
\end{align*}
A solution to the first equation provides a new $\dbar$-connection 
\[
\dbar_A = \dbar_P + [A,-] 
\]
on $P$.
A solution to the second equation is a $\dbar_A$-flat $(n,n-2)$-form-valued section of the coadjoint bundle.
This theory is often called {\em holomorphic BF theory}, due to the integrand of the action functional.

\begin{rmk}
Notice that this action functional resembles the $\beta\gamma$ system, 
with $A$ playing the role of $\gamma$ and $B$ playing the role of $\beta$.
One also recognizes this action functional as a holomorphic analog of topological $BF$ theory. 
\end{rmk}

It is natural to view two connections as equivalent if they are related by a gauge transformation
and to identify $B$-solutions that are related by a $\dbar_A$-exact term.
From the point of view of derived geometry,
we thus see holomorphic BF theory as describing the moduli space
\[
T^*[-1] \Bun_G(X),
\]
if we run over all principal bundles.

\begin{rmk}
Holomorphic BF theory exhibits some interesting properties at the quantum level,
at least perturbatively.
For example, in \cite[\S3.2]{BGKWWY} it is shown that on $X = \C^2$ 
the perturbative quantization of holomorphic BF theory is almost a topological field theory. 
More precisely, to first-order in perturbation theory, 
the divergence-free component of the stress tensor is rendered cohomologically trivial at one-loop in perturbation theory.
In other words, the quantum theory appears to be totally ``invariant'' under divergence-free vector fields,
although not all vector fields (as in the case of a TFT).
\end{rmk}

It is possible to formulate holomorphic versions of higher abelian gauge theories,
associated to moduli spaces of holomorphic $k$-gerbes with connection, see \cite{SWtensor,BWfive}.

\subsubsection{Interlude on holomorphic field theories as chiral sectors}

Every holomorphic theory we have introduced so far has a ``kinetic term'' (the quadratic term in the action functional) that only involves $\dbar$.
Thus the linearized equations of motion always have the form $\dbar \phi = 0$,
so that they manifestly relate to holomorphic functions (or sections).
There are, however, holomorphic theories with a different flavor.

We start with complex dimension~$1$.
The usual scalar field theory can be written as 
\[
\int_{\Sigma} \phi \wedge \del \dbar \phi, 
\]
where $\Sigma$ is a Riemann surface.
Using integration by parts, one can use instead the action functional
\[
\int_{\Sigma} \partial \phi \wedge \dbar \phi.
\]
This theory is not holomorphic, but we can force it to be holomorphic by demanding that $\dbar \phi = 0$.
Importantly, this constraint is \textit{not} an equation of motion.
When we impose this constraint, we can rewrite the action functional in terms of the $(1,0)$-form $\alpha \in \Omega^{1,0}(\Sigma)$ as
\beqn\label{eqn:scalar1}
\int_\Sigma \alpha \wedge \dbar (\del^{-1} \alpha) ,
\eeqn
where $\alpha$ plays the role of $\partial \phi$.
In this theory the kinetic term is the whole action,
and it has a rather different flavor from the examples we have seen so far.\footnote{This action functional is a precursor to the ``holomorphic gravity'' theory we meet next, 
known as Kodaira--Spencer theory.}
One might call this the {\em chiral sector} of the free scalar theory,
as it consists of the complex-valued scalar fields satisfying a ``chiral'' constraint.

There is a higher-dimensional analog of the chiral sector of the scalar field.
Let $X$ be a complex manifold of dimension $2n+1$, i.e., odd.
The fields are the $(n+1,n)$-forms $\alpha \in \Omega^{n+1,n}(X)$ that satisfy the constraint
\beqn
\del \alpha = 0.
\eeqn
(When $n=0$ and $X$ is simply a Riemann surface, 
this constraint is automatic.)
The action functional is $\int_X \alpha \dbar(\del^{-1} \alpha)$, the same as~\eqref{eqn:scalar1}.
Generally, this theory is the holomorphic part of a higher-form scalar theory, with action
\beqn
\int_X \|C\|^2 \op{dvol} = \int_X \d C \wedge \star \d C , 
\eeqn
where $C \in \Omega^{2n+1}(X)$.
For more details we refer to~\cite{GRW_WZW}.

The reader may (and should) be concerned with the apparent ``non-local'' form of the action functional in \eqref{eqn:scalar1}.
This can be made rigorous by cohomologically resolving the condition that $\alpha$ be $\del$-closed, rather than strictly imposing the constraint~\cite{CLbcov}.

\subsubsection{Holomorphic ``gravity''}

Any theory of gravity involves, in part, metric structure on a given manifold,
so that solutions to the equations of motion produce metrics satisfying a PDE. 
For a perturbative gravity theory,
one studies deformations of some metric structure on a manifold.
There is a holomorphic avatar of this situation: deformations of complex structure. 
The first appearance of such a theory was motivated by topological string theory and involves complex structures on Calabi--Yau threefolds~\cite{BCOV}.
It is known \textit{Kodaira--Spencer gravity} or the Kodaira-Spencer field theory,
due to its relationship with deformation of complex structures.\footnote{Kodaira and Spencer's work was an important initial clue in the relationship between deformation theory and dg Lie algebras.}
It has been pursued mathematically by Costello and Li as an approach to the quantum higher genus $B$-model, starting with \cite{CLbcov} and continuing in a series of papers.
In these works, Costello and Li also provide a generalization of the theory that makes sense in any complex dimension; 
here we focus just on the most familiar three-dimensional version.

Let $X$ be a complex threefold equipped with a Calabi--Yau structure; we denote by $\Omega$ the holomorphic volume form.
The primary field of Kodaira-Spencer theory is 
\beqn
\mu \in \Omega^{0,1}(X,\T_X),
\eeqn 
a $(0,1)$-form  on $X$ with values in the holomorphic tangent bundle~$\T_X$.
A constraint is imposed on these sections:
we assume that $\mu$ preserve the holomorphic volume form $\Omega$, meaning 
\[
{\rm div}_\Omega(\mu) = 0,
\]
i.e., its divergence vanishes.
In short, the fields are divergence-free vector fields.
The classical equation of motion is 
\beqn
\dbar \mu + \frac12 [\mu,\mu]_{SN} = 0 ,
\eeqn
where $[-,-]_{SN}$ denotes here the Schouten--Nijenhuis bracket.\footnote{On vector fields, $[-,-]_{SN}$ is simply the Lie bracket, and it is the natural multilinear extension to polyvector fields.}
The reader might recognize this equation as the Maurer--Cartan equation controlling deformations of complex structure $\dbar\rightsquigarrow \dbar + \mu$.
Formulating the action functional to obtain this equation of motion is a little involved, 
and we refer to \cite{CLbcov} for a systematic discussion.

\begin{rmk}
The kinetic term of this theory has the form
\[
\int_X \Omega \wedge \left[\Omega \vee \left(\frac12 \mu\, \dbar \del^{-1} \mu\right) \right]
\]
and hence resembles the chiral sector of the scalar, which we just discussed in the preceding example.
\end{rmk}

\printbibliography

\end{document}